\documentclass[10pt]{amsart}

\usepackage{graphicx}        % standard LaTeX graphics tool
\usepackage{multicol}        % used for the two-column index
\usepackage[bottom]{footmisc}% places footnotes at page bottom
\usepackage{amscd,amsmath,amssymb,amsfonts}
\usepackage[cmtip, all]{xy}

\usepackage{bbm}
\usepackage{tikz-cd}
\usepackage{enumerate}
\usepackage{mathrsfs}
\usepackage{amsthm}

\usepackage[plainpages, urlcolor=blue]{hyperref}

\newtheorem{theorem}{Theorem}[section]
\newtheorem{proposition}[theorem]{Proposition}
\newtheorem{lemma}[theorem]{Lemma}

\newtheorem{def/prop}[theorem]{Definition/Proposition}

\newtheorem{corollary}[theorem]{Corollary}

\newtheorem{IntroThm}{Theorem}

\theoremstyle{definition}
\newtheorem{definition}[theorem]{Definition}
\theoremstyle{remark}
\newtheorem{remark}[theorem]{Remark}

\numberwithin{equation}{section}

\newcommand{\CH}{{\rm CH}}

\newcommand{\Hom}{{\rm Hom}}

\newcommand{\Spec}{{\rm Spec\,}}

\newcommand{\0}{\emptyset}
\newcommand{\sE}{{\mathcal E}}
\newcommand{\sF}{{\mathcal F}}

\newcommand{\sH}{{\mathcal H}}
\newcommand{\sI}{{\mathcal I}}

\newcommand{\sK}{{\mathcal K}}
\newcommand{\sL}{{\mathcal L}}
\newcommand{\sM}{{\mathcal M}}
\newcommand{\sN}{{\mathcal N}}
\newcommand{\sO}{{\mathcal O}}

\newcommand{\sV}{{\mathcal V}}
\newcommand{\sW}{{\mathcal W}}

\newcommand{\A}{{\mathbb A}}

\newcommand{\C}{{\mathbb C}}

\newcommand{\G}{{\mathbb G}}

\renewcommand{\P}{{\mathbb P}}

\newcommand{\R}{{\mathbb R}}

\newcommand{\V}{{\mathbb V}}

\newcommand{\Z}{{\mathbb Z}}

\newcommand{\mov}{{\mathfrak{m}}}
\newcommand{\fix}{{\mathfrak{f}}}

\newcommand{\BM}{{\operatorname{B.M.}}}

\renewcommand{\det}{\operatorname{det}}

\newcommand{\id}{{\operatorname{\rm Id}}}

\newcommand{\Sch}{{\operatorname{\mathbf{Sch}}}}

\newcommand{\Def}{{\operatorname{Def}}} 
\renewcommand{\sp}{{\operatorname{sp}}}

\newcommand{\Ind}{{\operatorname{ind}}}
\newcommand{\fr}{{\operatorname{fr}}} 
\newcommand{\gen}{{\operatorname{gen}}}

\newcommand{\del}{\partial}

\renewcommand{\max}{{\operatorname{\rm max}}}

\newcommand{\Spc}{{\mathbf{Spc}}}

\newcommand{\Sm}{{\mathbf{Sm}}}

\newcommand{\Sym}{{\operatorname{Sym}}}

\newcommand{\rnk}{{\operatorname{\text{rnk}}}} 
\newcommand{\Bl}{\text{Bl}}

\newcommand{\SH}{{\operatorname{SH}}} 
\newcommand{\Th}{{\operatorname{Th}}} 
 
\renewcommand{\th}{{\operatorname{th}}} 
\newcommand{\sHom}{\mathcal{H}om}
 
\newcommand{\Aut}{{\operatorname{Aut}}}

\newcommand{\GL}{\operatorname{GL}}
\newcommand{\SL}{\operatorname{SL}}

\newcommand{\can}{\text{can}}

\newcommand{\vir}{\text{\it vir}}
\newcommand{\perf}{{\operatorname{\it perf}}}

\newcommand{\ind}[1]{}

\newcommand{\inp}[1]{}

\newcommand{\EM}{{\operatorname{EM}}}

\renewcommand{\:}{{\colon}}

\makeatletter
 \@namedef{subjclassname@2020}{\textup{2020} Mathematics Subject Classification}
\makeatother

\begin{document}

\title{Virtual Localization in equivariant Witt cohomology}

\date{November 10, 2024}

\author[M.~Levine]{Marc~Levine}
\address{Universit\"at Duisburg-Essen,
Fakult\"at Mathematik, Campus Essen, 45117 Essen, Germany}
\email{marc.levine@uni-due.de}

%\subjclass[2020]{
%Primary 14B05, 14C15
%Secondary 14C30,   14F42}
%\keywords{}

\setcounter{tocdepth}{2}

\begin{abstract}  We prove an analog of the virtual localization theorem of Graber-Pandharipande, in the setting of an action by the normalizer of the torus in $\SL_2$, and with the Chow groups replaced by the cohomology of a suitably twisted sheaf of Witt groups. 
\end{abstract}
\thanks{The author was partially supported by the DFG through the grant  LE 2259/7-2 and by the ERC through the project QUADAG.  This paper is part of a project that has received funding from the European Research Council (ERC) under the European Union's Horizon 2020 research and innovation programme (grant agreement No. 832833).  \\
\includegraphics[scale=0.08]{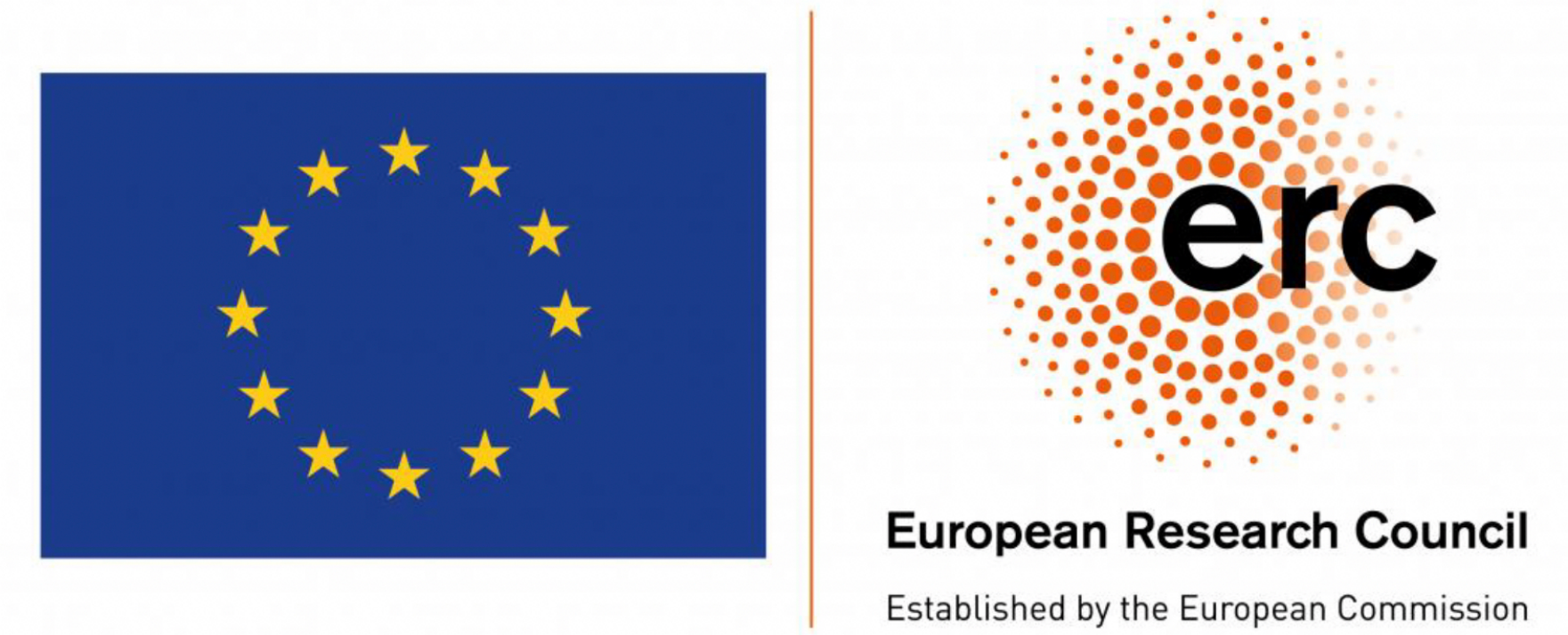}}

%\hskip 0.5\textwidth\begin{minipage}{0.5\textwidth} {\em Dedicated to the memory of Jacob Murre, for the joy and inspiration  his beautiful mathematics has given to me over my entire mathematical life.}
%\end{minipage}
%\\[0.5cm]

\maketitle

\tableofcontents

\section*{Introduction} The Bott residue theorem for a torus action on a variety $X$ over $\C$, as formulated by Edidin-Graham \cite{EGLoc}, gives an explicit formula for a class in the equivariant Chow groups of $X$ in terms of the  restriction of the class to the fixed point locus, after inverting a suitable product of equivariant Euler classes of line bundles in the equivariant Chow ring of a point. Graber and Pandharipande \cite{GP} have extended this to the case of a similar expression for the equivariant virtual fundamental class of $X$ associated to an equivariant perfect obstruction theory on $X$, in terms of the virtual fundamental class of the fixed point locus with respect to a certain  perfect obstruction theory induced by the one on $X$.

For applications to ``quadratically enhanced" theories, such as the Chow-Witt groups or cohomology of the sheaf of Witt groups,  torus localization is not useful, as inverting the Euler class of line bundles for such theories kills all the quadratic information. Instead, we
take $G=\SL_2^n$ or the normalizer $N$ of the diagonal torus $T_1$ in $\SL_2$.  We developed in \cite{LevineAtiyahBott} versions of Atiyah-Bott localization and a Bott residue formula for schemes with an action by these groups. Here the Euler class of equivariant line bundles gets replaced by the Euler classes of the  rank two bundles induced by  representations $G\to \GL_2$; localization with respect to these classes preserves at least some of the quadratic information.    

The main goal of this paper is to extend the virtual localization theorem of Graber-Pandharipande to this setting. We consider a $G$-linearized perfect obstruction theory $\phi_\bullet\:E_\bullet\to L_{X/B}$ on some $B$-scheme $X$ with $G$-action, and to this associate a virtual fundamental class $[X,\phi_\bullet]^\vir_{\sE, G}$ in the $G$-equivariant Borel-Moore homology $\sE^\BM_G(X/BG, E_\bullet)$. The equivariant Borel-Moore homology is defined as a limit, using the algebraic Borel construction of Totaro, Edidin-Graham and Morel-Voevodsky, and the virtual fundamental class is the one defined in \cite{LevineVirt} at each finite stage of the limit. 

For the virtual localization theorem, we restrict to the case of $G=N$ acting on a $k$-scheme $X$, $k$ a field. To explain the main players in the statement below, $\Sch^N/k$ is the category of quasi-projective $k$-schemes with an $N$-action, and $\Sm^N/k$ is the full subcategory of smooth $k$-schemes.  $|X|^N\subset X^{T_1}$ is the union of the $N$-stable irreducible components of $X^{T_1}$,  with  $i_j\:|X|^N_j\to X$ a certain union of connected components of $|X|^N$ determined by an equivariant closed immersion of $X$ in a smooth $Y$. The given $N$-linearized perfect obstruction theory 
$\phi\:E_\bullet\to  L_{X/k}$ on $X$ induces  a perfect obstruction theory $\phi_j\:
i_j^*E^\fix_\bullet\to L_{|X|^N_j/k}$ on each $|X|^N_j$,  and we have the Euler class of the virtual normal bundle,  $e_N(N^\vir_{i_j}):=e_N(\V(i_j^*E_0^\mov))\cdot e_N(\V(i_j^*E_1^\mov))^{-1}$. Here ${}^\fix$ and ${}^\mov$ refer to the ``fixed'' and ``moving'' parts of the corresponding sheaves with respect to the $T_1$-action. Also, the integer $M_0>0$ below is chosen so that the localization theorem \cite[Introduction, Theorem 5]{LevineAtiyahBott} for the $N$-action on $X$ holds after inverting $M_0e\in H^*(BN, \sW)$, where $e$ is the Euler class of the tautological rank two bundle on $BN$ corresponding to the inclusion $N\subset\SL_2\subset\GL_2$. See Definition~\ref{def:Strict} for the definition of a strict $N$-action; for the other notations  used in the statement, please see \S\ref{sec:VirLoc}. In the statement of our main theorem, we use the motivic spectrum $\EM(\sW_*)\in \SH(k)$ representing cohomology of the sheaf of Witt groups, via
\[
\EM(\sW_*)^{a,b}(X)=H^{a-b}(X, \sW)
\]
for $X\in \Sm/k$.

\begin{IntroThm}[Virtual localization, see \hbox{Theorem~\ref{thm:VirLoc}}]\label{IntroThm:Main} Let $k$ be a  field and let $\sE\in \SH(k)$ be the Eilenberg-MacLane spectrum $\sE:=\EM(\sW_*)$.  Let $i\:X\to Y$ be a closed immersion in $\Sch^N/k$, with $Y\in \Sm^N/k$, and let $\phi\:E_\bullet\to  L_{X/k}$ be an $N$-linearized perfect obstruction theory on $X$.  Suppose the $N$-action on $X$ is strict. We have the positive integer  $M=M_0\cdot \prod_{i,j}M^X_j \cdot M^Y_i$, where the $M_j^X, M_i^Y$ are as in 
Definition~\ref{def:NormalEulerClass} and $M_0$ is described above. Let $[|X|^N_j,\phi_j]_N^\vir\in \sE^\BM_N(|X|^N_j,  i_j^*E^\fix_\bullet)$ be the $N$-equivariant virtual fundamental class for the $N$-linearized perfect obstruction theory $\phi_j\: i_j^*E^\fix_\bullet\to L_{|X|^N_j/k}$ on $|X|^N_j$. Then
\[
[X,\phi]_N^\vir=\sum_{j=1}^s i_{j*}([|X|^N_j,\phi_j]_N^\vir\cap e_N(N^\vir_{i_j})^{-1})\in
\sE^\BM_N(X/BN,  E_\bullet)[1/Me].
\]
\end{IntroThm}
 
We recall the basic constructions and operations for cohomology and Borel-Moore homology in the setting of the theories represented by a motivic ring spectrum $\sE\in \SH(k)$ in  \S\ref{sec:Prelim}, following \cite{DJK} and \cite{LevineAtiyahBott}. The main technical heart of this paper is \S\ref{sec:Vistoli}, where we  prove a version (Proposition~\ref{prop:SpecializationComm}) of a result of Vistoli \cite[Lemma 3.16]{Vistoli}. Vistoli's lemma is about the Chow groups of DM stacks, whereas Proposition~\ref{prop:SpecializationComm} is suitable for application to the $G$-equivariant motivic theory.  We imagine that, with the help of the recent technology in the motivic stable homotopy theory of stacks provided by Chowdhury \cite{Chowdhury} and Khan-Ravi \cite{KR}, Proposition~\ref{prop:SpecializationComm} could be extended to motivic theories on a suitable type of stack. Kresch \cite[Proposition 4]{Kresch} has given a proof of Vistoli's lemma in the setting of the Chow groups of Artin stacks locally of finite type over a field, and perhaps his line of argument would be helpful in extending our result to the motivic setting for the type of Artin stacks considered by Chowdhury and Khan-Ravi.

One consequence of our motivic version of Vistoli's lemma is the commutativity of refined motivic Gysin pull-back, Corollary~\ref{cor:commutativity}, which as far as we are aware was not available in the literature up to now.

In \S\ref{sec:TEGEquivCoh} we recall the construction of equivariant cohomology and Borel-Moore homology in the style of Totaro \cite{Totaro},  Edidin-Graham \cite{EGEquivIntThy} and Morel-Voevodsky \cite{MorelVoevodsky}, following the work of Di Lorenzo and Mantovani \cite[\S 1.2]{DiLorenzoMantovani}. In \S\ref{sec:VirClass} we show how the construction of these classes as presented by Graber-Pandharipande also works in the motivic setting, and we apply this in 
\S\ref{sec:TEGVirClass} to the construction of virtual fundamental classes in equivariant 
Borel-Moore homology. We assemble all the ingredients   in \S\ref{sec:VirLoc}  with the statement and proof of our main result.

Theorem~\ref{IntroThm:Main} has been used by Anneloes Viergever \cite{AV} to compute some Witt ring-valued Donaldson-Thomas invariants for zero-dimensional subschemes of $\P^3$. Alessandro D'Angelo shows in his Ph.D. thesis \cite[Theorem 4.3.10]{DAngelo} how to extend Theorem~\ref{IntroThm:Main} from the case   $\sE=\EM(\sW_*)$ to an arbitrary $\SL$-oriented theory $\sE$ on which the algebraic Hopf map $\eta$ acts invertibly.

I wish to thank Fangzhou Jin for his very useful explanations of how to apply a number of results from \cite{DJK}. I would also like to thank the referee for their very helpful comments and suggestions, especially for pointing out similarities between parts of the proof of our version of Vistoli's lemma in \S\ref{sec:Vistoli}, and some of the arguments used by Rost \cite[\S10, 11, 13]{Rost} in proving the functoriality of pullback for his cycle complexes and the cohomology of cycle modules.

\section{Preliminaries and background}\label{sec:Prelim}

\subsection{Cohomology and Borel-Moore homology}\label{sec:CohBMHom} 
Fix a base-scheme $B$, quasi-projective over a noetherian ring $A$ of finite Krull dimension, and an  affine group-scheme $G$ over $B$.  We let $\Sch^G/B$ denote the category of quasi-projective $B$-schemes with $G$-action and let $\Sm^G/B$ be the full subcategory of smooth $B$-schemes  with $G$-action. Here we follow Hoyois \cite[\S1.2]{Hoyois} by defining the objects of $\Sch^G/B$ to be $B$-schemes $X$ endowed with a $G$-action, that admit a $G$-equivariant locally closed immersion into the projectivization of a $G$-equivariant vector bundle over $B$. The morphisms in $\Sch^G/B$ are just $G$-equivariant morphisms of $B$-schemes with $G$-action, ignoring the above-mentioned immersion.   We usually let $p_X:X\to B$ denote the structure morphism for $X\in \Sch^G/B$. 

In this section and in \S\ref{sec:Vistoli}, we assume that $G$ is {\em tame} \cite[Definition 2.26]{Hoyois}; for our ultimate purpose, the case of the trivial  $G$ will suffice, so this is not an essential limitation. We will drop the tame assumption in later sections.\footnote{The assumption that $B$ is quasi-projective over a noetherian ring of finite Krull dimension ensures that for tame $G$, $B$ has the $G$-resolution property of \cite[\S1.2]{Hoyois}}

Following Hoyois \cite{Hoyois}, for $X\in \Sch^G/B$, we have the  $G$-equivariant motivic stable homotopy category $\SH^G(X)$, together with a Grothendieck six-functor formalism. This says that $\SH^G(X)$ has a symmetric monoidal product $-\otimes_X-$ with adjoint  internal Hom  $\sHom(-,-)$, and for each morphism $f:Y\to X$ in $\Sch^G/B$, there are two pairs of adjoint functors
\[
f^*\:SH^G(X)\xymatrix{\ar@<3pt>[r]&\ar@<3pt>[l]}\SH^G(Y)\:f_*,\ 
f_!\:\SH^G(Y)\xymatrix{\ar@<3pt>[r]&\ar@<3pt>[l]}\SH^G(X)\:f^!,
\]
this package satisfying the properties described in \cite[Theorem 1.1]{Hoyois}. 

We let $\sH^G_*(X)$,  denote the pointed $G$-equivariant motivic unstable homotopy category, with infinite $\P^1$-suspension functor $\Sigma^\infty_{\P^1}: \sH^G_*(X)\to \SH^G(X)$.  Adjoining a disjoint base-point to $X\in \Sm^G/B$ gives the functor $\Sm^G/B \to \sH^G_*(X)$, $X\mapsto X_+$.  

We  refer the reader to \cite{LevineAtiyahBott} for our notation and basic facts about $\SH^G(-)$ extracted from  Hoyois, especially \cite[Theorem 1.1]{Hoyois}; we also use some constructions and results of D\'eglise-Jin-Khan \cite{DJK}. For use throughout the paper, we recall some of these facts here.

\begin{remark} The paper \cite{DJK} is set in the framework of the motivic stable category $\SH(-)$, not the $G$-equivariant version of Hoyois. However, all the constructions and results of \cite{DJK} carry over to the equivariant case (still for $G$ tame), and we will use this extension throughout this section.  In later sections, we only use  $\SH(-)$, so the cautious reader may safely restrict to that case.
\end{remark}

We have the $K$-theory space of $G$-linearized perfect complexes on $X$, $K^G(X)$, the associated fundamental groupoid $\sK^G(X)$ and the functor
\[
\Sigma^{(-)}:\sK^G(X)\to \Aut(\SH^G(X)) 
\]
which associated to   $v\in \sK^G(X)$, the suspension functor $\Sigma^v\:\SH^G(X)\to \SH^G(X)$. For $\sV$ a locally free $G$-linearized sheaf on $X$, we have the $G$-equivariant vector bundle  $\V(\sV):=\Spec_{\sO_X}(\Sym^*\sV)$ with projection $p\:\V(\sV)\to X$ and zero-section $s\:X\to \V(\sV)$.  Then
\[
\Sigma^\sV:=p_\#\circ s_*,\ \Sigma^{-\sV}:=s^!\circ p^*. 
\]
For each distinguished triangle in $D^\perf(X)$, $v'\to v\to v''\to v'[1]$, we have the canonical isomorphism
\[
\Sigma^v\cong \Sigma^{v'}\circ\Sigma^{v''}
\]
and for $v\in \sK^G(X)$, the canonical isomorphisms $\Sigma^{v[1]}\cong \Sigma^{-v}\cong (\Sigma^v)^{-1}$.

For each morphism $f\:Y\to X$ in $\Sch^G/B$ we have the natural transformation 
\begin{equation}\label{eqn:Lower!Lower*}
\alpha_f\:f_!\to f_*, 
\end{equation}
which is a natural isomorphism if $f$ is proper. 

Let $g:Y\to X$ be a smooth morphism in $\Sch^G/B$. Then $g^*:\SH^G(X)\to \SH^G(Y)$ has the left adjoint $g_\#:\SH^G(Y)\to \SH^G(X)$ and we have the {\em purity isomorphisms}
\begin{equation}\label{eqn:PurityIso}
g_!\cong g_\#\circ\Sigma^{-\Omega_g},\ g^!\cong \Sigma^{\Omega_g}\circ g^*,
\end{equation}
where $\Omega_g$ is the sheaf of relative K\"ahler differentials.

\begin{definition}[Twisted $\sE$-cohomology]
For $p_X:X\to B$ in $\Sch^G/B$,  $\sE\in \SH^G(B)$ and $v\in \sK^G(X)$, we have the twisted $\sE$-cohomology
%\footnote{Throughout this paper we use Hoyois' convention for the suspension functors, twisted cohomology and Borel-Moore homology,  using a locally free coherent sheaf for the twist. The D\'eglise-Jin-Khan convention in \cite{DJK} often uses the vector bundle $V=\V(\sV)$ instead.  For instance, we (and Hoyois) write the purity isomorphism for a smooth morphism $f:Y\to X$ as $f_!\cong f_\#\circ \Sigma^{-\Omega_f}$ which is written as $f_!\cong f_\#\circ \Sigma^{-T_f}$ ($T_f=\V(\Omega_f)$) in \cite{DJK}. }
\[
\sE^{a,b}(X, v):=\Hom_{\SH(X)}(1_X, \Sigma^{a,b}\Sigma^v p_X^*\sE).
\]
\end{definition}
Twisted $\sE$-cohomology extends to a functor on $\Sch^G/X$ by applying the functor $f^*:\SH^G(X)\to \SH^G(Y)$ for each $f:(Y'\xrightarrow{\pi_{Y'}}X)\to (Y\xrightarrow{\pi_{Y}}X)$ in $\Sch^G/X$, where we write  $\sE^{a,b}(Y, v)$ for $\sE^{a,b}(Y, \pi_Y^*v)$:
\begin{multline*}
\sE^{a,b}(Y, v)=\Hom_{\SH(Y)}(1_Y, \Sigma^{a,b}\Sigma^{\pi_Y^*v} p_Y^*\sE)\xrightarrow{f^*}
\Hom_{\SH(Y)}(f^*1_Y, f^*\Sigma^{a,b}\Sigma^{\pi_Y^*v} p_Y^*\sE)\\
=\Hom_{\SH(Y)}(1_{Y'}, \Sigma^{a,b}\Sigma^{\pi_{Y'}^*v} p_{Y'}^*\sE)=
\sE^{a,b}(Y', v)
\end{multline*}

For $i:Z\to Y$ a closed immersion in $\Sch^G/B$ with open complement $j:U\to Y$, we have the {\em localization distinguished triangle} in $\SH(Y)$,
\begin{equation}\label{eqn:LocTriangle}
j_!j^!\xrightarrow{v^!_!(j)}\id_{\SH(Y)}\xrightarrow{u^*_*(i)}i_*i^*\xrightarrow{\del}j_!j^![1],
\end{equation}
where $v^!_!(j)$ and $u^*_*(i)$ are the counit and unit of the respective adjunctions. Using the purity isomorphism, we can rewrite this as
\[
j_\#j^*\xrightarrow{v^!_!(j)}\id_{\SH(Y)}\xrightarrow{u^*_*(i)}i_*i^*\xrightarrow{\del}j_\#j^*[1].
\]
Applying this to $1_Y$ gives the canonical isomorphism
\begin{equation}\label{eqn:LocQuotient}
\Sigma^\infty_{\P^1}(Y/U)\cong i_*1_Z
\end{equation}
in $\SH^G(Y)$.  

Let $\sV$ be a locally free $G$-linearized coherent sheaf on $X\in \Sch^G/B$, giving the $G$-equivariant vector bundle $p_V:V:=\V(\sV)\to X$ with zero section $s_V:X\to V$, and open complement $j_V:V\setminus\{s(X)\}\to V$. The quotient $V/(V\setminus\{s_V(X)\}\in \sH_*^G(X)$ is the {\em Thom space} $\Th_X(V)$.

Applying \eqref{eqn:LocQuotient} to $s_V:X\to V$ and acting by $p_{V\#}:\SH^G(V)\to \SH^G(X)$ we have the canonical isomorphism
\begin{equation}\label{eqn:StableThomSpace}
\Sigma^\infty_{\P^1}\Th_X(V)\cong p_{V\#}s_{V*}1_X=\Sigma^\sV1_X
\end{equation}
in $\SH^G(X)$.  

 Composing $s_{V+}:X_+\to V_+$ with the quotient map $V_+\to \Th_X(V)$,  applying $\Sigma^\infty_{\P^1}$ and using the isomorphism \eqref{eqn:StableThomSpace}  gives the morphism $\bar{s}_V:1_X\to \Sigma^\sV1_X$ in $\SH^G(X)$.

If we have a morphism $f:Y\to X$, applying $f^*:\sH_*^G(X)\to:\sH_*^G(Y)$ to  the cofiber sequence
\[
V\setminus\{s_V(X)\}\to V\to \Th_X(V),
\]
and noting that $f^*(W)=W\times_XY$ for $W\in \Sm^G/X$,
gives the cofiber sequence 
\[
f^*V\setminus\{s_{f^*V}(X)\}\to f^*V\to f^* \Th_X(V)
\]
identifying $f^*\Th_X(V)$ with $\Th_X(f^*V)$. Passing to $f^*:\SH^G(X)\to \SH^G(Y)$ and using the canonical isomorphism $f^*(\Sigma^\sV1_X)=\Sigma^{f^*\sV}1_Y$, this gives
\begin{equation}\label{eqn:ThomSpaceId}
f^*(\bar{s}_V)=\bar{s}_{f^*V}:1_Y\to \Sigma^{f^*\sV}1_Y.
\end{equation}

We now define the Euler class of a $G$-vector bundle, following  \cite[\S3.1]{DJK}.
\begin{definition}[Euler class]\label{def:EulerClass} Let $\sV$ be a locally free $G$-linearized coherent sheaf on $X\in \Sch^G(B)$, giving the $G$-equivariant vector bundle $V:=\V(\sV)\to X$. Let $\sE\in \SH^G(B)$ be a commutative ring spectrum with unit map $u_\sE:1_B\to \sE$.  The composition
\[
1_X\xrightarrow{\bar{s}_V} \Sigma^{\sV}1_X\xrightarrow{\id\wedge p_X^*u_\sE}
\Sigma^{\sV}1_X\wedge p_X^*\sE=\Sigma^{\sV} p_X^*\sE
\]
in $\SH^G(X)$
defines the {\em cohomological Euler class} $e^\sE(V)\in \sE(X,\sV)$.
\end{definition}

\begin{remark}\label{rem:EulerNat}
Given $X$, $\sV$ and $V=\V(\sV)\to X$ as above, let $f:Y\to X$ be a morphism in $\Sch^G/B$, giving the locally free $G$-linearized coherent sheaf $f^*\sV$ on $Y$ with associated vector bundle $\V(f^*\sV)=f^*V\to Y$. By \eqref{eqn:ThomSpaceId}, we have
\begin{equation}\label{eqn:EulerClassNat}
f^*(e^\sE(V))=e^\sE(f^*V).
\end{equation} 
\end{remark}

\begin{definition}[Twisted $\sE$-Borel-Moore homology] Let $\sE\in \SH^G(B)$ be a commutive ring spectrum, and take $X\in \Sch^G/B$, $v\in \sK^G(X)$. The {\em $v$-twisted $\sE$-Borel-Moore homology of $X$} is defined as
\[
\sE^\BM_{a,b}(X/B, v):=\Hom_{\SH^G(B)}(\Sigma^{a,b}p_{X!}\Sigma^v 1_X, \sE)=
\Hom_{\SH^G(X)}(\Sigma^{a,b}\Sigma^v 1_X, p_X^!\sE).
\]
We will write $\sE^\BM(X/B,v)$ for $\sE^\BM_{0,0}(X/B,v)$; if we need to include mention of $G$, we will write this as $\sE^\BM_{a,b}(X/B,G, v)$.  For a morphism $f\:Y\to X$ and a $v\in K^G(X)$, we often write $\sE^\BM_{a,b}(Y/B,v)$  for $\sE^\BM_{a,b}(Y/B,f^*v)$  and often write $\sE^\BM_{a,b}(X,v)$ for $\sE^\BM_{a,b}(X/B,v)$ if the base-scheme $B$ is clear from the context. 
\end{definition}

The purity isomorphism \eqref{eqn:PurityIso} for $g=p_X:X\to B\in \Sm^G/B$ gives the {\em Poincar\'e duality isomorphism}
\begin{equation}\label{eqn:PDIso}
\sE^{a,b}(X, v)\cong \sE^\BM_{-a,-b}(X/B, \Omega_{X/B}-v).
\end{equation}

\subsection{Cup and cap products}\label{subsec:Products}  

We recall the following properties from the six-functor formalism (see \cite[Theorem 1.1]{Hoyois}). 

Let $f:Y\to X$ be a morphism in $\Sch^G/B$. Then  
\begin{align}
&\text{(Monoidality) There is a natural isomorphism $f^*((-)\otimes(-))\cong f^*(-)\otimes f^*(-)$} 
\label{align:Product}\\
&\text{(Projection formula) 
There is a natural isomorphism  $f_!(f^*(-)\otimes (-))\cong (-)\otimes f_!(-)$. }\label{align:ProjFormula}
 \end{align}
We have the natural map 
\begin{equation}\label{eqn:ExceptPullback}
f^*(-)\otimes f^!(-)\to f^!((-)\otimes(-))
\end{equation}
adjoint to the map $f_!(f^*(-)\otimes f^!(-))\to (-)\otimes(-)$ defined by
\[
f_!(f^*(-)\otimes f^!(-))\cong (-)\otimes f_!f^!(-)\xrightarrow{\id\otimes v^!_!(f)}(-)\otimes(-)
\]
where the isomorphism is \eqref{align:ProjFormula} and $v^!_!(f)$ is the counit of adjunction.

The cup product for $X\in \Sch^G/B$, $v,w\in \sK^G(X)$,
\[
\cup_X: \sE^{a,b}(X, v)\times \sE^{c,d}(X, w)\to \sE^{a+c, b+d}(X, v+w),
\]
is constructed by applying the monoidal product  
\begin{multline*}
\Hom_{\SH(X)}(1_X, \Sigma^{a,b}\Sigma^vp_X^*\sE)\times \Hom_{\SH(X)}(1_X, \Sigma^{c,d}\Sigma^wp_X^*\sE)\\\to \Hom_{\SH(X)}(1_X, \Sigma^vp_X^*\sE\otimes \Sigma^wp_X^*\sE)
\end{multline*}
followed by 
\begin{multline*}
\Sigma^{a,b}\Sigma^vp_X^*\sE\otimes \Sigma^{a,b}\Sigma^wp_X^*\sE=\Sigma^{a+c,b+d}\Sigma^{v+w}p_X^*\sE\otimes p_X^*\sE\\\cong \Sigma^{a+c,b+d}\Sigma^{v+w}p_X^*(\sE\otimes \sE)\xrightarrow{p_X^*(\mu_\sE)}\Sigma^{a+c,b+d}\Sigma^{v+w}p_X^*\sE,
\end{multline*}
where $\mu_\sE$ is the multiplication. 

Cap products
\[
\cap_X: \sE^{a,b}(X, v)\times \sE^\BM_{c,d}(X/B, w)\to \sE^{c-a, d-b}(X/B, w-v)
\]
are similarly constructed using the map
\begin{multline*}
\Sigma^{a,b}\Sigma^vp_X^*\sE\otimes \Sigma^{-c,-d}\Sigma^{-w}p_X^!\sE=\Sigma^{a-c,b-d}\Sigma^{v-w}p_X^*\sE\otimes p_X^!\sE\\\xrightarrow{\eqref{eqn:ProjFormula}} \Sigma^{a-c,b-d}\Sigma^{v-w}p_X^!(\sE\otimes \sE)\xrightarrow{p_X^!(\mu_\sE)}\Sigma^{a-c,b-d}\Sigma^{v-w}p_X^!\sE.
\end{multline*}

\begin{remark} For notational convenience, we sometime switch the order in the cap product, with the cohomology on the right and the Borel-Moore homology on the left, and do the same for the other products defined here and  below.
\end{remark}

\subsection{Functorialities}
A proper map $f\:Y\to X$ in $\Sch^G/B$ gives rise to the natural transformation $f^*_\BM\:p_{X!}\to p_{Y!}f^*$, defined as the composition
\[
p_{X!}\xrightarrow{u^*_*(f)}p_{X!}f_*f^*\xrightarrow{\alpha_f^{-1}}p_{X!}f_!f^*=p_{Y!}f^*,
\]
where $u^*_*(f):\id\to f_*f^*$ is the unit of   adjunction. 
This induces the {\em proper push-forward map}
\begin{equation}\label{eqn:ProperPushforward}
f_*\:\sE^\BM_{a,b}(Y, v)\to \sE^\BM_{a,b}(X,v)
\end{equation}
by pre-composition with $f^*_\BM$, applied to $\Sigma^{a,b}\Sigma^v1_X$.

For $g:Y\to X$ a smooth morphism in $\Sch^G/B$, we have the {\em smooth pullback map}
\begin{equation}\label{eqn:SmoothPullback}
g^!\:\sE^\BM_{a,b}(X, v)\to \sE^\BM_{a,b}(Y, v+\Omega_g)
\end{equation}
formed by applying the functor $g^!:\SH^G(X)\to \SH^G(Y)$ and using the  purity isomorphism $g^!1_X\cong \Sigma^{\Omega_g}g^*1_X=\Sigma^{\Omega_g}1_Y$:
\begin{multline*}
\sE^\BM_{a,b}(X, v)=\Hom_{\SH^G(X)}(\Sigma^{a,b}\Sigma^v 1_X, p_X^!\sE)\\
\xrightarrow{g^!}\Hom_{\SH^G(Y)}(g^!\Sigma^{a,b}\Sigma^v 1_X, g^!p_X^!\sE)=
\Hom_{\SH^G(Y)}(\Sigma^{a,b}\Sigma^{g^*v} g^!1_X,  p_Y^!\sE)\\
\cong \Hom_{\SH^G(Y)}(\Sigma^{a,b}\Sigma^{g^*v+\Omega_g}1_Y,  p_Y^!\sE)=
\sE^\BM_{a,b}(Y, v+\Omega_g).
\end{multline*}

\begin{proposition} \label{prop:PushPullComp}
1. Both proper pushforward and smooth pullback are functorial\\[2pt]
2. Let
\[
\xymatrix{
W\ar[r]^{f_W}\ar[d]^{g_Z}&Y\ar[d]^g\\
Z\ar[r]^{f}&X
}
\]
be a cartesian square in $\Sch^G/B$, with $f$ proper and $g$ smooth, and take  $\sE\in \SH^G(B)$ and  $v\in \sK^G(X)$. Then $f_W$ is proper, $g_Z$ is smooth, $f_W^*\Omega_g\cong \Omega_{g_Z}$,  and the diagram
\[
\xymatrix{
\sE^\BM_{a,b}(Z, v)\ar[r]^{f_*}\ar[d]^{g_Z^!}&\sE^\BM_{a,b}(X, v)\ar[d]^{g^!}\\
\sE^\BM_{a,b}(W, v+\Omega_{g_Z})\ar[r]^{f_{W*}}&\sE^\BM_{a,b}(Y, v+\Omega_g)
}
\]
commutes.
\end{proposition}

\begin{proof} Given proper morphisms $Z\xrightarrow{g}Y\xrightarrow{f}X$ in $\Sch^G/B$ functoriality of proper pushforward follows from the commutativity of the diagram
\[
\xymatrix{
p_{X!}\ar[r]^{u^*_*(f)}\ar[d]^{u^*_*(fg)}&p_{X!}f_*f^*\ar[d]^{u^*_*(g)}\ar[r]^{\alpha_f^{-1}}&p_{X!}f_!f^*\ar@{=}[r]\ar[d]^{u^*_*(g)}&p_{Y!}f^*\ar[d]^{u^*_*(g)}\\
p_{X!}(fg)_*(fg)^*\ar[d]^{\alpha_{fg}^{-1}}\ar@{=}[r]&p_{X!}f_*g_*g^*f^*\ar[r]^{\alpha^{-1}_f}&
p_{X!}f_!g_*g^*f^*\ar@{=}[r]\ar[d]^{\alpha^{-1}_g}&p_{Y!}g_*g^*f^*\ar[d]^{\alpha_g^{-1}}\\
p_{X!}(fg)_!(fg)^*\ar@{=}[rr]\ar@{=}[d]&&p_{X!}f_!g_!g^*f^*\ar@{=}[r]&p_{Y!}g_!g^*f^*\ar@{=}[d]\\
p_{Z!}(fg)^*\ar@{=}[rrr]&&&p_{Z!}(fg)^*
}
\]
Here $u^*_*(f):\id\to f_*f^*$, $u^*_*(g)$ and $u^*_*(fg)$ are the units of adjunction, and the identity in the upper right square is just $(f_*u^*_*(g)f^*)\circ u^*_*(f)=u^*_*(fg)$, which is a consequence of the adjunction. Similarly, the identity $(\id\times \alpha^{-1}_g)\circ( \alpha^{-1}_f\times \id)=\alpha^{-1}_{fg}$ (after applying the identifications $f_!g_!=(fg)_!$, $f_*g_*=(fg)_*$) follows  from the construction of $\alpha_{(-)}$ in \cite[\S6.2]{Hoyois}. 

Functoriality of smooth pullback follows easily from the six-functor formalism, in particular, the identity $(fg)^!=g^!f^!$ for composable morphisms $f, g$ in $\Sch^G/B$, plus the  exact sequence $0\to g^*\Omega_f\to \Omega_{fg}\to \Omega_g\to 0$ for composable smooth morphisms $f,g$.

(2) is proven in, e.g., \cite[Lemma 2.3]{LevineVirt}, with  a different notation for the smooth pullback. The definition of the smooth pullback map given there is not the same as what is given here, but the two maps are identical, which one sees using the adjunction $(-)_!\dashv (-)^!$.
\end{proof}

\begin{remark} For a smooth morphism $g:Y\to X$ in $\Sm^G/B$, the respective Poincar\'e duality isomorphisms for $X$ and $Y$ intertwine $g^*$ and $g^!$.
\end{remark}

For $p_Y:Y\to B$ in $\Sm^G/B$, we have the external cap product
\begin{equation}\label{eqn:ExtCapProd}
\cap_{Y,X}: \sE^{a,b}(Y, v)\times \sE^\BM_{c,d}(X/B, w)\to \sE^{c-a, d-b}(Y\times_BX/B, w-v+\Omega_{Y/B})
\end{equation}
defined by 
\[
a\cap_{Y,X} b:=p_1^*(a)\cap_{Y\times_BX} p_2^!(b),
\]
noting that $p_2$ is smooth with $L_{p_2}=p_1^*\Omega_{Y/B}$. Using the  purity isomorphism $\sE^\BM_{a,b}(Y, v)\cong \sE^{-a,-b}(Y, \Omega_{Y/B}-v)$, \eqref{eqn:ExtCapProd} gives the external product in Borel-Moore homology
\begin{equation}\label{eqn:ExternalProd}
\boxtimes_{Y,X}:\sE^\BM_{a,b}(Y, v)\times \sE^\BM_{c,d}(X/B, w)\to \sE^\BM_{a+c, b+d}(Y\times_BX/B, w+v). 
\end{equation}

Let $i:Z\to X$ be a closed immersion in $\Sch^G/B$ with open complement $j:U=X\setminus Z\to X$ and take $v\in \sK^G(X)$. Evaluating the localization distinguished triangle
\[
j_!j^!\to \id_{\SH^G(X)}\to i_*i^*\to j_!j^![1]
\]
at $\Sigma^{a,b}\Sigma^v1_X$, applying $p_{X!}$ and using the isomorphisms $i_!\cong i_*$, $j^!\cong j^*$ gives the localization triangle in $\SH(B)$
\[
p_{U!}\Sigma^{a,b}\Sigma^{j^*v}1_U\to p_{X!}\Sigma^{a,b}\Sigma^v1_X\to p_{Z!}\Sigma^{a,b}\Sigma{i^*v}1_Z\to p_{U!}\Sigma^{a,b}\Sigma^{j^*v}1_U[1].
\]
Applying $\Hom_{\SH^G(B)}(-, \sE)$ gives the long exact localization sequence
\begin{equation}\label{eqn:LocSeq}
\ldots\xrightarrow{i_*}\sE^\BM_{a,b}(X, v)\xrightarrow{j^!} \sE^\BM_{a,b}(U, v)
\xrightarrow{\del}\sE^\BM_{a+1,b}(Z, v)\xrightarrow{i_*}\sE^\BM_{a+1,b}(X, v)\xrightarrow{j^!} \ldots
\end{equation}

The localization sequence is compatible with smooth pullback.
\begin{proposition}\label{prop:LocSmoothComp} Let $i:Z\to X$ be a closed immersion in $\Sch^G/B$ with open complement $j:U\to X$. Let $g:X'\to X$ be a smooth morphism in $\Sch^G/B$, and let  $Z'=Z\times_XX'\xrightarrow{i'} X'$ and $U'=U\times_XX'\xrightarrow{j'}X'$ be the induced closed immersion and open complement and let  $f_Z:Z'\to Z$, $f_U:U'\to U$ be the projections. Take $\sE\in \SH^G(B)$ and $v\in \sK^G(X)$ and let $v'=f^*v+\Omega_f$. Then the diagram
\[
\xymatrix{
\sE^\BM_{a,b}(X, v)\ar[r]^{j^!}\ar[d]^{f^!}& \sE^\BM_{a,b}(U, v)
\ar[r]^{\del}\ar[d]^{f_U^!}&\sE^\BM_{a+1,b}(Z, v)\ar[r]^{i_*}\ar[d]^{f_Z^!}&
\sE^\BM_{a+1,b}(X, v)\ar[d]^{f^!}\\
\sE^\BM_{a,b}(X', v')\ar[r]^{j^{\prime!}}& \sE^\BM_{a,b}(U', v')
\ar[r]^{\del}&\sE^\BM_{a+1,b}(Z', v')\ar[r]^{i'_*}&\sE^\BM_{a+1,b}(X', v')
}
\]
commutes for all $a,b$.
\end{proposition}
\begin{proof} The commutativity of the left-had square is the functoriality of $(-)^!$ (Proposition~\ref{prop:PushPullComp}(1)), and the commutativity of the right-hand square is Proposition~\ref{prop:PushPullComp}(2). We then use \cite[Proposition 2.3.3]{CD} to take care of the middle square.
\end{proof}

\subsection{Specialization}\label{sec:specialization} 
We recall the construction of the {\em specialization map} (see e.g. \cite[Definition3.2.4]{DJK}).
Let $\iota\: X\hookrightarrow Y$ be a closed immersion of $B$-schemes with associated ideal sheaf $\sI_X\subset \sO_Y$. We have the blow-up $\Bl_{X\times0}Y\times\A^1\to Y\times\A^1$ and the {\em deformation space} $\Def_{X}Y:=\Bl_{X\times0}Y\times\A^1\setminus \Bl_{X\times0}Y\times0$, with open immersion $j\:Y\times \A^1\setminus\{0\}\to \Def_{X}Y$ and closed complement $i\:C_{X/Y}\to \Def_{X}Y$; here $C_{X/Y}$ is the normal cone of the immersion $\iota$, 
\[
C_{X/Y}:=\Spec_{\sO_X}\oplus_{n\ge0}\sI_X^n/\sI_X^{n+1}.
\]

Take $v$ in $\sK^G(Y)$. Let 
\[
\sp_{\iota/B,v}\:p_{C_{X/Y}!}\circ \Sigma^v\circ  p_{C_{X/Y}}^*\to  p_{Y!}\circ\Sigma^v\circ  p_Y^*
\]
be the natural transformation defined as follows. Start with the localization distinguished triangle
\[
j_!j^!\to \id_{\SH^G(\Def_{X}Y)}\to i_*i^*\xrightarrow{\del} j_!j^![1].
\]
Composing $\del$ with $p_{\Def_{X}Y!}$ on the left and  $\Sigma^v\circ  p_{\Def_{X}Y}^*$ on the right gives
\[
\del\:p_{C_{X/Y}!}\circ \Sigma^v\circ  p_{C_{X/Y}}^*\to
p_{Y\times(\A^1\setminus\{0\})!}\circ\Sigma_{S^1}\circ \Sigma^v\circ  p_{Y\times(\A^1\setminus\{0\})}^*,
\]
with $\del^*$ giving  the boundary map in the localization sequence for $C_{X/Y}\hookrightarrow \Def_{X}Y$.

We have $p_{Y\times(\A^1\setminus\{0\})!}=p_{Y!}\circ p_{1!}\cong p_{Y!}\circ p_{1\#}\circ\Sigma^{-\Omega_{p_1}}$, using the purity isomorphism \eqref{eqn:PurityIso}. Writing $Y\times (\A^1\setminus\{0\})=\Spec_{\sO_Y}\sO_Y[x^{\pm1}]$ gives us the generator $dx$ for $\Omega_{p_1}$, which gives us the canonical isomorphism $\Sigma^{-\Omega_{p_1}}\cong \Sigma^{-1}_{\P^1}$, and yields the identity
\begin{equation}\label{eqn:AnIdentity}
p_{Y\times(\A^1\setminus\{0\})!}\circ\Sigma_{S^1}\circ \Sigma^v\circ p_{Y\times(\A^1\setminus\{0\})}^*=
p_{Y!}\circ   (p_{1\#}\circ p_1^*\circ \Sigma^{-1}_{\P^1})\circ \Sigma_{S^1}\circ  \Sigma^v\circ p_Y^*.
\end{equation}
The quotient map $(\A^1\setminus\{0\})_+\to \G_m:=(\A^1\setminus\{0\})/\{1\}$ induces the map  $p_{1\#}\circ p_1^*\cong \Sigma_{(\A^1\setminus\{0\})_+}\to \Sigma_{\G_m}$. Since 
$\Sigma_{\P^1}\cong \Sigma_{S^1}\circ\Sigma_{\G_m}$,
this gives  the map 
\begin{equation}\label{eqn:Vartheta}
\vartheta_{Y,v}\: p_{Y\times\A^1\setminus\{0\}!}\circ\Sigma_{S^1}\circ\Sigma^v\circ   p_{Y\times\A^1\setminus\{0\}}^*\to p_{Y!}\circ \Sigma_{S^1}^{-1}\circ\Sigma_{S^1}\circ \Sigma^v\circ   p_Y^*=
 p_{Y!}\circ \Sigma^v\circ   p_Y^*,
\end{equation}
and we set  
\[
\sp_{\iota/B,v}:=\vartheta_{Y,v}\circ\del\:p_{C_{X/Y}!}\circ \Sigma^v\circ  p_{C_{X/Y}}^*\to  p_{Y!}\circ \Sigma^v\circ   p_Y^*.
\]
Applying this to $1_Y$ and taking $\Hom_{\SH^G(B)}(\Sigma^{**}(-),\sE)$ gives the specialization map in $\sE$-Borel-Moore homology,
\[
\sp^*_{\iota/B,v}:=\del^*\circ \vartheta_{Y,v}^\:\sE^\BM_{**}(Y/B, v)\to \sE^\BM_{**}(C_{X/Y}/B,v),
\]

We will often drop the $-/B$ in the notation if the base-scheme is clear, writing 
$\sp_{\iota,v}$, $\sp^*_{\iota,v}$ for $\sp_{\iota/B,v}$, $\sp^*_{\iota/B,v}$. 

The specialization map commutes with smooth pull-back in cartesian squares.  

\begin{proposition}\label{prop:SpecSmoothCommute} Let $\iota\:X\to Y$ be a closed immersion and let $f\:Y'\to Y$ be a smooth morphism, both  in $\Sch^G/B$. Form the cartesian diagram
\[
\xymatrix{
X'\ar[r]^{f_X}\ar[d]^{\iota'}&X\ar[d]^\iota\\
Y'\ar[r]^f&Y,
}
\]
inducing the smooth morphism $f_C:C_{X'/Y'}\to C_{X/Y}$.
Then for $v\in \sK^G(Y)$ and $\sE\in \SH^G(B)$, the diagram
\[
\xymatrixcolsep{70pt}
\xymatrix{
\sE^\BM_{a,b}(Y, v)\ar[r]^{\sp^*_{\iota,v}}\ar[d]^{f^!}& \sE^\BM_{a,b}(C_{X/Y},v)\ar[d]^{f^!_C}\\
\sE^\BM_{a,b}(Y', v+\Omega_f)\ar[r]^{\sp^*_{\iota',v+\Omega_f}}& \sE^\BM_{a,b}(C_{X'/Y'},v+\Omega_f)
}
\]
commutes.
\end{proposition}

\begin{proof} If we apply the adjunctions 
\begin{gather*}
\Hom_{\SH(Y)}(\Sigma^{**}\Sigma^v1_Y, p_Y^!\sE)\cong \Hom_{\SH(B)}(p_{Y!}\Sigma^{**}\Sigma^v1_Y, \sE)\\
\Hom_{\SH(Y')}(\Sigma^{**}\Sigma^{v+\Omega_f}1_{Y'}, p_{Y'}^!\sE)\cong \Hom_{\SH(B)}(p_{Y'!}\Sigma^{**}\Sigma^{v+\Omega_f}1_{Y'}, \sE),
\end{gather*}
the smooth pullback by $f$ is given by applying the natural transformation 
\[
p_{Y'!}\Sigma^{**}\Sigma^{v+\Omega_f}f^*
\cong p_{Y!}f_!f^!\Sigma^{**}\Sigma^v\xrightarrow{v^!_!(f)}p_{Y!}\Sigma^{**}\Sigma^v
\]
to $1_Y$, and then precomposing with the resulting map $p_{Y'!}\Sigma^{**}\Sigma^{v+\Omega_f}1_{Y'}\to p_{Y!}\Sigma^{**}\Sigma^v1_Y$. The smooth pullback with respect to $f$ and $f\times\id_{\A^1\setminus\{0\}}$  intertwines the two maps
$\vartheta_{Y,v}$ and $\vartheta_{Y',v+\Omega_f}$, as one sees from the commutative diagram
\[
\xymatrix{
p_{Y'\times(\A^1\setminus\{0\})!}\Sigma_{S^1}\Sigma^{v+\Omega_f}p_{Y'\times(\A^1\setminus\{0\})}^*\ar[r]^{\vartheta_{Y', v+\Omega_f}}\ar@{=}[d]&p_{Y'!}\Sigma^{v+\Omega_f}p_{Y'}^*\ar@{=}[d]\\
p_{Y\times(\A^1\setminus\{0\})!}\Sigma_{S^1}(f\times\id_{\A^1\setminus\{0\}})_!\Sigma^{\Omega_f}(f\times\id_{\A^1\setminus\{0\}})^*\Sigma^vp_{Y\times(\A^1\setminus\{0\})}^*\ar@{=}[d]\ar[r]^-{\vartheta_{Y,v}'}&
p_{Y!}f_!\Sigma^{\Omega_f}f^*\Sigma^vp_Y^*\ar@{=}[d]\\
p_{Y\times(\A^1\setminus\{0\})!}\Sigma_{S^1}(f\times\id_{\A^1\setminus\{0\}})_!(f\times\id_{\A^1\setminus\{0\}})^!\Sigma^vp_{Y\times(\A^1\setminus\{0\})}^*\ar[d]^{v^!_!(f\times\id_{\A^1\setminus\{0\}})}\ar[r]^-{\vartheta_{Y,v}''}&
p_{Y!}f_!f^!\Sigma^vp_Y^*\ar[d]^{v^!_!(f)}\\
p_{Y\times(\A^1\setminus\{0\})!}\Sigma_{S^1}\Sigma^vp_{Y\times(\A^1\setminus\{0\})}^*\ar[r]^{\vartheta_{Y,v}} &
p_{Y!}\Sigma^vp_Y^*
}
\]
Here the map $\vartheta_{Y,v}'$  is defined similarly to $\vartheta_{Y,v}$, using the fact that  $f\circ p_1=p_1\circ (f\times \id_{\A^1\setminus\{0\}})$, together with functoriality of $(-)^*$ and smooth base-change for $(-)^*$ and $(-)_!$, to give
\[
p_{1!}\circ (f\times\id_{\A^1\setminus\{0\}})_!\Sigma^{\Omega_f}(f\times\id_{\A^1\setminus\{0\}})^*\circ \Sigma^v p_1^*
=(p_{1!}\circ p_1^*)\circ f_!\Sigma^{\Omega_f}f^*\Sigma^v.
\]
The map $\vartheta_{Y,v}''$ is defined similarly, using  
\[
p_{1!}\circ (f\times\id_{\A^1\setminus\{0\}})_!(f\times\id_{\A^1\setminus\{0\}})^!\circ\Sigma^v p_1^*
=(p_{1!}\circ p_1^*)\circ f_!f^!\Sigma^v.
\]
By Proposition~\ref{prop:LocSmoothComp}, the  diagram
\[
\xymatrix{
\sE^\BM_{a+1,b}(Y\times\A^1\setminus\{0\}/B, v)\ar[r]^-{\del}\ar[d]^{(f\times \id)^!}&\sE^\BM_{a, b}(C_{X/Y}/B,v)\ar[d]^{f^!_C}\\
\sE^\BM_{a+1,b}(Y'\times\A^1\setminus\{0\}/B, v)\ar[r]^-{\del}&\sE^\BM_{a, b}(C_{X/Y}/B,v)
}
\]
commutes,  giving
\[
f^!_C\circ \sp^*_{\iota,v}=f^{\prime !}_C\circ \del\circ (\vartheta_{Y,v})^*=
\del\circ (\vartheta_{Y',v+\Omega_f})^*\circ f^! =   \sp^*_{\iota',v+\Omega_f}\circ f^!.
\]
\end{proof}

\begin{remark} For $\sN$ a locally free sheaf on some $Y$ with associated vector bundle  $p\:N\to Y$, $N:=\V(\sN)$,   the smooth pull-back
\[
p^!\:\sE^\BM_{a,b}(Y,v)\to \sE^\BM_{a,b}(N, v+\sN)
\]
is an isomorphism and one defines the pull-back by the zero section
\[
0_N^!\: \sE^\BM_{a,b}(N, v+\sN)\to \sE^\BM_{a,b}(Y,v)
\]
to be the inverse of $p^!$. Thus, specialization also commutes with pull-back by zero-sections in cartesian squares. 
\end{remark}

\subsection{Gysin morphisms}\label{sec:Gysin}  We   work over the fixed base-scheme $B$, unless  explicitly mentioned to the contrary, and we fix a commutative ring spectrum $\sE\in \SH^G(B)$.

For $X\in \Sch^G/B$, let $D_G(X)$ denote the derived category of $G$-linearized complexes of quasi-coherent sheaves on $X$. For a morphism $f\:X\to Y$ in $\Sch^G/B$, we have the cotangent complex $L_f\in D_G(X)$. $L_f$ agrees with the class of the relative differentials  $\Omega_f$ in case $f$ is smooth, and with $\sN_f[1]$ if $f$ is a regular immersion with  conormal sheaf $\sN_f:=\sI_X/\sI_X^2$.  

Let $\iota\:Y_0\to Y$ be a regular immersion in $\Sch^G/B$. Then the normal cone $C_{Y_0/Y}$ is the normal bundle $N_\iota:=\V(\sN_\iota)$ and one defines the Gysin map
\[
\iota^!\:\sE^\BM(Y, v)\to \sE^\BM(Y_0, v+L_\iota)
\]
as the composition
\[
\sE^\BM(Y, v)\xrightarrow{\sp^*_{\iota,v}}\sE^\BM(N_\iota, v)
\xrightarrow{0^!}
\sE^\BM(Y_0, v-\sN_\iota)=\sE^\BM(Y_0, v+L_\iota).
\]

More generally, suppose we are given a cartesian square
\[
\xymatrix{
X_0\ar[r]^{\iota'}\ar[d]_q&X\ar[d]^p\\
Y_0\ar[r]^\iota&Y
}
\]
in $\Sch^G/B$, with $\iota$ a regular immersion. Letting $N_\iota\to Y_0$ be the normal bundle, this diagram defines a closed immersion $\alpha\:C_{X_0/X}\to q^*N_\iota$. For $v\in \sK^G(X)$, we have the proper push-forward map
\[
\alpha_*\:\sE^\BM(C_{X_0/X}, v)\to \sE^\BM(q^*N_\iota, v)
\]
and the Gysin pull-back by the zero-section
\[
0_{q^*N_\iota}^!\:\sE^\BM(q^*N_\iota, v)\to \sE^\BM(X_0, v-\sN_\iota)=\sE^\BM(X_0, v+L_\iota).
\]
 
Define
\[
\iota_p^!\:\sE^\BM(X,v)\to \sE^\BM(X_0, v+L_\iota)
\]
to be the composition
\[
\sE^\BM(X,v)\xrightarrow{\sp^*_{\iota'}}\sE^\BM(C_{X_0/X},v)\\\xrightarrow{\alpha_*}
\sE^\BM(q^*N_\iota,v)\xrightarrow{0_{q^*N_\iota}^!}\sE^\BM(X_0, v+L_\iota).
\]
Note that $\iota_{\id_Y}^!$ is the Gysin map $\iota^!$. 

\subsection{Lci pullback} Recall that a morphism $f:Y\to X$ in $\Sch^G/B$ is an {\em lci morphism} if we can factor $f$ as $p\circ i$ where $i$ is a regular immersion and $p$ is a smooth morphism, both in $\Sch^G/B$\footnote{Since objects in $\Sch^G/B$ are quasi-projective over $B$, an lci morphism is a {\em smoothable} lci morphism in the sense of \cite{DJK}}

Given an lci morphism $f:Y\to X$, factored as $Y\xrightarrow{i}P\xrightarrow{p}X$ as above, define
\[
f^!:\sE^\BM_{a,b}(X, v)\to \sE^\BM_{a,b}(Y, v+L_f)
\]
as the composition
\[
\sE^\BM_{a,b}(X, v)\xrightarrow{i^!}:\sE^\BM_{a,b}(P, v+L_p)\xrightarrow{p^!}
\sE^\BM_{a,b}(Y, v+i^*L_p+L_i)\cong \sE^\BM_{a,b}(Y, v+L_f)
\]
where the last isomorphism uses the distinguished triangle $i^*L_p\to L_f\to L_i\to i^*L_p[1]$.
The arguments of \cite[Theorem 3.3.2]{DJK} adapted to the $G$-equivariant setting show that $f^!$ is well-defined, independent of the choice of factorization.  

\begin{definition}[Fundamental class]\label{def:FundClass} 1. Let $f:Y\to X$ be an lci morphism in $\Sch^G/B$, giving the lci pullback map $f^!:(1_X)^\BM(X/X)\to (1_X)^\BM(Y/X, L_f)$. Define
\[
\eta_f:=f^!(\id_{1_X})\in (1_X)^\BM(Y/X, L_f)
\]
considered as a morphism $\Sigma^{L_f}1_Y=\Sigma^{L_f}f^*1_X\to f^!1_X$. For $\sE\in \SH^G(B)$ a commutative ring spectrum with unit map $u_\sE:1_B \to \sE$, and let $\eta_f^\sE\in \sE^\BM(Y/X, L_f)$ be defined as $f^!(p_X^*(u_\sE))\circ \eta^f$.  For the structure morphism $p_Y\to B\in \Sch^G/B$ of an lci scheme over $B$, we write $[Y]_\sE$ for $\eta^\sE_{p_Y}$, and call $[Y]_\sE$ the {\em fundamental class} of the lci scheme $Y$.\\[2pt]
2. Given a cartesian diagram
\begin{equation}\label{eqn:LciCartSq}
\xymatrix{
Y'\ar[r]\ar[d]&X'\ar[d]^q\\
Y\ar[r]^f&X
}
\end{equation}
with $f$ lci, factor $f$ as $p\circ i$ as above, giving the diagram with all squares cartesian
\[
\xymatrix{
Y'\ar[r]^{\tilde{i}}\ar[d]&P'\ar[r]^{\tilde{p}}\ar[d]^r&X'\ar[d]^q\\
Y\ar[r]^i&P\ar[r]^p&X
}
\]
and define the {\em refined lci pullback}
\[
f^!_q:\sE^\BM_{a,b}(X', v)\to \sE^\BM_{a,b}(Y', v+L_f)
\]
as the composition
\[
\sE^\BM_{a,b}(X', v)\xrightarrow{\tilde{p}^!}\sE^\BM_{a,b}(P', v+L_p)\xrightarrow{i_r^!}
\sE^\BM_{a,b}(P', v+L_p+L_i)\cong \sE^\BM_{a,b}(P', v+L_f).
\]
\end{definition}

\begin{remark}\label{rem:CompDJK} 1. The classes $\eta_f\in (1_X)^\BM(Y/X, L_f)$ forms the {\em system of fundamental classes} constructed in \cite{DJK}, extended to the $G$-equivariant setting.  In particular, for composable lci morphisms $f,g$ in $\Sch^G/B$, we have $f^!(\eta^\sE_g)=\eta^\sE_{gf}$. For $f:Y\to X$ an lci morphism of lci schemes in $\Sch^G/B$, we similarly have 
\begin{equation}\label{eqn:LciPullbackFundClass}
f^!([X]_\sE)=f^!p_X^!([B]_\sE)=p_Y^!([B]_\sE)=[Y]_\sE
\end{equation}
in $\sE^\BM(Y, L_{Y/B})$.\\[2pt]
2. The lci pullback $f^!:\sE^\BM(X, v)\to \sE^\BM(Y, v+L_f)$ for an lci morphism $f:Y\to X$ in $\Sch^G/B$ is  
 the map sending $\alpha:\Sigma^v1_X\to p_X^!\sE$ to $f^!(\alpha)\circ \eta_f$, where  $f^!(-)$ is the functor $f^!:\SH(X)\to \SH(Y)$ from the six-functor formalism. \\[2pt]
3. Using the above notation, the refined lci pullback $f_q^!$  agrees with the refined Gysin map 
of \cite[Definition 4.2.5]{DJK} for the cartesian diagram \eqref{eqn:LciCartSq} and  extended to the $G$-equivariant setting.\\[2pt]
4. For $f$ a regular immersion, we often refer to $f^!$ as the Gysin pullback, and $f_q^!$ as the refined Gysin pullback.\\[2pt]
5. For $f$ smooth, $\eta_f:\Sigma^{\Omega_{Y/X}}1_Y\to f^!1_X$ is just the purity isomorphism \eqref{eqn:PurityIso} evaluated at $1_X$. 
\end{remark}

\begin{proposition}\label{prop:LciProperties} Take $Y\in \Sch^G/B$ with $v\in \sK^G(Y)$.\\[2pt]
 1 (compatibility with proper push-forward and smooth pull-back). Let
\[
\xymatrix{
W_0\ar[r]^-{f''}\ar[d]_{q_2}\ar@/_20pt/[dd]_q&W\ar[d]^{p_2}\ar@/^20pt/[dd]^p\\
X_0\ar[r]^-{f'}\ar[d]_{q_1}&X\ar[d]^{p_1}\\
Y_0\ar[r]^f&Y
}
\]
be a commutative diagram in $\Sch^G/B$, with both squares cartesian, with $f$ an lci morphism.\\[2pt]
1a. Suppose $p_2$ is proper. Then $q_2$ is also proper and the diagram
\[
\xymatrix{
\sE^\BM(W,v)\ar[r]^-{f^!_p}\ar[d]^{p_{2*}}& \sE^\BM(W_0, v+L_f)\ar[d]^{q_{2*}}\\
\sE^\BM(X,v)\ar[r]^-{f^!_{p_1}}& \sE^\BM(X_0, v+L_f)
}
\]
commutes.\\[2pt]
1b. Suppose $p_2$ is smooth. Then $q_2$ is also smooth and the diagram
\[
\xymatrix{
\sE^\BM(W,v+L_{p_2})\ar[r]^-{f^!_p}& \sE^\BM(W_0, v+L_f+L_{q_2})\\
\sE^\BM(X,v)\ar[u]^{p_{2}^!}\ar[r]^-{f^!_{p_1}}& \sE^\BM(X_0, v+L_f)\ar[u]^{q^!_{2}}
}
\]
commutes. \\[2pt]
2 (functoriality).  Let
\[
\xymatrix{
X_1\ar[r]^{f_1'}\ar[d]_{p_1}&X_0\ar[r]^{f_0'}\ar[d]_{p_0}&X\ar[d]^{p}\\
Y_1\ar[r]^{f_1}&Y_0\ar[r]^{f_0}&Y
}
\]
be a commutative diagram in $\Sch^G/B$, with all squares cartesian, and with $f_0$ and $f_1$ lci. Then
\[
f_1^!\circ f_0^!=(f_0\circ f_1)^!\:
\sE^\BM(X,v)\to \sE^\BM(X_1, v+L_{f_0\circ f_1}).
\]
Here we use the distinguished triangle
\[
f_1^*L_{f_0}\to L_{f_0\circ f_1}\to L_{f_1}\to f_1^*L_{f_0}[1]
\]
to identify $\sE^\BM(X_1, v+L_{f_1}+L_{f_0})$ with $\sE^\BM(X_1, v+L_{f_0\circ f_1})$.\\[2pt]
3 (excess intersection formula).  Let
\[
\xymatrix{
W_0\ar[r]^{\iota_2}\ar[d]_{q_2}&P_W\ar[d]^{p_2}\ar[r]^{g_2}&W\ar[d]^{r_2}\\
X_0\ar[r]^{\iota_1}\ar[d]_{q_1}&P_X\ar[d]^{p_1}\ar[r]^{g_1}&X\ar[d]^{r_1}\\
Y_0\ar[r]^\iota&P\ar[r]^g&Y
}
\]
be a commutative diagram in $\Sch^G/B$, with all squares cartesian, and with $\iota$ and $\iota'$ both regular immersions and $g$ smooth; let $f=g\circ\iota$, $f_1=g_1\circ \iota_1$.  We have a natural surjection  of locally free sheaves $q_1^*\sN_\iota\to \sN_{\iota'}$ with kernel $\sV$. Letting $V:=\V(\sV)$, this gives the Euler class $e^\sE(q^*V)\in \sE(W_0, q^*\sV)$ (Definition~\ref{def:EulerClass}).  We have the cap product  (\S\ref{subsec:Products})
\[
\sE(W_0, q^*\sV)\times  \sE^\BM(W_0, v+L_{\iota_1})\to
\sE^\BM(W_0, v+L_{\iota_1}-q^*\sV).
\]
Then for $\alpha\in \sE^\BM(W, v)$, we have
\[
f^!_{p_1p_2}(\alpha)=e^\sE(q^*V)\cap (f_1)^!_{p_2}(\alpha) \in \sE^\BM(W_0, v+L_\iota)=
\sE^\BM(W_0, v+L_{\iota_1}-q^*\sV),
\]
where we use the exact sequence
\[
0\to \sN_{\iota_1}\to q_1^*\sN_\iota\to \sV\to 0
\]
to give the canonical isomorphism 
\begin{multline*}
\sE^\BM(W_0, v+L_\iota):=\sE^\BM(W_0, v-\sN_\iota)\\\cong 
\sE^\BM(W_0, v-(\sN_{\iota_1}+q^*\sV))=:\sE^\BM(W_0, v+L_{\iota_1}-q^*\sV).
\end{multline*}
In addition, suppose the bottom square is Tor-independent. Then $\sV=0$ and 
\[
f^!_{p_1p_2} = (f_1)^!_{p_2}.
\]
\end{proposition}

\begin{proof} For (1) and (2), see \cite[Lemma 2.3.14, Proposition 2.4.2, Theorem 4.2.1, Proposition 4.2.2, Proposition 4.2.6]{DJK}; this is in the non-equivariant setting, but the arguments generalize without difficulty. For (3), we have $f^!_{p_1p_2}=\iota^!_{p_1p_2}\circ g_2^!$ and 
 $(f_1)^!_{p_2}=(\iota_1)^!_{p_2}\circ g_2^!$, so we reduce to the case of a regular immersion, which follows from \cite[Proposition 3.2.8]{DJK}. \end{proof}

\subsection{Refined intersection product}\label{subsec:RefinedProduct} We use the external product \eqref{eqn:ExternalProd} and refined lci pullback to define the refined intersection product.

Suppose we have a $Y_0\in \Sm^G/B$ and morphisms $f\:Y\to Y_0$ in $\Sm^G/B$, $g\:X\to Y_0$ in $\Sch^G/B$. Let $\delta_{Y_0}\:Y_0\to Y_0\times_BY_0$ be the diagonal. This gives us the cartesian diagram
\begin{equation}\label{eqn:RefinedIntersectionDiag}
\xymatrix{
Y\times_{Y_0}X\ar[r]^{\delta'}\ar[d]^p&Y\times_BX\ar[d]^{f\times g}\\
Y_0\ar[r]^{\delta_{Y_0}}&Y_0\times_BY_0.
}
\end{equation}
The refined intersection product
\[
(-)\cdot_{f,g}(-)\:\sE^\BM_{a,b}(Y/B, v)\times \sE^\BM_{c,d}(X/B, w)\to
\sE^\BM_{a+c,b+d}(Y\times_{Y_0}X/B, v+w-\Omega_{Y_0/B})
\]
is defined as the composition
\begin{multline*}
\sE^\BM_{a,b}(Y/B, v)\times \sE^\BM_{c,d}(X/B, w)\xrightarrow{\boxtimes_{Y,X}}
\sE^\BM_{a+c,b+d}(Y\times_BX/B, v+w)\\
\xrightarrow{\delta_{Y_0}^!}\sE^\BM_{a+c,b+d}(Y\times_{Y_0}X/B, v+w-\Omega_{Y_0/B}).
\end{multline*}

\begin{lemma}\label{lem:ProductComp} 
1.  Given $f:Z\to X$ in $\Sch^G/B$, $\alpha\in \sE^{**}(X,v)$, $\beta\in \sE^{**}(X,w)$,  we have
\begin{equation}\label{eqn:ProductCohoFunct}
f^*(\alpha \cup_X\beta)=f^*(\alpha)\cup_Zf^*(\beta)
\end{equation}
If $f$ is proper, then for $\alpha\in \sE^{**}(X,v)$, $\gamma\in \sE^\BM_{**}(Z,w)$
\begin{equation}\label{eqn:ProjFormula}
f_*(f^*(\alpha)\cap_Z \gamma)=\alpha\cap_X f_*(\gamma)
\end{equation}
If $f$ is lci, then for $\alpha\in \sE^{**}(X,v)$, $\beta\in \sE^\BM_{**}(X,w)$
\begin{equation}\label{eqn:ProductFunct}
f^!(\alpha \cap_X\beta)=f^*(\alpha)\cap_Zf^!(\beta)
\end{equation}
2. Let $i_{Y'}\:Y'\to Y$ be a closed immersion in $\Sm^G/B$ and let $i_{X'}\:X'\to X$ be a closed immersion in $\Sch^G/B$.  For $\alpha\in 
\sE^\BM_{a,b}(Y', v)$, $\beta\in  \sE^\BM_{c,d}(X', w)$,  we have 
\[
i_{Y'*}(\alpha)\boxtimes_{Y,X} i_{X'*}(\beta)=(i_{Y'}\times i_{X'})_*(\alpha\boxtimes_{Y', X'}\beta).
\]
3. Suppose  that  $i_{X'}$ is a regular immersion and that the cartesian squares
\[
\xymatrix{
Y\times_BX'\ar[r]\ar[d]&X'\ar[d]\\
Y\times_BX\ar[r]&X
}
\]
and
\[
\xymatrix{
Y'\times_BX'\ar[r]\ar[d]&Y'\ar[d]\\
Y\times_BX'\ar[r]&Y
}
\]
are both Tor-independent (for example if $Y\to B$ and $X'\to B$ are both flat). Then $i_{Y'}\times i_{X'}:Y'\times_BX'\to Y\times_B X$ is a regular immersion and for $\alpha\in 
\sE^\BM_{a,b}(Y, v)$, $\beta\in  \sE^\BM_{c,d}(X, w)$,  we have 
\[
i_{Y'}^!(\alpha)\boxtimes_{Y',X'} i_{X'}^!(\beta)=(i_{Y'}\times i_{X'})^!(\alpha\boxtimes_{X,Y}\beta).
\]
4. (a) Given $Y_0\in\Sm^G/B$ and morphisms $f:Y\to Y_0$, $g:X\to Y_0$ in $\Sch^G/B$, we have
\[
i_{Y'*}(\alpha)\cdot_{f,g}i_{X'*}(\beta)=i_{Y'\times X'*}(\alpha\cdot_{fi_{Y'},gi_{X'}}\beta).
\]
(b)  If the diagram \eqref{eqn:RefinedIntersectionDiag} is Tor-independent, then $\delta'$ is a regular immersion and 
\[
\alpha\cdot_{f,g}\beta=\delta^{\prime!}(\alpha\boxtimes_{Y',X'} \beta).
\]
\end{lemma}

\begin{proof} (2) and (3) follow from (1) and Proposition~\ref{prop:LciProperties}.

For (1), the formula \eqref{eqn:ProductCohoFunct} is just the monoidality \eqref{align:Product}.  
The projection formula \eqref{eqn:ProjFormula}  follows from \eqref{align:ProjFormula}, Propostion~\ref{prop:PushPullComp} and Proposition~\ref{prop:LciProperties}. 

For the formula \eqref{eqn:ProductFunct}, we take $a=b=0$ to simplify the notation. Given $\alpha:1_X\to\Sigma^vp_X^*\sE$, $\beta:1_X\to \Sigma^{-w}p_X^!\sE$, $\alpha\cap_X\beta$ is the composition
\[
1_X=1_X\otimes_X1_X\xrightarrow{\alpha\otimes\beta}\Sigma^{v-w}p_X^*\sE\otimes p_X^!\sE\to \Sigma^{v-w}p_X^!(\sE\otimes\sE)\to \Sigma^{v-w}p_X^!\sE.
\]
Applying $f^!$ gives (using Remark~\ref{rem:CompDJK} (2))
\[
\Sigma^{L_f}1_Y=\Sigma^{L_f}f^*1_X\xrightarrow{\eta_f}f^!(1_X) \xrightarrow{f^!(\alpha\cup_X\beta} \Sigma^{v-w}f^!p_X^!\sE
\]
representing $f^!(\alpha\cap_X\beta)\in \sE^\BM(Y, w-v+L_f)$, 
and \eqref{eqn:ProductFunct} follows from the commutativity of the following diagram.
\[
\xymatrix{
f^*1_X\otimes f^!1_X\ar[r]\ar[d]^{f^*\alpha\otimes f^!\beta}&f^!(1_X\otimes 1_X)\ar[d]^{f^!(\alpha\otimes\beta)}\\
f^*p_X^*\sE\otimes f^!p_X^!\sE\ar@{=}[dd]\ar[r]&f^!(p_X^*\sE\otimes p_X^!\sE)\ar[d]\\
&f^!p_X^!(\sE\otimes\sE)\ar@{=}[d]\\
(p_Xf)^*\sE\otimes(p_Xf)^!\sE\ar[r]&(p_Xf)^!(\sE\otimes\sE)
}
\]
where the unmarked  arrows are instances of \eqref{eqn:ExceptPullback}.

See also \cite[Proposition 2.2.12, Theorem 4.2.1]{DJK} 

For (4),  (a) follows from Proposition~\ref{lem:ProductComp}(2)  and Proposition~\ref{prop:LciProperties}(1a).  (b) follows (a) and the excess intersection formula, noting that $p^*\sN_{\delta_{Y_0}}\to \sN_{\delta'}$ is an isomorphism, so the Euler class in Proposition~\ref{prop:LciProperties}(3) is 1. 
 \end{proof}

\subsection{Relative refined lci pullback}\label{subsec:RelPull-back}
Let $p:B_0\to B_1$ be a morphism in $\Sm^G/B$; since $B_0$ and $B_1$ are smooth over $B$, $p$ is an lci morphism. We define the {\em relative purity isomorphism}
\begin{equation}\label{eqn:RelPurity}
\phi_p:p_{B_0\#}\xrightarrow{\sim} p_{B_1\#}p_!\Sigma^{L_p}
\end{equation}
as the composition
\[
p_{B_0\#}\cong p_{B_0!}\Sigma^{\Omega_{p_{B_0}}}
\cong p_{B_1!}p_!\Sigma^{p^*\Omega_{p_{B_1}}+L_p}
= p_{B_1!}\Sigma^{\Omega_{p_{B_1}}}p_!\Sigma^{L_p}
\cong p_{B_1\#}p_!\Sigma^{L_p},
\]
where we use the purity isomorphisms $p_{B_i!}\cong p_{B_i\#}\circ\Sigma^{-\Omega_{p_{B_i}}}$, $i=0,1$, and the isomorphism $\Sigma^{p^*\Omega_{p_{B_1}}+L_p}\cong
\Sigma^{\Omega_{p_{B_0}}}$ arising from the distinguished triangle $p^*\Omega_{p_{B_1}}\to \Omega_{p_{B_0}}\to L_p$.

Next, consider a commutative diagram in $\Sch^G/B$
\begin{equation}\label{eqn:RelPullbackSq}
\xymatrix{
X_0\ar[d]^{f}\ar[dr]^g\\
B_0\ar[r]^p&B_1
}
\end{equation}
with $p:B_0\to B_1$ in $\Sm^G/B$ as above. For $\sE\in \SH^G(B)$, $v\in \sK(X_0)$, \eqref{eqn:RelPurity} induces the natural isomorphism
\[
\phi_p^*:\sE^\BM(X_0/B_1, v+L_p)\xrightarrow{\sim} \sE^\BM(X_0/B_0, v)
\]
as the composition
\begin{align*}
\sE^\BM(X_0/B_1, v+L_p)&=\Hom_{\SH(B_1)}(g_!\Sigma^{v+L_p}1_{X_0}, p_{B_1}^*\sE)\\
&=\Hom_{\SH(B)}(p_{B_1\#}p_!\Sigma^{L_p}f_!\Sigma^{v}1_{X_0}, \sE)\\
&\xymatrix{\ar[r]^{\phi_p^*}_\sim&}
\Hom_{\SH(B)}(p_{B_0\#}f_!\Sigma^{v}1_{X_0}, \sE)\\
&=\Hom_{\SH(B_0)}(f_!\Sigma^{v}1_{X_0}, p_{B_0}^*\sE)=\sE^\BM(X_0/B_0, v).
\end{align*}

Let
\begin{equation}\label{eqn:RelPullbackSq}
\xymatrix{
X_0\ar[r]^q\ar[d]^{f_0}&X_1\ar[d]^{f_1}\\
B_0\ar[r]^p&B_1
}
\end{equation}
be a cartesian square in $\Sch^G/B$, with $B_0$ and $B_1$ smooth over $B$. Take $v\in \sK^G(X)$ and $\sE\in \SH^G(B)$.  

\begin{definition}\label{def:RelLciPullback}
We define the {\em relative refined  pull-back}
\[
p^!_{f_1/B_0}\:\sE^\BM(X_1/B_1, v)\to \sE^\BM(X_0/B_0,v)
\]
as the composition
\[
\sE^\BM(X_1/B_1, v)\xrightarrow{p^!_{f_1}}  \sE^\BM(X_0/B_1,v+L_p)
\xymatrix{\ar[r]^{\phi_p^*}_\sim&} \sE^\BM(X_0/B_0,v).
\]
\end{definition}

\begin{remark}\label{rem:RelPullback} Since the relative refined lci pullback is just the refined lci pullback composed with the relative purity isomorphism $\phi_p^*$, all the properties of refined lci pullback, as described in Proposition~\ref{prop:LciProperties} and elsewhere, hold in the corresponding fashion for the relative refined lci pullback, when the latter is defined.

If the morphism $p$ is smooth, then the additional properties that hold for smooth pullback, such as commuting with specialization in cartesian squares (Proposition~\ref{prop:SpecSmoothCommute})  or the compatibility with localization sequences (Proposition~\ref{prop:LocSmoothComp}) also hold in the relative setting. 

If the context makes the meaning clear, we will sometimes write $p^!_{f_1}$, or even just $p^!$,  for $p^!_{f_1/B_0}$.
\end{remark}

\subsection{Fundamental classes, cone classes and Euler classes}\label{sec:FundamentalClasses}  Let $\sE\in\SH(B)$ be a commutative ring spetrum  Recall (Definition~\ref{def:FundClass}) the fundamental class $[X]_\sE\in \sE^\BM(X/B, L_{p_X})$  of an lci scheme $p_X:X\to B\in \Sch^G/B$.  Since $\sE$ is fixed in this section, we will drop the subscript ${}_\sE$ and write $[X]$ for $[X]_\sE$.

For $p_X\:X\to B$ lci, and $v\in \sK^G(X)$,  cap product with $[X]$ gives the map
\[
-\cap_X[X]\:\sE^{a,b}(X, v)\to \sE^\BM_{-a,-b}(X, L_{X/B}-v);
\]
in case $X$ is smooth over $B$, this is the same as the   Poincar\'e duality isomorphism \eqref{eqn:PDIso}. 

\begin{lemma}\label{lem:PDpull-back} Let $f\:Y\to X$ be a morphism in $\Sm/B$ ($f$ is automatically an lci morphism). Then the diagram
\[
\xymatrix{
\sE^{a,b}(X, v)\ar[r]_-\sim^-{-\cap_X[X]}\ar[d]^{f^*}& \sE^\BM_{-a,-b}(X, L_{X/B}-v)\ar[d]^{f^!}\\
\sE^{a,b}(Y,f^*v)\ar[r]_-\sim^-{-\cap_Y[Y]}& \sE^\BM_{-a,-b}(Y, L_{Y/B}-v)
}
\]
commutes.
\end{lemma}
This follows from \eqref{eqn:ProductFunct}.

\begin{definition}[Normal cone class]\label{def:FundClassCone}
Let $\iota\:X\to Y$ be a  closed immersion in $\Sch^G/B$, with $Y$ smooth over $B$. We have the cone $C_{X/Y}:=\Spec \oplus_{n\ge0}\sI_X^n/\sI_X^{n+1}$. Define the fundamental class
$[C_{X/Y}]$
by 
\[
[C_{X/Y}]=\sp^*_{\iota,\Omega_{Y/B}}([Y])\in \sE^\BM(C_{X/Y},\Omega_{Y/B}).
\]
If we need to keep track of the base-scheme $B$, we write this as $[C_{X/Y}/B]$.
\end{definition}

\begin{remark} Let  $\iota\:X\to Y$ be a regular closed immersion in $\Sch/B$ with $Y$ smooth over $B$. We have the conormal sheaf $\sN_\iota:=\sI_X/\sI_X^2$ and corresponding normal bundle $N_\iota:=\V(\sN_\iota)$. Then   $C_{X/Y}=N_\iota$. Since $X$ is lci over $B$,  so is $N_\iota$ and $[C_{X/Y}]=[N_\iota]$. Indeed, 
\[
[X]=\iota^!([Y])=0^!(\sp_{\iota}^*([Y]))=0^!([C_{X/Y}]),
\]
\[
0^!([N_\iota])=0^!(p_{N_\iota}^!([B]))=p_X^!([B])=[X]=0^!([C_{X/Y}])
\]
and $0^!\:\sE^\BM(N_\iota, \Omega_{N_\iota/B})\to \sE^\BM(X, \Omega_{X/B})$ is an isomorphism.
\end{remark}

\begin{remark}
It was shown in \cite[Theorem 3.2]{LevineVirt} that, given two $G$-equivariant closed immersions $\iota_j\:X\to Y_j$, $j=1,2$, $Y_j\in \Sm^G/B$, there is a canonical isomorphism
\[
\psi_{i_2, i_1}\: \sE^\BM_{**}(C_{X/Y_1}, \Omega_{Y_1/B})\xrightarrow{\sim}
 \sE^\BM_{**}(C_{X/Y_2}, \Omega_{Y_2/B})
 \]
satisfying $\psi_{i_2, i_1}[C_{X/Y_1}]=[C_{X/Y_2}]$ and  $\psi_{i_3, i_1}\circ \psi_{i_2, i_1}=\psi_{i_3, i_1}$. This says that $[C_{X/Y}]\in \sE^\BM(C_{X/Y}, \Omega_{Y/B})$ is canonically defined, independent of the choice of immersion $\iota$.
\end{remark}

We conclude with the definition of the Euler class in $\sE$-Borel-Moore homology for a $G$-vector bundle $V\to X$ on an lci scheme $X$. 

\begin{definition}[Borel-Moore Euler class]\label{def:BMEulerClass} Let $\sV$ be a locally free $G$-linearized coherent sheaf on $X\in \Sch^G(B)$, giving the $G$-equivariant vector bundle $V:=\V(\sV)\to X$. Let $\sE\in \SH^G(B)$ be a commutative ring spectrum with unit map $u_\sE:1_B\to \sE$. Suppose that $X$ is an lci scheme over $B$. Then we have the  fundamental class $[X]_\sE\in \sE^\BM(X, L_{X/B})$, and the zero-section $s_V:X\to V$ is a regular immersion. We define the {\em Borel-Moore Euler class} $e^\BM_{\sE}(V)$ by
\[
e^\BM_{\sE}(V):=s_V^!s_{V*}([X]_\sE)\in \sE^\BM(X/B, L_{X/B}-\sV).
 \]
\end{definition}

\begin{remark}\label{rem:BMEulerNat} For
 $X$ lci, the Borel-Moore Euler class is related to the cohomological Euler class by 
 \[
 e^\BM_{\sE}(V)=e^\sE(V)\cap_X [X]_\sE,
 \]
and the naturality \eqref{eqn:EulerClassNat} together with \eqref{eqn:ProductFunct} and \eqref{eqn:LciPullbackFundClass} gives a similar identity for $f:Y\to X$ an lci morphism of lci schemes over $B$,
 \begin{equation}\label{eqn:BMEulerClassNat}
f^!(e^\BM_{\sE}(V))=e^\BM_{\sE}(f^*V)\in \sE^\BM(Y,  L_{Y/B}-f^*\sV).
\end{equation}
\end{remark}

 \subsection{Oriented and $\SL$-oriented spectra} We briefly leave the $G$-equivariant setting. 
 Let $\sE$ be a commutative ring spectrum in $\SH(B)$. There is a well-establised theory of oriented and $\SL$-oriented spectra \cite{Anan21, Anan19, Panin}; for us, the relevant information is the following.
 
An orientation for   $\sE$ assigns, for each   $p_X:X\to B\in \Sch/B$ and each $v\in \sK(X)$,  of   {\em Thom isomorphisms} in $\SH(X)$,
\[
\th_{v,*}:\Sigma^{v-\rnk(v)\sO_X}p_X^*\sE\xrightarrow{\sim}p_X^*\sE,\ \th_{v,!}:\Sigma^{v-\rnk(v)\sO_X}p_X^!\sE\xrightarrow{\sim}p_X^!\sE,
\]
natural in $v$ and $X$; if $p_X$ is smooth then $\th_{v,*}$ and $\th_{v,!}$ are induced from one another by the purity isomorphism $p_X^!\cong \Sigma^{\Omega_{X/B}}p_X^*$.  In addition, for each distinguished triangle $v'\to v\to v''\to v[1]$ in $D^\perf(X)$,   the following diagram commutes
\[
\xymatrixcolsep{20pt}
\xymatrix{
\Sigma^{v-\rnk(\sV)\sO_X}p_X^*\sE\ar[rrrr]^{\th_{v,*}}_{\sim}\ar[d]^\wr&&&&p_X^*\sE\\
\Sigma^{v'-\rnk(v')\sO_X}(\Sigma^{v''-\rnk(v'')\sO_X}p_X^*\sE)\ar[rrr]^-{\Sigma^{v'-\rnk(v')\sO_X}\th_{v'',*}}_-{\sim}&&&\Sigma^{v'-\rnk(v')\sO_X}p_X^*\sE\ar[r]^-{\th_{v',*}}_-{\sim}&p_X^*\sE\ar@{=}[u]
}
\]
where the vertical isomorphism arises from the natural isomorphism $\Sigma^v\cong
\Sigma^{v'}\circ\Sigma^{v''}$; we have a corresponding commutative diagram for the Thom isomorphisms $\th_{-,!}$.

An $\SL$-orientation for  $\sE$ assigns,  for  each $p_X:X\to B\in \Sch/B$ and each $v\in \sK(X)$, of  Thom isomorphisms in $\SH(X)$,
\[
\th_{v,*}:\Sigma^{v-\rnk(v)\sO_X}p_X^*\sE\xrightarrow{\sim}\Sigma^{\det v-\sO_X}p_X^*\sE, 
\]
\[
\th_{v,!}:\Sigma^{v-\rnk(v)\sO_X}p_X^!\sE\xrightarrow{\sim}\Sigma^{\det v-\sO_X}p_X^!\sE, 
\]
natural in $v$ and $X$, with a similar compatibility with respect to  distinguished triangles in $D^\perf(X)$, and similarly induced from one another by the purity isomorphism if $p_X$ is smooth.

In addition,   in the $\SL$-oriented case, we have natural canonical isomorphisms
\[
\Sigma^{\sL\otimes\sM^{\otimes 2}}p_X^*\sE\xymatrix{\ar[r]^{\can}_\sim&}
\Sigma^{\sL}p_X^*\sE,\ \Sigma^{\sL\otimes\sM^{\otimes 2}}p_X^!\sE\xymatrix{\ar[r]^{\can}_\sim&}
\Sigma^{\sL}p_X^!\sE
\]
for each pair of invertible sheaves $\sL, \sM$ on $X\in \Sch/B$; these isomorphisms satisfy the evident associativity and again are induced from one another by purity if $p_X$ is smooth.

In both cases, these isomorphism arise from a theory of Thom classes;  we refer the reader to  \cite[\S 3.3]{LevineAtiyahBott} for a detailed overview of  Thom classes and Thom isomorphisms for oriented $\sE$ and $\SL$-oriented $\sE$. 

For  $X\in \Sch/B$ and $v\in \sK(X)$, the Thom isomorphisms $\th_{v,!}$ induce the Thom isomorphism
\[
\sE^\BM_{a,b}(X, v-\rnk(v)\sO_X)\xymatrix{\ar[r]^{\th_{v,!}}_\sim&} \sE^\BM_{a,b}(X, \det(v)-\sO_X)
\]
\begin{definition}\label{def:(Rel)Orientation}
An {\em orientation} for $v\in \sK(X)$ is a pair $(\sL,\rho)$ with $\sL$ an invertible sheaf on $X$ and $\rho$ an isomorphism $\rho:\det(v)\xrightarrow{\sim}\sL^{\otimes 2}$. Given an orientation for $v$, we thus have the isomorphism
\begin{multline*}
\sE^\BM_{a,b}(X, v-\rnk(v)\sO_X)\xymatrix{\ar[r]^{\th_{v,!}}_\sim&} \sE^\BM_{a,b}(X, \det(v)-\sO_X)\\
\xrightarrow{\rho} \sE^\BM_{a,b}(X, \sL^{\otimes 2}-\sO_X)
\xymatrix{\ar[r]^{\can}_\sim&} \sE^\BM_{a,b}(X).
\end{multline*}

Let $V\to X$ is a vector bundle of rank $d$ with sheaf of sections $\sV$ on an lci scheme $X\in \Sch^G/B$. A {\em relative orientation} for $V$ is a pair $(\sL,\rho)$, with $\sL$ an invertible sheaf on $X$ and $\rho$ an isomorphism $\rho:\det L_{X/B}\otimes \det(\sV)^\vee\xrightarrow{\sim}\sL^{\otimes 2}$, that is, an orientation for $L_{X/B}+\sV[1]$.
\end{definition}

 \begin{remark} \label{rem:EulerClassDeg}  Suppose $\sE\in \SH(B)$ is $\SL$-oriented, $X\in \Sch/B$ is lci and proper over $B$  of pure dimension $d$ over $B$ and $V\to X$ is a vector bundle of rank $d$ with sheaf of sections $\sV$, and $(\sL,\rho)$ a relative orientation for $V$.  Then we have the Borel-Moore Euler class 
\begin{multline*}
e^\BM_{\sE}(V)\in \sE^\BM(X/B, L_{X/B}-\sV)\\\cong 
\sE^\BM(X/B, \det L_{X/B}\otimes\det(\sV)^\vee){\xymatrix{\ar[r]^\rho_\sim&}}\sE^\BM(X/B),
\end{multline*}
where the first isomorphism arises from the $\SL$-orientation and the second from the relative orientation $\rho$.  We may then apply the degree map
\[
\deg_B^\sE:=p_{X*}:\sE^\BM(X/B)\to \sE^\BM(B/B)=\sE^{0,0}(B),
\]
giving the element $\deg_B^\sE(e^\BM_{\sE}(V))\in\sE^{0,0}(B)$.
\end{remark}

\section{Vistoli's lemma}\label{sec:Vistoli} 

In this section we state and prove our  version of Vistoli's lemma in the $G$-equivariant motivic setting, for tame $G$ over a base-scheme $B$, quasi-projective over a noetherian ring of finite Krull dimension. We fix a commutative monoid $\sE\in\SH^G(B)$ and use the $\sE$-valued twisted Borel-Moore homology $\sE^\BM_{*,*}(-/B,-)$. To simplify the notation, we sometimes state results for $\sE^\BM(-/B,-)$, but the results are valid for $\sE^\BM_{*,*}(-/B,-)$  by applying a suitable twist. As we fix the base-scheme $B$ for the entire section, we will drop the $-/B$ from the notation for Borel-Moore homology, unless we are using additional base-schemes.

Consider a cartesian square in $\Sch^G/B$
\[
\xymatrix{
Z\ar[r]\ar[d]&X_1\ar[d]^{i_1}\\
X_2\ar[r]^{i_2}&Y
}
\]
with $i_1, i_2$ closed immersions. This gives rise to the commutative diagram with all squares cartesian
\begin{equation}\label{eqn:DoubleConeDiagram}
\xymatrix{
&C_{\alpha_1}\ar[d]_{\beta_1}\ar[dr]\\
C_{\alpha_2}\ar[dr]\ar[r]^-{\beta_2}&q^*N_f\ar[r]\ar[d]
&i_2^{*}C_{i_1}\ar[r]^{\alpha_1}\ar[d]&C_{i_1}\ar[d]\\
&i_1^{*}C_{i_2}\ar[r]\ar[d]^{\alpha_2}&Z\ar[r]^{i_2'}\ar[d]^{i_1'}&X_1\ar[d]^{i_1}\\
&C_{i_2}\ar[r]&X_2\ar[r]^{i_2}&Y.
}
\end{equation}
To explain the notation: given a pair of maps $f\:S\to T$, $U\to T$ in $\Sch^G/B$ we write $f^*U$ for the fiber product $S\times_TU$. If $f$ is a closed immersion, we write $C_f\to T$ or $C_{S/T}\to T$ for the normal cone of $f\:S\to T$.  

The maps $\beta_1, \beta_2$ are closed immersions defined as follows. If $X_j\subset Y$ is defined by an ideal sheaf $\sI_j$, then $C_{i_j}=\Spec_{\sO_Y} \oplus_n\sI_j^n/\sI_j^{n+1}$. For $j\neq j'$, $i^{*}_{j'}C_{i_j}$ is $\Spec_{\sO_Y}  \oplus_n (\sI_j^n/\sI_j^{n+1})\otimes_{\sO_Y} \sO_Y/I_{j'}$ and 
\[
C_{\alpha_j}=\Spec_{\sO_Y} \oplus_m \sI_{j'}^m\cdot(\oplus_n \sI_j^n/I_j^{n+1})/\sI_{j'}^{m+1}\cdot(\oplus_n \sI_j^n/\sI_j^{n+1}), 
\]
while 
\[
i_2^*C_{i_1}\times_Yi_1^*C_{i_2}=\Spec_{\sO_Y}\oplus_{n,m} \sI_{j'}^m/\sI_{j'}^{m+1} \otimes_{\sO_Y}\sI_j^n/\sI_j^{n+1}.
\]
The sheaf of $\sO_Y$-algebras $\oplus_m \sI_{j'}^m\cdot(\oplus_n \sI_j^n/I_j^{n+1})/\sI_{j'}^{m+1}\cdot(\oplus_n (\sI_j^n/\sI_j^{n+1})$ is in an evident way a quotient of  $\oplus_{n,m} \sI_{j'}^m/\sI_{j'}^{m+1} \otimes_{\sO_Y}\sI_j^n/\sI_j^{n+1}$, which defines the closed immersion $\beta_j$. 

Here is our version of Vistoli's lemma.
\begin{proposition}\label{prop:SpecializationComm}  For $v\in K^G(Y)$, we have
\[
\beta_{1*}\circ \sp^*_{\alpha_1,v}\circ\sp^*_{i_1,v}=
\beta_{2*}\circ \sp^*_{\alpha_2,v}\circ\sp^*_{i_2,v}
\]
as maps $\sE^\BM(Y,v)\to \sE^\BM(i_2^*C_{i_1}\times_Yi_1^*C_{i_2},v)$.
\end{proposition}

\begin{remark}\label{rem:CommentsAndHistory}
Vistoli \cite[Lemma 3.16]{Vistoli} proves this in the case of the Chow groups by constructing for $y\in \CH_*(Y)=H\Z^\BM_{2*,*}(Y)$ an explicit rational equivalence 
\[
\beta_{1*}\circ \sp^*_{\alpha_1}\circ\sp^*_{i_1}(y)\sim
\beta_{2*}\circ \sp^*_{\alpha_2}\circ\sp^*_{i_2}(y).
\]
on $i_2^*C_{i_1}\times_Yi_1^*C_{i_2}$. He uses this to extend his result to the case of Deligne-Mumford stacks.
Our proof is different from Vistoli's and relies on a geometric construction  \eqref{eqn:VistoliGeomConst}, combined with the formal properties of localization sequences for Borel-Moore homology; we do not consider the extension to stacks here, but we believe that our essentially formal approach should make this extension feasible.

Our proof relies on the {\em double deformation space} \eqref{eqn:DDDiag} and the associated {\em double deformation diagram}  \eqref{eqn:VistoliGeomConst}  and constructing an element in the appropriate Borel-Moore homology of an open subscheme of $i_2^*C_{i_1}\times_Yi_1^*C_{i_2}$, whose boundary in a suitable localization sequence gives the relation asserted in Proposition~\ref{prop:SpecializationComm}. A similar approach was used by Rost \cite[Lemmas 11.6 and 11.7, Theorem 13.1]{Rost} to establish a functoriality of pull-back maps for his cycle complexes. The diagram \cite[(10.5)]{Rost} used by Rost is a bit different from ours, as it serves a somewhat different purpose: Rost was proving an ``associativity'' result with his diagram, whereas we are proving what amounts to a ``commutativity'' result (see Corollary~\ref{cor:commutativity} below). Rost also mentions as an aside a construction of another double deformation space \cite[(10.6)]{Rost} which appears to be the same, or at least closely related, to our construction, but he does not use this in his paper.  We thank the referee for pointing out to us these constructions and arguments of Rost. Rost's construction in  \cite[(10.6)]{Rost} was generalized to a multi-deformation diagram by Ivorra \cite[\S3.1.3]{Ivorra}.
\end{remark}

Before we give the proof of Proposition~\ref{prop:SpecializationComm}, we derive as a corollary of Proposition~\ref{prop:SpecializationComm} the commutativity of refined Gysin pull-back. 

\begin{corollary}[Commutativity]\label{cor:commutativity} Let 
\[
\xymatrix{
Z\ar[r]^{i_2'}\ar[d]^{i_1'}&X_1\ar[r]^{q'}\ar[d]^{i_1}&S\ar[d]^i\\
X_2\ar[r]^{i_2}\ar[d]_{p'}&Y\ar[d]_p\ar[r]^q&T\\
S'\ar[r]^j&T'
}
\]
be a commutative diagram in $\Sch^G/B$ with all squares cartesian. Suppose that $i$ and $j$ are regular immersions. Then the diagram
\[
\xymatrix{
\sE^\BM(Y, v)\ar[r]^-{i^!_q}\ar[d]^{j^!_p}&\sE^\BM(X_1, v+L_i)\ar[d]^{j^!_{pi_1}}\\
\sE^\BM(X_2, v+L_j)\ar[r]^-{i^!_{qi_2}}&\sE^\BM(Z, v+L_i+L_j)
}
\]
commutes.
\end{corollary}

\begin{proof} We  refer to the diagram \eqref{eqn:DoubleConeDiagram}, and let $t_2:Z\to S$, $t_1:Z\to T$ be the respective compositions, $t_2:=q'\circ i_2'$, $t_1:=p'\circ i_1'$, giving the commutative diagram
\[
\xymatrix{
Z\ar[r]^{i_2'}\ar[d]^{i_1'}\ar@/^20pt/[rr]^{t_2}\ar@/_20pt/[dd]_{t_1}&X_1\ar[r]^{q'}\ar[d]^{i_1}&S\ar[d]^i\\
X_2\ar[r]^{i_2}\ar[d]_{p'}&Y\ar[d]_p\ar[r]^q&T\\
S'\ar[r]^j&T'
}
\]

Let $\sI_i\subset \sO_T$ be the ideal sheaf of the closed immersion $i$, and let $N_i\to S$ be the normal bundle, $N_i=\Spec_{\sO_S}\Sym^*\sI_i/\sI_i^2$. Letting $\sI_{i_1}\subset \sO_S$ be the ideal sheaf of the closed immersion $i_1$, we have the surjection $\sO_S\otimes_{\sO_T}\sI_i\to \sI_{i_1}$, which induces the closed immersion
\[
\gamma_{i_1}\:C_{i_1}=\Spec_{\sO_X}\oplus_n \sI_{i_1}^n/\sI_{i_1}^{n+1}\hookrightarrow
\Spec_{\sO_X}\Sym^*\sI_i/\sI_i^2\otimes_{\sO_S}\sO_{X_1}=
q^{\prime*}N_i
\]
over $X_1$.  The closed immersion $i_2'\:Z\to X_1$ induces the closed immersions  $i_2^{\prime\prime}\:t_2^*N_i\to q^{\prime*}N_i$  and $\alpha_1:i_2^{\prime*}C_{i_1}\to C_{i_1}$; let   $C_{\alpha_1}$ be the cone of $\alpha_1$. This gives us the
 the commutative diagram
\begin{equation}\label{eqn:BigComDiag1}
\xymatrix{
C_{\alpha_1}\ar@{_(->}[d]^{\gamma_{\alpha_1}}\ar[r]&i_2^{\prime*}C_{i_1} \ar@/_30pt/[ddddr]_{h_1}\ar[r]^{\alpha_1}\ar[ddr]_{g_1}\ar@{^(->}[dr]^{\gamma'_{i_1}}&C_{i_1}\ar[ddr]|-\hole\ar@{^(->}[dr]^{\gamma_{i_1}}\\
h_1^*N_j\ar@<3pt>[ur]\ar@/_20pt/[dddr]&&t_2^*N_i\ar[d]\ar[r]_-{i_2''}&q^{\prime*}N_i\ar[d]\\
& &Z\ar@/_13pt/[dd]_{t_1}\ar@/_13pt/[rr]_(.28){t_2}\ar[r]^{i_2'}\ar[d]^{i_1'}&X_1\ar[r]^{q'}\ar[d]^{i_1}|(.32)\hole&S\ar[d]^i\\
&&X_2\ar[r]^{i_2}\ar[d]^{p'}&Y\ar[d]_p\ar[r]^q&T\\
&N_j\ar[r]&S'\ar[r]^j&T',
}
\end{equation}
with all squares cartesian; here the closed immersion $\gamma_{\alpha_1}$ is defined similarly to $\gamma_{i_1}$, using  the cartesian diagram
\[
\xymatrix{
i_2^{\prime*}C_{i_1}\ar[r]\ar[d]^{h_1} &C_{i_1}\ar[d]^h\\
S'\ar[r]^j&T
}
\]
instead of the cartesian diagram
\[
\xymatrix{
X_1\ar[r]\ar[d]&S\ar[d]^j\\
Y\ar[r]&T
}
\]
used to construct $\gamma_{i_1}$. Let $C_{i_2^{\prime\prime}}\to t_2^*N_i$ be the corresponding cone.

To avoid overloading the notation, we omit the subscripts on the refined Gysin maps $i^!_?$, $j^!_?$, etc., leaving the context to make the meaning clear. 

Let $0_{q',i}\:X_1\to q^{\prime*}N_i$, $0_{t_2,i}\:Z\to t_2^*N_i$, and
$0_{h_1,j}\:i_2^{\prime*}C_{i_1}\to h_1^*N_j$ be the respective 0-sections. For $y\in \sE^\BM_{a,b}(Y,v)$ we have
\[
j^!i^!(y)=j^![0^!_{q^{\prime},i}(\gamma_{i_1*}\sp_{i_1}^*(y))].
\]
Since the refined Gysin map commutes with proper push-forward, smooth pull-back, and pull-back with respect to 0-sections (Proposition~\ref{prop:LciProperties}) we have
\begin{align*}
j^!i^!(y)&=j^![0^!_{q',i}(\gamma_{i_1*}\sp_{i_1}^*(y))]\\
&=0^!_{t_2,i}j^!(\gamma_{i_1*}\sp_{i_1}^*(y))\\
&=0^!_{t_2,i}\gamma'_{i_1*}j^!\sp_{i_1}^*(y))\\
&=0^!_{t_2,i}\gamma'_{i_1*}0_{h_1,j}^!\gamma_{\alpha_1*} \sp_{\alpha_1}^*\sp_{i_1}^*(y).
\end{align*}

 Exchanging the role of $i$ and $j$, $p$ and $q$, etc.,  in \eqref{eqn:BigComDiag1} gives us the corresponding diagram
\[
\xymatrix{
C_{\alpha_2}\ar[d]\ar@{^(->}[r]^{\gamma_{\alpha_2}}&h_2^*N_i\ar[dl]\ar@/^20pt/[drrr]\\
i_1^{\prime*}C_{i_2}\ar[d]_{\alpha_2}\ar@{_(->}[dr]_{\gamma'_{i_2}}\ar[drr]^{g_2}
\ar@/^30pt/[drrrr]^{h_2}&&&&N_i\ar[d]\\
C_{i_2}\ar[drr]|-\hole\ar@{_(->}[dr]_{\gamma'_{i_2}}&t_1^*N_j\ar[r]\ar[d]^(.35){i_1^{\prime\prime}} &Z\ar@/^16pt/[dd]^(.28){t_1}\ar@/^16pt/[rr]^{t_2}\ar[r]^{i_2'}\ar[d]^{i_1'}&X_1\ar[r]^{q'}\ar[d]^{i_1}&S\ar[d]^i\\
&p^{\prime*}N_j\ar[r]&X_2\ar[r]^{i_2}|(.37)\hole\ar[d]^{p'}&Y\ar[d]_p\ar[r]^q&T\\
&&S'\ar[r]^j&T',
}
\]
and we have the respective 0-sections $0_{p',j}\:X_2\to p^{\prime*}N_j$, $0_{t_1,j}\:Z\to t_1^*N_j$, and
$0_{h_2,i}\:i_1^{\prime*}C_{i_2}\to h_2^*N_i$. 

We have the commutative diagram
\[
\xymatrix{
C_{\alpha_1}\ar[r]^{\gamma_{\alpha_1}}\ar[dr]_{\gamma_1}&h_1^*N_j\ar[r]\ar@{^(->}[d]^{\gamma_1''}&i_2^{\prime*}C_{i_1}\ar@{^(->}[d]^{\gamma'_{i_1}}\\
&t_2^*N_i\times_Zt_1^*N_j\ar[r]^-{p_1}&t_2^*N_i,
}
\]
with the square cartesian,  $\gamma_1''$ the pull-back of the closed immersion $\gamma'_{i_1}$, and $\gamma_1:=\gamma_1''\circ\gamma_{\alpha_1}$. Let $0_1\:t_2^*N_i\to 
t_2^*N_i\times_Zt_1^*N_j$ be the 0-section to $p_1$. Since the pull-back by a 0-section commutes with proper push-forward in a cartesian square, we have 
\begin{align*}
j^!i^!(y)&=0^!_{t_2,i}\gamma'_{i_1*}0_{h_1,j}^!\gamma_{\alpha_1*} \sp_{\alpha_1}^*\sp_{i_1}^*(y)\\
&=0^!_{t_2,i}0_1^!\gamma_{1*}\sp_{\alpha_1}^*\sp_{i_1}^*(y)\\
&=0_{1,2}^!\gamma_{1*}\sp_{\alpha_1}^*\sp_{i_1}^*(y),
\end{align*}
where $0_{1,2}=0_1\circ 0_{t_2,i}\:Z\to t_2^*N_i\times_Zt_1^*N_j$ is the zero-section.

By symmetry, we have the commutative diagram, with the square being cartesian,
\[
\xymatrix{
C_{\alpha_2}\ar[r]^{\gamma_{\alpha_2}}\ar[dr]_{\gamma_2}&h_2^*N_i\ar[r]\ar@{^(->}[d]^{\gamma_2''}&i_1^{\prime*}C_{i_2}\ar@{^(->}[d]^{\gamma'_{i_2}}\\
&t_2^*N_i\times_Zt_1^*N_j\ar[r]^-{p_1}&t_1^*N_j,
}
\]
and with $\gamma_2:=\gamma_2''\circ  \gamma_{\alpha_2}$. Letting $0_2\:t_1^*N_j\to t_2^*N_i\times_Zt_1^*N_j$ be the 0-section to $p_2$, we have the commutative diagram
\[
\xymatrix{
Z\ar[r]^-{0_{t_1,j}}\ar[d]_{0_{t_2,i}}\ar[dr]^{0_{1,2}}&t_1^*N_j\ar[d]^{0_2}\\
t_2^*N_i\ar[r]^-{0_1}&t_2^*N_i\times_Zt_1^*N_j, 
}
\]
that is,
\[
0_{1,2}=0_1\circ 0_{t_2,i} = 0_2\circ 0_{t_1,j}.
\]
Arguing as above gives us the identity
\[
i^!j^!(y)=0^!_{1,2}\gamma_{2*}\sp_{\alpha_2}^*\sp_{i_2}^*(y).
\]

We have the closed immersion 
\[
\gamma_{1,2}:=\gamma_{i_1}'\times \gamma_{i_2}'\:i_2^{\prime*}C_{i_1}\times_Zi_1^{\prime*}C_{i_2}\to t_2^*N_i\times_Zt_1^*N_j.
\]
Recalling the closed immersions
\[
\beta_1\:C_{\alpha_1}\to  i_2^{\prime*}C_{i_1}\times_Zi_1^{\prime*}C_{i_2},\ \beta_2\: C_{\alpha_2}\to i_2^*C_{i_1}\times_Zi_1^*C_{i_2}
\]
one checks that 
\[
\gamma_i=\gamma_{12}\circ\beta_i,\ i=1,2,
\]
so we have
\[
j^!i^!(y)=0^!_{1,2}\gamma_{12*}\beta_{1*}\sp_{\alpha_1}^*\sp_{i_1}^*(y)
\]
and
\[
i^!j^!(y)=0^!_{1,2}\gamma_{12*}\beta_{2*}\sp_{\alpha_2}^*\sp_{i_2}^*(y),
\]
so the identity $i^!j^!(y)=j^!i^!(y)$ follows from Proposition~\ref{prop:SpecializationComm}.
\end{proof}

For the proof of Proposition~\ref{prop:SpecializationComm}, we use a ``double deformation diagram'', 
\begin{equation}\label{eqn:DDDiag}
\pi_{12}\:\Def_{12}\to Y\times\A^1\times\A^1,
\end{equation}
which we now proceed to construct.

 Let $\Def_1=(\Def_{X_1}Y)\times\A^1$ and let $\Def_2=\Def_{X_2\times\A^1}(Y\times\A^1)$. $\Def_1$ is an open subscheme of $\Bl_{X_1\times0\times\A^1}Y\times\A^1\times\A^1$ and 
$\Def_2$ is an open subscheme of $\Bl_{X_2\times\A^1\times0}Y\times\A^1\times\A^1$, giving   structure maps $\pi_i\:\Def_i\to Y\times\A^1\times\A^1$, $i=1,2$. Let $\Def_{12}$ be the fiber product
\[
\Def_{12}:=\Def_1\times_{Y\times\A^1\times\A^1}\Def_2
\]
with structure morphism
\[
\pi_{12}\:\Def_{12}\to Y\times\A^1\times\A^1,
\]
and projections
\[
p_j\:\Def_{12}\to \Def_j.
\]

We have $\pi_1^{-1}(Y\times0\times\A^1)=C_{i_1}\times\A^1$ and $\pi_2^{-1}(Y\times\A^1\times0)=C_{i_2}\times\A^1$, giving closed immersions
\[
\sigma_j\:C_{i_j}\times\A^1\to \Def_j
\]
with respective open complements
\[
\eta_1\:Y\times(\A^1\setminus\{0\})\times\A^1\to \Def_1,\
\eta_2\:Y\times\A^1\times(\A^1\setminus\{0\})\to \Def_2.
\]

The proper transform of $\sigma_1(C_{i_1}\times\A^1)$ to $\Def_{12}$ is  an open subscheme of $\Bl_{i_2^*C_{i_1}\times0}C_{i_1}\times\A^1$, namely $\Def_{i_2^*C_{i_1}}C_{i_1}$, giving the closed immersion
\[
\hat\sigma_1\: \Def_{i_2^*C_{i_1}}C_{i_1}\to \Def_{12}
\]
over the closed immersion $Y\times0\times\A^1\to Y\times\A^1\times\A^1$. This gives us the commutative diagram
\[
\xymatrix{
\Def_{i_2^*C_{i_1}}C_{i_1}\setminus C_{\alpha_1}\ar@{^(->}[r]\ar[d]^{\bar{p}_1^\circ}_\wr&\Def_{i_2^*C_{i_1}}C_{i_1}\ar[r]^{\hat\sigma_1}\ar[d]^{\bar{p}_1}& \Def_{12}\ar[d]^{p_1}\ar@/^30pt/[dd]^{\pi_{12}}\\
C_{i_1}\times(\A^1\setminus\{0\})\ar[d]^{\bar{\pi}_1^\circ}\ar@{^(->}[r]&C_{i_1}\times\A^1\ar[r]^{\sigma_1}\ar[d]^{\bar{\pi}_1}&\Def_1\ar[d]^{\pi_1}\\
Y\times0\times(\A^1\setminus\{0\})\ar@{^(->}[r]&Y\times0\times\A^1\ar[r]&Y\times\A^1\times\A^1,
}
\]
with $\hat\sigma_1$ the canonical morphism.

Similarly, we have the closed immersion 
\[
\hat\sigma_2\: \Def_{i_1^*C_{i_2}}C_{i_2}\to \Def_{12}
\]
over the closed immersion $Y\times\A^1\times0\to Y\times\A^1\times\A^1$, which gives us the commutative diagram
\[
\xymatrix{
\Def_{i_1^*C_{i_2}}C_{i_2}\setminus C_{\alpha_2}\ar@{^(->}[r]\ar[d]^{\bar{p}_2^\circ}_\wr&\Def_{i_1^*C_{i_2}}C_{i_2}\ar[r]^{\hat\sigma_2}\ar[d]^{\bar{p}_2}& \Def_{12}\ar[d]^{p_2}\ar@/^30pt/[dd]^{\pi_{12}}\\
C_{i_2}\times(\A^1\setminus\{0\})\ar[d]^{\bar{\pi}_2^\circ}\ar@{^(->}[r]&C_{i_2}\times\A^1\ar[r]^{\sigma_2}\ar[d]^{\bar{\pi}_2}&\Def_2\ar[d]^{\pi_2}\\
Y\times(\A^1\setminus\{0\})\times0\ar@{^(->}[r]&Y\times\A^1\times0\ar[r]&Y\times\A^1\times\A^1,
}
\]
with $\hat\sigma_2$ the canonical morphism. 

Thus, the closed immersions
\[
\Def_{i_2^*C_{i_1}}C_{i_1}\hookrightarrow\pi_{12}^{-1}(Y\times0\times\A^1),\ 
\Def_{i_1^*C_{i_2}}C_{i_2}\hookrightarrow\pi_{12}^{-1}(Y\times\A^1\times0)
\]
are isomorphisms over $Y\times\A^1\times\A^1\setminus Y\times0\times0$.

We have
\[
\pi_{12}^{-1}(Y\times0\times0)=\pi_1^{-1}(Y\times0\times0)\times_Y\pi_2^{-1}(Y\times0\times0)=
i_2^*C_{i_1}\times_Yi_1^*C_{i_2}
\]
so
\[
\pi_{12}^{-1}(Y\times0\times\A^1)=\Def_{i_2^*C_{i_1}}C_{i_1}\cup i_2^*C_{i_1}\times_Yi_1^*C_{i_2}
\]
and
\[
\pi_{12}^{-1}(Y\times\A^1\times0)=\Def_{i_1^*C_{i_2}}C_{i_2}\cup i_2^*C_{i_1}\times_Yi_1^*C_{i_2}.
\]

The fiber of $\Def_{i_2^*C_{i_1}}C_{i_1}\to Y\times\A^1\times\A^1$ over $Y\times0\times0$ is $C_{\alpha_1}\subset \Def_{i_2^*C_{i_1}}C_{i_1}$ and the inclusion of this fiber in $\pi_{12}^{-1}(Y\times0\times0)$ is the closed immersion $\beta_1$. Similarly, the inclusion of $C_{\alpha_2}\subset \Def_{i_1^*C_{i_2}}C_{i_2}$ in $\pi_{12}^{-1}(Y\times0\times0)$ is the closed immersion $\beta_2$.

We fit this all together in the following commutative diagram
\begin{equation}\label{eqn:VistoliGeomConst}
\xymatrixcolsep{5pt}
\xymatrixrowsep{10pt}
\xymatrix{
&\Def_{i_2^*C_{i_1}}C_{i_1}\ar[rrr]^{\hat\sigma_1}\ar@{-}[d]&&&\Def_{12}\ar[dddd]^{\pi_{12}}\\
C_{\alpha_1}\ar@/_90pt/[ddddrr]\ar@{^(->}[rr]^{\beta_1}\ar@{^(->}[ur]^{\sigma_{\alpha_1}}&\ar@{-}[d]&i_2^*C_{i_1}\times_Yi_1^*C_{i_2}\ar[rru]\ar@{}[d]|(.15){}="a"\ar@{}[d]|(.9){}="b"\ar@{-}"a";"b"\\
&C_{\alpha_2}\ar@/_70pt/[dddr]\ar@{^(->}[ur]^(.35){\beta_2}\ar@{^(->}[rr]^(.63){\sigma_{\alpha_2}}\ar[dd]&\ar[ddd] &\Def_{i_1^*C_{i_2}}C_{i_2}\ar[ruu]_{\hat{\sigma}_2}\ar[dddd]\\
\\
&Y\times0\times\A^1\ar@{-}[rr]|(.5)\hole&&\ar[r]&Y\times\A^1\times\A^1\\
&&Y\times0\times0\ar[rru]|(.52)\hole\ar[ul]\ar[dr]\\
&&&Y\times\A^1\times0\ar[ruu]
}
\end{equation}
Moreover, $\pi_{12}$ restricts to an isomorphism
\[
\Def_{12}\setminus\pi_{12}^{-1}(Y\times0\times\A^1\cup Y\times\A^1\times0)\xymatrix{\ar[r]^{\pi_{12}}_\sim&}Y\times\A^1\times\A^1 \setminus(Y\times0\times\A^1\cup Y\times\A^1\times0),
\]
we have
\[
\pi_{12}^{-1}(Y\times0\times\A^1\cup Y\times\A^1\times0)=
\hat{\sigma}_1(\Def_{i_2^*C_{i_1}}C_{i_1})\cup  i_2^*C_{i_1}\times_Yi_1^*C_{i_2}\cup
\hat{\sigma}_2(\Def_{i_1^*C_{i_2}}C_{i_2}),
\]
and
\begin{gather*}
i_2^*C_{i_1}\times_Yi_1^*C_{i_2}\cap \hat{\sigma}_1(\Def_{i_2^*C_{i_1}}C_{i_1})=\beta_1(C_{\alpha_1}),\\
i_2^*C_{i_1}\times_Yi_1^*C_{i_2}\cap \hat{\sigma}_2(\Def_{i_1^*C_{i_2}}C_{i_2})=\beta_2(C_{\alpha_2}).
\end{gather*}

With these preparations, we now proceed to the proof of Proposition~\ref{prop:SpecializationComm}. As above, for $(W\to B)\in \Sch^G/B$, we let $p_W\:W\to B$ denote the structure morphism.

\begin{proof}[Proof of Proposition~\ref{prop:SpecializationComm}]
We first define a 2-variable version of the map $\vartheta_{Y,v}$ \eqref{eqn:Vartheta}:
 \[
 \vartheta_{Y,v,2}\:p_{Y\times(\A^1\setminus\{0\})^2!}\circ  \Sigma_{S^1}^2\circ \Sigma^v\circ p_{Y\times(\A^1\setminus\{0\})^2}^*\to
 p_{Y!}\circ  \Sigma^v\circ p_Y^*.
 \]
 Here it is important to keep track of the order of the two factors in the $Y$-scheme $Y\times (\A^1\setminus\{0\})^2$, so we will write this as
 \[
 Y\times (\A^1\setminus\{0\})^2=\Spec_{\sO_Y}\sO_Y[x_1^{\pm1}, x_2^{\pm1}],
 \]
 and we fix an isomorphism $\Omega_{ (\A^1\setminus\{0\})^2/k}\cong \sO_{ (\A^1\setminus\{0\})^2}^2$ by using the (ordered) basis $dx_1, dx_2$. We will distinguish the two factors by writing
 \[
 \A^1_i=\Spec_{\sO_Y} \sO_Y[x_i],\quad i=1,2.
 \]
 and define the pointed $Y$-scheme $\G_{m,i}$ by
 \[
\G_{m,i}: =(\A^1_i\setminus\{0\}, \{1\})
\]

We have the commutative diagram of projections
\[
\xymatrix{
&Y\times(\A^1_1\setminus\{0\})\times (\A^1_2\setminus\{0\})\ar[dl]_{p^{12}_1}\ar[dd]^{p^{12}_Y}\ar[dr]^{p^{12}_2}\\
Y\times(\A^1_1\setminus\{0\})\ar[dr]_{p^1_Y}&&Y\times(\A^1_2\setminus\{0\})\ar[dl]_{p^2_Y}\\
&Y
}
\]

We let 
\[
\vartheta_{p^{12}_2,v[1]}\:
p_{Y\times (\A^1\setminus\{0\})^2!}\Sigma_{S^1}^2\Sigma^v p_{Y\times (\A^1\setminus\{0\})^2}^*\to p_{Y\times (\A^1_2\setminus\{0\})!}\Sigma_{S^1}\Sigma^v p_{Y\times (\A^1_2\setminus\{0\})}^*
\]
be the composition
\begin{align*}
p_{Y\times (\A^1\setminus\{0\})^2!}&\Sigma_{S^1}^2\Sigma^v p_{Y\times (\A^1\setminus\{0\})^2}^*\\
&=
p_{Y\times (\A^1_2\setminus\{0\})!}p^{12}_{2!}\Sigma_{S^1}^2\Sigma^v p^{12*}_{2}p_{Y\times (\A^1_2\setminus\{0\})}^*\\
&=p_{Y\times (\A^1_2\setminus\{0\})!}p^{12}_{2\#}\Sigma^{-\sO_{Y\times(\A^1\setminus\{0\})^2}dx_1}\Sigma_{S^1}p^{12*}_{2}\Sigma_{S^1}\Sigma^v p_{Y\times (\A^1_2\setminus\{0\})}^*\\
&=p_{Y\times (\A^1_2\setminus\{0\})!}p^{12}_{2\#}p^{12*}_{2}\Sigma^{-\sO_{Y\times(\A^1_2\setminus\{0\})}dx_1}\Sigma_{S^1}\Sigma_{S^1}\Sigma^v p_{Y\times (\A^1_2\setminus\{0\})}^*\\
&=p_{Y\times (\A^1_2\setminus\{0\})!}\Sigma_{(\A^1_1\setminus\{0\})_+}\Sigma^{-\sO_{Y\times(\A^1_2\setminus\{0\})}dx_1}\Sigma_{S^1}\Sigma_{S^1}\Sigma^v p_{Y\times (\A^1_2\setminus\{0\})}^*\\
&\to p_{Y\times (\A^1_2\setminus\{0\})!}\Sigma_{\G_{m,1}}[\Sigma_{S^1}\Sigma_{\G_{m,1}}]^{-1}\Sigma_{S^1}\Sigma_{S^1}\Sigma^v p_{Y\times (\A^1_2\setminus\{0\})}^*\\
&= p_{Y\times (\A^1_2\setminus\{0\})!}\Sigma_{S^1}\Sigma^v p_{Y\times (\A^1_2\setminus\{0\})}^*.
\end{align*}
We define 
\[
\vartheta_{p^2_Y,v}\:p_{Y\times (\A^1_2\setminus\{0\})!}\Sigma_{S^1}\Sigma^v p_{Y\times (\A^1_2\setminus\{0\})}^*\to p_{Y!}\Sigma^vp_Y^*
\]
as $\vartheta_{Y,v}$ where we use our coordinate $x_2$ for  $\A^1\setminus\{0\}$. We let
\[
\vartheta_{Y,v,2}\:p_{Y\times (\A^1\setminus\{0\})^2!}\Sigma_{S^1}^2\Sigma^v p_{Y\times (\A^1\setminus\{0\})^2}^*\to p_{Y!}\Sigma^vp_Y^*
\]
be the composition
\[
\vartheta_{Y,v,2}:=\vartheta_{p^2_Y,v}\circ \vartheta_{p^{12}_2,v[1]}.
\]

Writing $p^{12}_Y=p^1_Y\circ p^{12}_1$ and proceeding as above gives us the maps
\[
\vartheta_{p^{12}_1,v[1]}\:
p_{Y\times (\A^1\setminus\{0\})^2!}\Sigma_{S^1}^2\Sigma^v p_{Y\times (\A^1\setminus\{0\})^2}^*\to p_{Y\times (\A^1_1\setminus\{0\})!}\Sigma_{S^1}\Sigma^v p_{Y\times (\A^1_1\setminus\{0\})}^*
\]
and
\[
\vartheta_{p^1_Y,v}\:p_{Y\times (\A^1_1\setminus\{0\})!}\Sigma_{S^1}\Sigma^v p_{Y\times (\A^1_1\setminus\{0\})}^*\to p_{Y!}\Sigma^vp_Y^*,
\]
and we define
\[
\vartheta_{Y,v,1}:=\vartheta_{p^1_Y,v}\circ \vartheta_{p^{12}_1,v[1]}.
\]

\begin{lemma}  \label{lem:Sign}  $\vartheta_{Y,v,1}=-\vartheta_{Y,v,2}$.
\end{lemma}

\begin{proof} We can write  $\vartheta_{Y,v,2}$ as the composition 
\begin{multline*}
p_{Y\times (\A^1\setminus\{0\})^2!}\Sigma_{S^1}^2\Sigma^v p_{Y\times (\A^1\setminus\{0\})^2}^*\\
\to p_{Y!}[\Sigma_{\G_{m,1}}\circ(\Sigma^{\sO_Ydx_1})^{-1}]\circ[\Sigma_{\G_{m,2}}\circ(\Sigma^{\sO_Ydx_2})^{-1}]\Sigma_{S^1}^2\Sigma^vp_Y^*\\
\cong 
p_{Y!}\Sigma_{S^1}^{-1}\circ\Sigma_{S^1}^{-1}\Sigma_{S^1}^2\Sigma^vp_Y^*
\cong 
p_{Y!}\Sigma^vp_Y^*
\end{multline*}
and $\vartheta_{Y,v,1}$ similarly as
\begin{multline*}
p_{Y\times (\A^1\setminus\{0\})^2!}\Sigma_{S^1}^2\Sigma^v p_{Y\times (\A^1\setminus\{0\})^2}^*\\
\to p_{Y!}[\Sigma_{\G_{m,2}}(\Sigma^{\sO_Ydx_2})^{-1}]\circ[\Sigma_{\G_{m,1}}(\Sigma^{\sO_Ydx_1})^{-1}]\Sigma_{S^1}^2\Sigma^vp_Y^*\\
\cong 
p_{Y!}\Sigma_{S^1}^{-1}\circ\Sigma_{S^1}^{-1}\Sigma_{S^1}^2\Sigma^vp_Y^*
\cong 
p_{Y!}\Sigma^vp_Y^*.
\end{multline*}
Acting  by the respective exchange of factors maps gives a commutative diagram
\[
\xymatrix{
[\Sigma_{\G_{m,1}}(\Sigma^{\sO_Ydx_1})^{-1}]\circ[\Sigma_{\G_{m,2}}(\Sigma^{\sO_Ydx_2})^{-1}]\ar[d]^{\tau}\ar[r]^-\sim&\Sigma_{S^1}^{-1}\circ\Sigma_{S^1}^{-1}\ar[d]^\tau\\
[\Sigma_{\G_{m,2}}(\Sigma^{\sO_Ydx_2})^{-1}]\circ[\Sigma_{\G_{m,1}}(\Sigma^{\sO_Ydx_1})^{-1}]\ar[r]^-\sim&
\Sigma_{S^1}^{-1}\circ\Sigma_{S^1}^{-1}
}
\]
and transforms $\vartheta_{Y,v,2}$ into $\vartheta_{Y,v,1}$. As $\tau:S^1\wedge S^1\to 
S^1\wedge S^1$ is multiplication by $-1$, this proves the lemma.
\end{proof}

Let $\Def_{12}^\circ=\Def_{12}\setminus \pi_{12}^{-1}(Y\times0\times0)$, let 
$\Def^\circ_{i_2^*C_{i_1}}C_{i_1}=\Def_{i_2^*C_{i_1}}C_{i_1}\setminus C_{\alpha_1}$ and let 
$\Def^\circ_{i_1^*C_{i_2}}C_{i_2}=\Def_{i_1^*C_{i_2}}C_{i_2}\setminus C_{\alpha_2}$. This gives us the closed immersion
\[
i^\circ\:\Def^\circ_{i_2^*C_{i_1}}C_{i_1}\amalg \Def^\circ_{i_1^*C_{i_2}}C_{i_2}\to \Def_{12}^\circ
\]
with open complement
\[
\eta^\circ\:Y\times(\A^1\setminus\{0\})^2\to \Def_{12}^\circ.
\]

The localization sequence for
\[
\Def^\circ_{i_2^*C_{i_1}}C_{i_1}\amalg \Def^\circ_{i_1^*C_{i_2}}C_{i_2}\xrightarrow{i^\circ} \Def_{12}^\circ\xleftarrow{\eta^\circ}Y\times(\A^1\setminus\{0\})^2
\]
gives the boundary map
\[
(\del_1^*,\del_2^*)\:\sE^\BM_{a+2,b}(Y\times(\A^1\setminus\{0\})^2,v)\to
\sE^\BM_{a+1,b}(\Def^\circ_{i_2^*C_{i_1}}C_{i_1},v)\oplus
\sE^\BM_{a+1,b}(\Def^\circ_{i_1^*C_{i_2}}C_{i_2}, v).
\]
We claim that  diagram
\begin{equation}\label{eqn:Sign2}
\xymatrix{
\sE^\BM_{a,b}(Y,v)\ar[dd]^{\vartheta^*_{Y,v,2}} \ar[rd]^-{(\sp^*_{i_1,v}, \sp^*_{i_2,v} )}\\&
\hbox to 100pt{\hss $\sE^\BM_{a,b}(C_{i_1},v)\oplus \sE^\BM_{a,b}(C_{i_2},v)$}\ar[dd]^{(-\bar{p}_1^{\circ*}\circ\vartheta^*_{C_{i_1},v}, \bar{p}_2^{\circ*}\circ\vartheta^*_{C_{i_2},v})}\\
\sE^\BM_{a+2,b}(Y\times(\A^1\setminus\{0\})^2,v)\ar[rd]^-{(\del_1^*,\del_2^*)}\\&
\hbox to 180pt{\hss $\sE^\BM_{a+1,b}(\Def^\circ_{i_2^*C_{i_1}}C_{i_1},v)\oplus
\sE^\BM_{a+1,b}(\Def^\circ_{i_1^*C_{i_2}}C_{i_2}, v)$}
}
\end{equation}
commutes.

To see this,  let 
\[
\del_{i_1}\:p_{C_{i_1}!}\Sigma^vp_{C_{i_1}}^*\to p_{Y\times(\A^1_1\setminus\{0\})!}\circ \Sigma_{S^1}\Sigma^vp_{Y\times(\A^1_1\setminus\{0\}}^*
\]
be the map induced from the boundary map in the localization triangle for 
\[
C_{i_1}\hookrightarrow \Def_{i_1}Y\hookleftarrow Y\times (\A^1\setminus\{0\}).
\]
We recall that $\pi_{12}^{-1}(Y\times\A^1\times (\A^1-\{0\})\cong \Def_{i_1}Y\times(\A^1_2\setminus\{0\})$ and this isomorphism restricts to the isomorphism $\bar{p}_1^\circ\:\Def^\circ_{i_2^*C_{i_1}}C_{i_1}\to C_{i_1}\times (\A^1_2\setminus\{0\})$, lying over $Y\times0\times(\A^1_2\setminus\{0\})$. It then follows directly from the definitions    that the diagram
\begin{equation}\label{eqn:ThetaComDiag}
\xymatrixcolsep{50pt}
\xymatrix{
p_{\Def^\circ_{i_2^*C_{i_1}}C_{i_1}!}\Sigma_{S^1}\Sigma^vp_{\Def^\circ_{i_2^*C_{i_1}}C_{i_1}}^*\ar[r]^-{\vartheta_{C_{i_1}}\circ \bar{p}_1^\circ}\ar[d]^{\del_1}&p_{C_{i_1}!}\Sigma^vp_{C_{i_1}}^*\ar[d]^{\del_{i_1}}\\
p_{Y\times(\A^1\setminus\{0\})^2!}\circ \Sigma_{S^1}^2\Sigma^vp_{Y\times(\A^1\setminus\{0\})^2}^*\ar[r]^-{\vartheta_{p^1_Y,v[1]}}&
p_{Y\times(\A^1_1\setminus\{0\})!}\circ \Sigma_{S^1}\Sigma^vp_{Y\times(\A^1_1\setminus\{0\}}^*
}
\end{equation}
commutes.

Recall that the specialization map $\sp_{i_1,v}\:p_{C_{i_1}!}\Sigma^vp_{C_{i_1}}^*\to 
p_{Y!}\Sigma^vp_{Y}^*$ is defined as
\[
\sp_{i_1,v}:=\vartheta_{Y,v}\circ \del_{i_1}.
\]
Using Lemma~\ref{lem:Sign} and the commutativity of \eqref{eqn:ThetaComDiag}, we have
\begin{align*}
\vartheta_{Y,v,2}\circ\del_1&=-\vartheta_{Y,v,1}\circ\del_1\\
&=-\vartheta_{Y,v}\circ \vartheta_{p^1_Y,v[1]}\circ \del_1\\
&=-\vartheta_{Y,v}\circ \del_{i_1}\circ \vartheta_{C_{i_1},v}\circ \bar{p}_1^{\circ}\\
&=-\sp_{i_1,v}\circ (\vartheta_{C_{i_1},v}\circ \bar{p}_1^{\circ}).
\end{align*}

The argument for  $C_{i_2}$ is similar, where just use the definition of $\vartheta_{Y,v,2}$ as 
$\vartheta_{Y,v}\circ \vartheta_{p^2_Y,v[1]}$, and the identity
\[
\vartheta_{p^2_Y,v[1]}\circ \del_2=\del_{i_2}\circ \vartheta_{C_{i_2},v}\circ \bar{p}_2^{\circ},
\]
analogous to the commutativity of \eqref{eqn:ThetaComDiag}, where $\del_{i_2}$ is induced by the boundary map in the localization triangle for 
\[
C_{i_2}\hookrightarrow \Def_{i_2}Y\hookleftarrow Y\times (\A^1\setminus\{0\}).
\]
This gives
\begin{align*}
\vartheta_{Y,v,2}\circ\del_2&=
\vartheta_{Y,v}\circ \vartheta_{p^2_Y,v[1]}\circ \del_2\\
&=\vartheta_{Y,v}\circ \del_{i_2}\circ \vartheta_{C_{i_2},v}\circ \bar{p}_2^{\circ}\\
&=\sp_{i_2,v}\circ (\vartheta_{C_{i_2},v}\circ \bar{p}_2^{\circ}),
\end{align*}
which completes the proof of the commutativity of \eqref{eqn:Sign2}.

Now let  $\Def^{\circ'}_{12}=\Def_{12}\setminus(\beta_1(C_{\alpha_1})\cup \beta_2(C_{\alpha_2})$ 
and let  $i_2^*C_{i_1}\times_Yi_1^*C_{i_2}^{\circ'}=i_2^*C_{i_1}\times_Yi_1^*C_{i_2}\setminus(\beta_1(C_{\alpha_1})\cup \beta_2(C_{\alpha_2})$. The closed immersion $i^\circ$ extends to a closed immersion
\[
i^{\circ'}\: \Def^\circ_{i_2^*C_{i_1}}C_{i_1}\amalg \Def^\circ_{i_1^*C_{i_2}}C_{i_2}\amalg 
i_2^*C_{i_1}\times_Yi_1^*C_{i_2}^{\circ'}\to \Def_{12}^{\circ'}
\]
with open complement
\[
\eta^{\circ'}\:Y\times(\A^1\setminus\{0\})^2\to  \Def_{12}^{\circ'}.
\]
The corresponding localization sequence gives us the boundary map
\begin{multline*}
\sE^\BM_{a+2,b}(Y\times(\A^1\setminus\{0\})^2,v)
\xrightarrow{(\del_1,\del_2,\del_3')}
\left[\raise40pt\hbox{$
\xymatrixrowsep{2pt}
\xymatrix{\sE^\BM_{a+1,b}(\Def^\circ_{i_2^*C_{i_1}}C_{i_1},v)\\\oplus\\
\sE^\BM_{a+1,b}(\Def^\circ_{i_1^*C_{i_2}}C_{i_2},v)\\\oplus \\
\sE^\BM_{a+1,b}(i_2^*C_{i_1}\times_Yi_1^*C_{i_2}^{\circ'},v).}$}\right].
\end{multline*}

We have the closed immersion
\[
i\:\Def_{i_2^*C_{i_1}}C_{i_1}\cup \Def_{i_1^*C_{i_2}}C_{i_2}\cup
i_2^*C_{i_1}\times_Yi_1^*C_{i_2}\to \Def_{12},
\]
with open complement
\[
\eta\:Y\times(\A^1\setminus\{0\})^2\to \Def_{12}
\]
and the closed immersion
\[
\bar{i}\:\beta_1(C_{\alpha_1})\cup \beta_2(C_{\alpha_2})\to 
\Def_{i_2^*C_{i_1}}C_{i_1}\cup \Def_{i_1^*C_{i_2}}C_{i_2}\cup
i_2^*C_{i_1}\times_Yi_1^*C_{i_2}
\]
with open complement
\begin{multline*}
\Def^\circ_{i_2^*C_{i_1}}C_{i_1}\amalg \Def^\circ_{i_1^*C_{i_2}}C_{i_2}\amalg 
i_2^*C_{i_1}\times_Yi_1^*C_{i_2}^{\circ'}\\\xrightarrow{\bar{\eta}}
\Def_{i_2^*C_{i_1}}C_{i_1}\cup \Def_{i_1^*C_{i_2}}C_{i_2}\cup
i_2^*C_{i_1}\times_Yi_1^*C_{i_2}.
\end{multline*}
This gives us the boundary maps
\begin{multline*}
\sE^\BM_{a+2,b}(Y\times(\A^1\setminus\{0\})^2, v)\\\xrightarrow{\del}
\sE^\BM_{a+1,b}(\Def_{i_2^*C_{i_1}}C_{i_1}\cup \Def_{i_1^*C_{i_2}}C_{i_2}\cup
i_2^*C_{i_1}\times_Yi_1^*C_{i_2},v)
\end{multline*}
and
\[
\left[\hbox{$\begin{matrix}
\sE^\BM_{a+1,b}(\Def^\circ_{i_2^*C_{i_1}}C_{i_1},v)\\\oplus\\
\sE^\BM_{a+1,b}(\Def^\circ_{i_1^*C_{i_2}}C_{i_2},v)\\\oplus \\
\sE^\BM_{a+1,b}(i_2^*C_{i_1}\times_Yi_1^*C_{i_2}^{\circ'},v)\end{matrix}$}\right]
\quad \xrightarrow{\bar{\del}}\quad
\sE^\BM_{a,b}(\beta_1(C_{\alpha_1})\cup \beta_2(C_{\alpha_2}),v),
\]
fitting into their respective long exact localization sequences.

Now take $y\in \sE^\BM_{a+2,b}(Y\times(\A^1\setminus\{0\})^2, v)$. We claim that
\begin{equation}\label{eqn:MainSpecializationIdentity}
\bar{\del}((\del_1,\del_2, \del_3')(y))=0\in \sE^\BM_{a,b}(\beta_1(C_{\alpha_1})\cup \beta_2(C_{\alpha_2}),v).
\end{equation}
Indeed, we have the commutative diagram
\[
\xymatrix{
\hbox to 100pt{$\sE^\BM_{a+2,b}(Y\times(\A^1\setminus\{0\})^2, v)$\hss}\ar[dr]^\del\ar@/_50pt/[ddr]_-{(\del_1,\del_2, \del_3')}\\
&\hbox to 220pt{\hss$\sE^\BM_{a+1,b}(\Def_{i_2^*C_{i_1}}C_{i_1}\cup \Def_{i_1^*C_{i_2}}C_{i_2}\cup
i_2^*C_{i_1}\times_Yi_1^*C_{i_2},v)$}\ar[d]^-{\bar{\eta}^*}&\\
&\left[\hbox{$\begin{matrix}\sE^\BM_{a+1,b}(\Def^\circ_{i_2^*C_{i_1}}C_{i_1},v)\\ 
\oplus\\
\sE^\BM_{a+1,b}(\Def^\circ_{i_1^*C_{i_2}}C_{i_2},v)\\
\oplus\\ 
\sE^\BM_{a+1,b}(i_2^*C_{i_1}\times_Yi_1^*C_{i_2}^{\circ'},v)\end{matrix}$}\right]
\ar[d]^-{\bar{\del}}\\
&\sE^\BM_{a,b}(\beta_1(C_{\alpha_1})\cup \beta_2(C_{\alpha_2}),v)
}
\]
with the right-hand column a part of the long exact localization sequence for $\bar{i}, \bar{\eta}$. 

Let
\[
\bar{\del}_1\:\sE^\BM_{a+1,b}(\Def^\circ_{i_1^*C_{i_2}}C_{i_2},v)\to \sE^\BM_{a,b}(C_{\alpha_1},v)
\]
be the boundary map in the localization sequence for
\[
C_{\alpha_1}\hookrightarrow \Def_{i_1^*C_{i_2}}C_{i_2}\hookleftarrow
\Def^\circ_{i_1^*C_{i_2}}C_{i_2}
\]
and define
\[
\bar{\del}_2\:\sE^\BM_{a+1,b}(\Def^\circ_{i_2^*C_{i_1}}C_{i_1},v)\to \sE^\BM_{a,b}(C_{\alpha_2},v)
\]
similarly. Let
\[
\bar{\del}_3\:\sE^\BM_{a+1,b}(i_2^*C_{i_1}\times_Yi_1^*C_{i_2}^{\circ'},v)\to
\sE^\BM_{a,b}(\beta_1(C_{\alpha_1})\cup \beta_2(C_{\alpha_2}),v)
\]
be the boundary map for the localization sequence for
\[
\beta_1(C_{\alpha_1})\cup \beta_2(C_{\alpha_2})\hookrightarrow
i_2^*C_{i_1}\times_Yi_1^*C_{i_2}\hookleftarrow
i_2^*C_{i_1}\times_Yi_1^*C_{i_2}^{\circ'}
\]

Let $\bar{\beta}_j\:C_{\alpha_j}\to \beta_1(C_{\alpha_1})\cup \beta_2(C_{\alpha_2})$ be the closed immersion induced by $\beta_j$. Then
\[
\bar{\beta}_{1*}(\bar{\del}_1(x_1))=\bar{\del}(x_1,0,0),\
\bar{\beta}_{2*}(\bar{\del}_2(x_2))=\bar{\del}(0, x_2,0),\ 
\bar{\del}_3(x_3)=\bar{\del}(0,0,x_3)
\]
for 
\begin{gather*}
x_1\in \sE^\BM_{a+1,b}(\Def^\circ_{i_1^*C_{i_2}}C_{i_2},v),\ x_2\in \sE^\BM_{a+1,b}(\Def^\circ_{i_2^*C_{i_1}}C_{i_1},v),\\
x_3\in  
\sE^\BM_{a+1,b}((i_2^*C_{i_1}\times_Yi_1^*C_{i_2})^{\circ'},v).
\end{gather*}

Take $y\in \sE^\BM_{a+2, b}(Y,v)$. Putting this all together with \eqref{eqn:MainSpecializationIdentity} and the commutativity of \eqref{eqn:Sign2},  and inserting the definition of the various specialization maps,  we arrive at the identity
\[
\bar{\beta}_{2*}(\sp_{\alpha_2}^*(\sp_{i_2}^*(y)))-\bar{\beta}_{1*}(\sp_{\alpha_1}^*(\sp_{i_1}^*(y)))
=\bar{\del}_3(\del_3'(y))
\]
in $\sE^\BM_{a,b}(\beta_1(C_{\alpha_1})\cup \beta_2(C_{\alpha_2}),v)$. Pushing forward to
$\sE^\BM_{a,b}(i_2^*C_{i_1}\times_Yi_1^*C_{i_2},v)$, and applying the localization sequence for $\beta_1(C_{\alpha_1})\cup \beta_2(C_{\alpha_2})\hookrightarrow
i_2^*C_{i_1}\times_Yi_1^*C_{i_2}$ yields the identity
\[
\beta_{1*}(\sp_{\alpha_1}^*(\sp_{i_1}^*(y)))=\beta_{2*}(\sp_{\alpha_2}^*(\sp_{i_2}^*(y)))
\]
in $\sE^\BM_{a,b}(i_2^*C_{i_1}\times_Yi_1^*C_{i_2},v)$, which completes the proof.
\end{proof}

We conclude this section with a computation of the relative refined Gysin pull-back the fundamental class of a normal cone (Definition~\ref{def:FundClassCone}). Since this involves a change of base-scheme, we reintroduce the $-/B$ to the notation. 

Consider a commutative diagram in $\Sch^G/B$ with all squares cartesian
\begin{equation}\label{eqn:CartDiagTorInd}
\xymatrix{
X_0\ar[r]^{\bar{f}}\ar[d]_{i_0}&X_1\ar[d]^{i_1}\\
Y_0\ar[r]^f\ar[d]^{\pi_0}&Y_1\ar[d]^{\pi_1}\\
B_0\ar[r]^{\hat{f}}&B_1,
}
\end{equation}
with $i_1$ a closed immersion and with $\pi_1$ and $f_0$ in $\Sm^G/B_1$.  Thus $f$ and $\pi_0$ are in $\Sm^G/B$ and $i_0$ is a closed immersion.

This gives the diagram with all squares cartesian
\begin{equation}\label{eqnBigCartesianConeDiag}
\xymatrix{
f^*C_{X_1/Y_1}\ar[r]^{\bar{\bar{f}}}\ar[d]_{\bar{q}}\ar@/_20pt/[dd]_q&C_{X_1/Y_1}\ar[d]^{\bar{p}}\ar@/^20pt/[dd]^p\\
X_0\ar[r]^{\bar{f}}\ar[d]_{i_0}&X_1\ar[d]^{i_1}\\
Y_0\ar[r]^f\ar[d]^{\pi_0}&Y_1\ar[d]^{\pi_1}\\
B_0\ar[r]^{\hat{f}}&B_1\hbox to0pt{\,.}}
\end{equation}

We have the relative refined lci pull-back (\S\ref{subsec:RelPull-back})
\[
\hat{f}^!_{\pi_1\circ p/B_0}\:\sE^\BM(C_{X_1/Y_1}/B_1, \Omega_{Y_1/B_1})\to \sE^\BM(f^*C_{X_1/Y_1}/B_0, \Omega_{Y_0/B_0}),
\]
the fundamental classes
\[
[C_{X_1/Y_1}/B_1]\in \sE^\BM(C_{X_1/Y_1}/B_1, \Omega_{Y_1/B_1}),\ [C_{X_0/Y_0}/B_0]\in \sE^\BM(C_{X_0/Y_0}/B_0, \Omega_{Y_0/B_0}),
\]
and the closed immersion $\beta\:C_{X_0/Y_0}\to f^*C_{X_1/Y_1}$, inducing the map
\[
\beta_*\:\sE^\BM(C_{X_0/Y_0}/B_0, \Omega_{Y_0/B_0})\to \sE^\BM(f^*C_{X_1/Y_1}/B_0, \Omega_{Y_0/B_0}).
\]
 
\begin{proposition}\label{prop:FundClassIdentity} We have the identity
\[
\beta_*[C_{X_0/Y_0}/B_0]=\hat{f}^!_{\pi_1\circ p/B_0}[C_{X_1/Y_1}/B_1]
\]
in $\sE^\BM(\iota^*C_{X/Y}/B_0, \Omega_{Y_0/B_0})$.
\end{proposition}

\begin{proof} Suppose first that $\hat{f}$ is a smooth morphism. Then in the diagram \eqref{eqnBigCartesianConeDiag}, the maps $f$, $\bar{f}$ and $\overline{\bar{f}}$ are all smooth,   all squares are Tor-independent, and $\beta:C_{X_0/Y_0}\to f^*C_{X_1/Y_1}$ is an isomorphism. Using the relative purity isomorphism,  this reduces us to showing that $\bar{\bar{f}}^![C_{X_1/Y_1}/B]=[C_{X_0/Y_0}/B]$, which is an easy consequence of the definition of the fundamental classes 
$[C_{X_0/Y_0}/B]$ and $[C_{X_1/Y_1}/B]$, and the compatibility of specialization and localization sequences with smooth  pullback (Propositions~\ref{prop:SpecSmoothCommute})  and \ref{prop:LocSmoothComp}). 

This reduces us to the case of $\hat{f}$ a regular immersion. Since $\pi_1$ is smooth, $f$ is also a regular immersion and, changing the relative refined Gysin pullback back to the usual refined Gysin pullback via the relative purity isomorphism, we  reduce to showing that $\beta_*[C_{X_0/Y_0}/B]=f^!_p[C_{X_1/Y_1}/B]$. We can ignore the last row in the diagram \eqref{eqnBigCartesianConeDiag} and, since we are now working over $B$, we drop the $-/B$ from the notation

This comes down to a special case of Proposition~\ref{prop:SpecializationComm}. In our situation, the diagram \eqref{eqn:DoubleConeDiagram} is replaced with the following; here $N_f\to Y_0$ is the normal bundle of $f$ and  $\gamma:C_{X_0/Y_0}\to X_0$ is the structure morphism.
\[
\xymatrix{
&C_{\bar{\bar{f}}}\ar[d]_{\beta_1}\ar[dr]\\
\bar{r}^*C_{X_0/Y_0}\ar[dr]_{\bar{r}^*(\gamma)}\ar[r]^-{\beta_2}&q^*N_f\ar[r]^{\bar{\bar{r}}}\ar[d]
&f^{*}C_{X_1/Y_1}\ar[r]^{\bar{\bar{f}}}\ar[d]^{\bar{q}}&C_{X_1/Y_1}\ar[d]^{\bar{p}}\\
&i_0^{*}N_f\ar[r]^{\bar{r}}\ar[d]^{\alpha_2}&X_0\ar[r]^{\bar{f}}\ar[d]^{i_0}&X_1\ar[d]^{i_1}\\
&N_f\ar[r]^r&Y_0\ar[r]^f&Y_1.
}
\]
The map $\beta_2$ is the pullback $\bar{\bar{r}}^*(\beta)$ via the identities $q^*N_f=\bar{q}^*i_0^*N_f=\bar{r}^*f^*C_{X_1/Y_1}$ and 
\begin{equation}\label{eqn:VBId}
\bar{r}^*C_{X_0/Y_0}=\beta^*q^*N_f,  
\end{equation}
giving the the commutative diagram with both squares cartesian
\[
\xymatrix{
q^*N_f\ar[d]^{\bar{\bar{r}}}&\bar{r}^*C_{X_0/Y_0}\ar[d]^{\tilde{r}}\ar[r]^{\bar{r}^*(\gamma)}\ar[l]_{\beta_2}&i_0^*N_f\ar[d]^{\bar{r}}\\
f^*C_{X_1/Y_1}&C_{X_0/Y_0}\ar[r]^\gamma\ar[l]_\beta&X_0
}
\]
The map $\beta_1\:C_{\bar{\bar{f}}}\to q^*N_f$ is the map used to define $f^!_p$:
\begin{equation}\label{eqn:f!Def}
f^!_p=0^!_{q^*N_f}\circ \beta_{1*}\circ \sp_{\bar{\bar{f}}}^*.
\end{equation}
 
Via the Poincar\'e duality isomorphism, $\sE^{0,0}(Y_1)\cong \sE^\BM(Y_1/B, \Omega_{Y_1/B})$, 
the unit $1_{Y_1}$ gives the fundamental class $[Y_1]\in \sE^\BM(Y_1/B, \Omega_{Y_1/B})$.  Since $[C_{X_1/Y_1}]=\sp_{i_1}^*[Y_1]$, applying Proposition~\ref{prop:SpecializationComm} to  $[Y_1]$ yields the identity
 \begin{equation}\label{eqn:SpecRel2}
 \beta_{1*}(\sp_{\bar{\bar{f}}}^*[C_{X_1/Y_1}])=
 \beta_{2*}(\sp_{\alpha_2}^*(\sp_{f}^*[Y_1]))
 \end{equation}
 in $\sE^\BM(q^*N_\iota/B, \Omega_{Y_1/B})$. 
 
 Since $f$ is a regular immersion, we have $\sp_f^*[Y_1]=[N_f]$. Using Proposition~\ref{prop:SpecSmoothCommute} yields
 \begin{equation}\label{eqn:SpecRel3}
 \sp_{\alpha_2}^*[N_f]=\sp_{\alpha_2}^*(r^![Y_0])=\tilde{r}^!(\sp_{i_0}^*[Y_0])=\tilde{r}^!([C_{X_0/Y_0}]).
 \end{equation}
 We then have
 \begin{align*}
 f^!_p[C_{X_1/Y_1}]&=0^!_{q^*N_f}\beta_{1*}(\sp_{\bar{\bar{f}}}^*[C_{X_1/Y_1}])&\text{ \eqref{eqn:f!Def}}\\
 &=
 0^!_{q^*N_f}\beta_{2*}\tilde{r}^!([C_{X_0/Y_0}])&\text{\eqref{eqn:SpecRel2}, \eqref{eqn:SpecRel3}}\\
 &=\beta_*0^!_{\beta^*q^*N_f}\tilde{r}^!([C_{X_0/Y_0}])&\text{\eqref{eqn:VBId}, Proposition~\ref{prop:LciProperties}(1,3)}\\
 &=\beta_*[C_{X_0/Y_0}].
 \end{align*}
\end{proof}

\begin{remark}\label{rem:FundClassIdentity} Suppose that the top square in \eqref{eqn:CartDiagTorInd} is also Tor-independent. Then
\[
[C_{X_0/Y_0}/B_0]=\hat{f}_{\pi_1 p}^![C_{X_1/Y_1}/B_1]
\]
Indeed, this additional Tor-independence implies that the map $\beta\:C_{X_0/Y_0}\to f^*C_{X_1/Y_1}:=C_{X_1/Y_1}\times_{X_1}X_0$ is an isomorphism.
\end{remark}

\begin{remark}\label{rem:VistoliLemChangeBase} We may apply Proposition~\ref{prop:FundClassIdentity} in the case $B=Y_1$, giving the identity
\[
\beta_*[C_{X_0/Y_0}/Y_0]=f_{p/Y_0}^![C_{X_1/Y_1}/Y_1]
\]
in $\sE^\BM(f^*C_{X_1/Y_1}/Y_0)$. Similarly, under the hypotheses of Proposition~\ref{prop:FundClassIdentity}, applying the relative purity isomorphism gives the identity 
\[
\beta_*[C_{X_0/Y_0}/B]=\hat{f}_{\pi_1p}^![C_{X_1/Y_1}/B]
\]
in $\sE^\BM(f^*C_{X_1/Y_1}/B, \Omega_{Y_0/B})$.
\end{remark}

\section{Equivariant (co)homology via the algebraic Borel construction}\label{sec:TEGEquivCoh}  For   $G\to B$ tame,  $\sE\in \SH^G(B)$, $X\in \Sch^G/B$, we have been using the $\sE$-valued Borel-Moore homology defined  by
\[
\sE^\BM_{a,b}(X/B,G,v):=\Hom_{\SH^G(B)}(\Sigma^{a,b}p_{X!}(\Sigma^v1_X),\sE).
\]
For $\sE\in \SH(B)$, we also  have the equivariant Borel-Moore homology defined by the algebraic version of the Borel construction. We have discussed this construction in \cite[\S 4]{LevineAtiyahBott}, following the construction in \cite[\S 1]{DiLorenzoMantovani}, yielding a theory $\sE^\BM_G(X/BG, v)$; we use the same construction here, which we briefly recall.

We assume that $G$ is a closed subgroup scheme of $\GL_n/B$ for some $n$ and that $G$ is smooth over $B$; we no longer need to assume that $G$ is tame. However, we prefer to work in a certain full subcategory $\Sch^G_q/B$ of $\Sch^G/B$, for which the various quotients that we use are representable as quasi-projective $B$-schemes. This is probably not a necessary restriction, as the fppf sheaf quotient of a free $G$-action give a well-defined object in $\Spc(B)$, and this is the approach used, e.g., in \cite[\S4]{MorelVoevodsky}. Alternatively, the quotient will exist as a quasi-separated algebraic space, and one could use the extension of the six-functor formalism to quasi-separated algebraic spaces furnished by \cite{Chowdhury} or \cite{KR}. However, we prefer to keep things concrete and so we  will restrict to representable quotients. As we note below,  $\Sch^G_q/B=\Sch^G/B$ for $B=\Spec k$, $k$ a field, which is the setting for which we are ultimately interested.

\subsection{Construction and first properties of equivariant cohomology and equivariant Borel-Moore homology}

 Let $E_j\GL_n$ be the $B$-scheme of $n\times n+j$ matrices of rank $n$, with the action of $\GL_n$ by left multiplication. We let $E_jG$ denote the scheme $E_j\GL_n$ with the $G$-action by restriction with respect to a given embedding (that is, a homomorphism of group schemes that is a closed immersion of $B$-schemes) $G\hookrightarrow \GL_n/B$, which we fix throughout this section.
 
 \begin{definition}
We let $\Sch^G_q/B$ be the full subcategory of $\Sch^G/B$ consisting of those $X$ for which all the fppf quotients  $G\backslash X\times_BE_jG$ are represented in $\Sch/B$; here $G$ acts diagonally on $X\times_BE_jG$, and we usually denote the representing scheme by $X\times^GE_jG$. We similarly define $\Sm^G_q/B:=\Sm^G/B\cap \Sch^G_q/B$.
\end{definition}
\begin{remark} 1. In case $B=\Spec k$, $k$ a field, we have $\Sch^G_q/B=\Sch^G/B$, by \cite[Lemma 9, Proposition 23]{EGEquivIntThy}.\\[2pt]
2. Suppose $X$ is in $\Sch^G_q/B$ and $U\subset X$ is a $G$-stable open subscheme. Then $U$ is in $\Sch^G_q/B$ and $G\backslash U\times_BE_jG$ is represented by an open 
subscheme $U\times^GE_jG$ of $X\times^GE_jG$.\\[2pt]
3. Since $G$ is assumed to be smooth over $B$ and the action of $G$ on $X\times_BE_jG$ is free,  the quotient map $X\times_BE_jG\to  X\times^GE_jG$ defines a $G$-torsor for the \'etale topology for each $X\in \Sch^G_q/B$.
\end{remark}

We assume  throughout that $B$ is in $\Sch^G_q/B$, that is, the quotient  $B_jG:=G\backslash E_jG$ exists as a quasi-projective $B$-scheme and the quotient map $E_jG\to B_jG$ is a $G$-torsor for the \'etale topology. Let $F=\A^n_B$ with the left action of $\GL_n$ via the usual matrix multiplication; we consider $F$ as a $G$-representation via the inclusion $G\hookrightarrow \GL_n$. 

We have the $G$-equivariant closed immersion $i_j\: E_jG\to E_{j+1}G$ defined by adding a $j+1$st column of zeroes, and the $G$-equivariant  open immersion $\eta_j\:E_jG\times F\to E_{j+1}G$ by considering $F$ as the $j+1$st column. Let $p_j\:E_jG\times F\to E_jG$ be the projection with 0-section $i^0_j$. 

Let $B_jG=B\times^GE_jG$, $N_jG= B\times^G (E_jG\times F)$.  We have the corresponding regular immersions  $i_{G,j}\:B_jG\to B_{j+1}G$, $i^0_{G,j}\:B_jG\to N_jG$,  open immersion $\eta_{G,j}\:N_jG\to B_{j+1}G$  and projection $p_{j,G}\:N_jG\to B_jG$, making $N_jG$ a vector bundle over $B_jG$ with zero-section $i^0_{G,j}$, isomorphic to the normal bundle of $i_{G,j}$. This gives us the Ind-objects $EG:=\{E_jG\}_j$, $BG:=\{B_jG\}_j$  in $\Sm^G/B$, with transition maps $i_j$, $i_{G,j}$ all regular immersions, and the morphism of ind-objects $p_G\:EG\to BG$, making $EG$ an ind-$G$-torsor over $BG$.  

Take $X\in \Sch^G_q/B$. We have the maps $i_{X,j}\:X\times_BE_jG\to X\times_BE_{j+1}G$,  $i^0_{X,j}\:X\times_BE_jG\to X\times_B E_jG\times F$, $\eta_{X,j}\:X\times_BE_{j}G\times F\to X\times_BE_{j+1}G$, $p_{X,j}\:X\times_BE_{j}G\times F\to
X\times_BE_jG$, with $i_{X,j}=\id_X\times_B i_j$, etc. These yield  the corresponding maps on the quotients,
$i_{X,G,j}\:X\times^GE_jG\to X\times^GE_{j+1}G$, $i^0_{X,G,j}$, $\eta_{X, G,j}$ and
$p_{X,G,j}$.  The projections $\pi_{X,j}\:X\times_BE_jG\to E_jG$ induce the projections
$\pi_{X,G,j}\:X\times^GE_jG\to B_jG$. We denote the ind-object $\{X\times^GE_jG, i_{X,G,j}\}_j$  by $X\times^GEG$. The maps   $\pi_{X,G,j}\:X\times^GE_jG\to B_jG$ give the corresponding map of ind-objects $\pi_{X,G}\:X\times^GEG\to BG$. 

We have the cartesian, Tor-independent diagram
\begin{equation}\label{eqn:StructureDiagr}
\xymatrix{
X\times^GE_jG\ar[r]^{i_{X,G,j}}\ar[d]_{\pi_{X,G,j}}&X\times^GE_{j+1}G \ar[d]_{\pi_{X,G,j+1}}\\
B_jG\ar[r]^{i_{G,j}}& B_{j+1}G.
}
\end{equation}

Given $w\in \sK^G(X)$, we have the corresponding perfect complex $w_j$ on $X\times^GE_j$ defined by applying $G$-descent to the $G$-linearized perfect complex $p_1^*(w)$ on $X\times_BE_jG$. This gives us the  map of groupoids
\[
(-)_j\:\sK^G(X)\to \sK(X\times^GE_jG)
\]
with canonical natural isomorphisms $w_j\cong i_{X,G,j}^*(w_{j+1})$. Given an $\sE\in \SH(B)$, the regular immersion $i_{G,j}$ thus gives us the relative pullback map (\S\ref{subsec:RelPull-back})
\begin{equation}\label{eqn:pull-back1}
i^!_{G,j}\:\sE^\BM_{**}(X\times^GE_{j+1}G/B_{j+1}G, w_{j+1})\to
\sE^\BM_{**}(X\times^GE_jG/B_jG, w_j-\sN_{i_{G,j}}).
\end{equation}

This gives us the pro-system $\{\sE^\BM_{a,b}(X\times^GE_jG/B_jG, w_j), i_{X,G,j}^!\}_j$. Similarly,  the maps $i_{X,G,j}^*\:\sE^{a,b}(X\times^GE_{j+1}G, w_{j+1})\to \sE^{a,b}(X\times^GE_jG, w_j)$ define the pro-system $\{\sE^{a,b}(X\times^GE_jG, w_j), i_{X,G,j}^*\}_j$. This justifies the following definition.
\begin{definition}
1. For $X\in \Sch^G/B$ and $w\in \sK^G(X)$, define the $w$-twisted $G$-equivariant Borel-Moore homology by
\[
\sE_{G, a,b}^\BM(X/BG, w):=\lim_{\substack{\leftarrow\\ j}}\sE^\BM_{a,b}(X\times^GE_jG/B_jG,  w_j),
\]
and define the $w$-twisted $G$-equivariant cohomology as
\[
\sE^{a,b}_G(X, w):=\lim_{\substack{\leftarrow\\j}}\sE^{a,b}(X\times^GE_jG, w_j).
\]
2.  Given a motivic ring spectrum $\sE\in\SH(B)$, we call $\sE$ {\em bounded} if, for each $X\in \Sch^G_q/B$, each $w\in \sK^G(X)$ and each $(a,b)$, the pro-system
$\{\sE^\BM_{a,b}(X\times^GE_jG/B_jG, w_j)\}_j$
is eventually constant.
\end{definition}
\begin{remark} Taking $B=\Spec k$,  $k$ a perfect field, if $\sE=\EM(\sM_*)$ for $\sM_*$ a homotopy module \cite[Definition 5.2.4]{MorelICTP}, the description of $\EM(\sM_*)^\BM(-,-)$ using the Rost-Schmid complex implies that $\EM(\sM_*)$ is bounded. See \cite{DiLorenzoMantovani}, or \cite[\S 4]{LevineAtiyahBott} for details of the proof.
\end{remark}

\begin{remark}\label{rem:BMEquivHomOps}  As we have discussed in \cite[Proposition 4.11, Remark 4.12]{LevineAtiyahBott},  the operations of refined lci pull-back and proper push-forward for Borel-Moore homology, and the usual pull-back for cohomology,  as well as the corresponding identities relating these operations for twisted $\sE$-Borel-Moore homology and twisted $\sE$-cohomology described in the previous sections, extend to the equivariant setting defined using the algebraic Borel construction described above; essentially, one just applies the functor $X\mapsto X\times^GE_jG$ for $j\gg 0$, and then check that the operation is compatible with the transition maps going from $j+1$ to $j$ in the pro-system.   See  Proposition~\ref{prop:EquivPushforwardPullback} for central argument and its application in the cases of proper pushforward and refined lci pullback.  We also recall the extension of the various product structures below.
 
The exact localization sequence for a closed immersion $i\:Z\to X$ with open complement $j\:X\setminus U\to X$ gives a corresponding exact localization sequence for $Z\times^GE_jG\to X\times^GE_jG$ and its open complement $(X\setminus Z)\times^GE_jG\to X\times^GE_jG$. In case $\sE$ is bounded, these yield an exact localization sequence on the limit.  

However, if the given theory is not bounded, the possible lack of an exact localization sequence is a serious drawback. Similarly, one cannot use the homotopical approach of Morel-Voevodsky \cite[\S4]{MorelVoevodsky} to show the independence of the $G$-equivariant theory on the choice of embedding $G\hookrightarrow \GL_n$; instead we need to use (in \S\ref{subsec:Independence}) the method of Totaro \cite{Totaro} and Edidin-Graham \cite{EGEquivIntThy} . This has been rectified in the work of D'Angelo \cite{DAngelo}, which replaces the limit of $\sE$-valued Borel-Moore homology with the $\sE$-valued cohomology of a ``Borel-Moore motive'' formed by taking an appropriate colimit of  the objects $\pi_{X,G,j!}(\Sigma^{w_j}1_{X\times^GE_jG})$. This preserves all the exactness properties one would like to have, and agrees with the naive theory presented here in the case of a bounded theory. 
 \end{remark}

\begin{proposition} \label{prop:EquivPushforwardPullback}
1.  Let $f:Y\to X$ be a proper morphism in $\Sch^G_q/B$ and take $v\in \sK^G(X)$.  Then for each $j$, the map $f\times^G\id_{E_jG}:Y\times^G{E_jG}\to
X\times^G{E_jG}$ is proper, and the maps $(f\times^G\id_{E_jG})_*:\sE^\BM(Y\times^G{E_jG}/B_jG, v_j)\to \sE^\BM(X\times^G{E_jG}/B_jG, v_j)$ define a map of pro-systems
$\{(f\times^G\id_{E_jG})_*\}$, giving a well-defined map on the limit
\[
f_*:\sE^\BM_G(Y/BG,v)\to \sE^\BM(X/BG,v).
\]
2. Let
\[
\xymatrix{
X_0\ar[r]^{f'}\ar[d]_q&X\ar[d]^p\\
Y_0\ar[r]^f&Y
}
\]
be a cartesian square 
in $\Sch^G/B$,   and take $v\in \sK^G(Y)$.   Then for each $j$, the map
\[
f\times^G\id_{E_jG}\:Y_0\times^GE_jG\to Y\times^GE_jG
\]
is  lci,   with $L_{f\times^G\id_{E_jG}}\cong(L_f)_j$, the diagram
\[
\xymatrix{
X_0\times^GE_jG\ar[r]^{f'\times^G\id}\ar[d]_{q\times^G\id}&X\times^GE_jG\ar[d]^{p\times^G\id}\\
Y_0\times^GE_jG\ar[r]^{f\times^G\id}&Y\times^GE_jG
}
\]
is cartesian, and the sequence of refined lci pullback maps
\[
\sE^\BM_{**}(X\times^GE_jG/B_jG, v_j)
\xrightarrow{(f\times^G\id_{E_jG})^!_{p\times^G\id}}
\sE^\BM_{**}(X_0\times^GE_jG/B_jG, v_j-(\sN_\iota)_j)
\]
gives rise to a well-defined map of pro-systems
\[
\{\sE^\BM_{**}(X\times^GE_jG/B_jG, v_j)\}_j\to \{\sE^\BM_{**}(X_0\times^GE_jG/B_jG, v_j\-(L_j)\}_j.
\]
This in turn gives a well-defined map
\[
f^!_p\:\sE^\BM_{G,\ **}(X/BG, v)\to \sE^\BM_{G,\ **}(X_0/BG, v+L_f).
\]
\end{proposition}

\begin{proof} For (1),  the properness of $f\times^G\id_{E_jG}$ is easy to see.  To show that the maps $(f\times^G\id_{E_jG})_*$ give rise to a well-defined map of pro-systems, 
we need show that that
\[
(f\times^G\id_{E_{j+1}G})_*\circ i_{X, G,j}^!= i_{X, G,j}^!\circ(f\times^G\id_{E_jG})_*
\]
Factoring the transition maps $i_{X,G,j}^!$ as $i_{X,G,j}^!=i^{0j}_{X,G,j}\circ\eta_{X,G,j}^!=
 (p_{X,G,j}^!)^{-1}\circ\eta_{X,G,j}^!$, where $\eta_{X,G,j}$ is smooth and $p_{X,G,j}$ is the (smooth) projection for a vector bundle,  we see that this follows from the compatibility of proper pushforward with smooth pullback (Proposition~\ref{prop:PushPullComp})
 
 For (2), the first two assertions are easily verified and are left to the reader. To show that the maps $(f\times^G\id_{E_jG})^!_{p\times^G\id}$ give rise to a well-defined map of pro-systems, 
we need show that that
\[
(f\times^G\id_{E_jG})^!_{p\times^G\id}\circ i_{X, G,j}^!= i_{X, G,j}^!\circ(f\times^G\id_{E_{j+1}G})^!_{p\times^G\id}
\]
as maps
\[
\sE^\BM_{**}(X\times^GE_{j+1}G/B, v_{j+1}+\Omega_{B_{j+1}G/B})
\to \sE^\BM_{**}(X_0\times^GE_jG/B, v_j+\Omega_{B_jG/B}-(\sN_\iota)_j).
\]
Factoring the transition maps $i_{X,G,j}^!$ as in (1),  we see that this follows from the compatibility of refined lci pullback with smooth pullback (Proposition~\ref{prop:LciProperties}(1b)).
 \end{proof}

\begin{remark}\label{rem:EquivariantExtension} 1. As noted above, once we have a well-defined proper pushforward and  refined lci pullback on equivariant Borel-Moore homology, all the relations satisfied by proper pushforward and refined lci (respectively, Gysin, respectively smooth) pullback for Borel-Moore homology, as in Proposition~\ref{prop:PushPullComp}, Proposition~\ref{prop:LciProperties},  and Corollary~\ref{cor:commutativity}, are also satisfied in the equivariant setting. The same holds for the Poincar\'e duality isomorphism \eqref{eqn:PDIso}. For a bounded theory, we also have an equivariant localization sequence, with the same compatibilities with the smooth pullback as expressed in Proposition~\ref{prop:LocSmoothComp}.
\\[5pt]
2. The argument used for Proposition~\ref{prop:EquivPushforwardPullback} also shows that cup products, cap products, external cup and cap products and refined intersection products also extend to the equivariant setting, and satisfy the same relations as stated in Lemma~\ref{lem:ProductComp}. \\[5pt]
3. Let $p_X:X\to B\in \Sch^G_q/B$ be an lci scheme over $B$.  Let $[B]_{G,\sE}\in \sE^\BM_G(B/B)\cong \sE^{0,0}_G(B)$ correspond to the unit and let $[X]_{G,\sE}=p_X^![B]_{G,\sE}\in \sE^\BM_G(X/BG,L_{X/k})$. It is easy to see that $[X]_{G,\sE}$ is given by the family of fundamental classes $[X\times^GE_jG/B_jG]_\sE\in \sE^\BM(X\times^GE_jG/B_jG,L_{X\times^GE_jG/B_jG})$, defined by considering $B_jG$ as the base-scheme.  We have the fundamental class $\eta_{G,f}\in \sE^\BM_{G_{Y_0}}(Y/B_jG_{Y_0}, L_{Y/Y_0})$ for  $f:Y\to Y_0$ an lci morphism in $\Sch^G_q/Y_0$, defined similarly, where $G_{Y_0}$ is the group-scheme $G\times_BY_0$ over $Y_0$.
\\[5pt]
4. Vistoli's lemma Proposition~\ref{prop:SpecializationComm} and  all of its consequences as detailed in \S\ref{sec:Vistoli} similarly extend to the equivariant setting.
\end{remark}

\begin{remark}[Equivariant Euler classes]\label{rem:EquivEuleClass} Let $\sE\in \SH(B)$ be a commutative ring spectrum. For $X\in \Sch^G/B$, let $V\to X$ be a $G$-equivariant vector bundle. By applying descent to the vector bundle $p_1^*V\to X\times_BE_jG$, we have the vector bundle $V_j\to X\times^GE_jG$ for every $j$, giving the Euler class
\[
 e^\sE(V_j)\in \sE^{0,0}(X\times_BE_jG, \sV_j)
 \]
 where $\sV_j$ is the sheaf of sections of $V_j$ (see \S\ref{sec:FundamentalClasses}). 
 
The maps  $i_{X,G,j}^*$ commute with pullback  and cup products in cohomology. Thus 
 \[
 i_{X,G,j}^*(e^\sE(V_{j+1}))=e^\sE(V_j), 
 \]
so  the family $\{e^\sE(V_j)\in \sE^{0,0}(X\times_BE_jG, \sV_j)\}_j$ defines the equivariant Euler class $e^\sE_G(V)\in \sE^{0,0}_G(X, \sV)$.
 
 If $X$ is lci, then taking $B_jG$ as our base-scheme, we have  the  Borel-Moore Euler class
 \[
 e^\BM_\sE(V_j/B_jG)\in \sE^\BM(X\times_BE_jG/B_jG, L_j-\sV_j)
 \]
 where $L_j=L_{X\times_BE_jG/B_jG}$.  As we have seen in the proof of Proposition~\ref{prop:EquivPushforwardPullback}, the transition maps  $i_{X,G,j}^!$ commute with proper pushforward and refined lci pullback. Since
 \[
 e^\BM_\sE(V_j/B_jG)=s_j^!s_{j*}([X\times_BE_jG/B_jG/B_jG]_\sE)=e^\sE(V_j)\cap [X\times_BE_jG/B_jG/B_jG]_\sE,
 \]
 where $s_j$ is the zero-section for $V_j$, it follows that the family $\{e^\BM_\sE(V_j/B_jG)\}_j$ defines an element $e^\BM_{G,\sE}(V)\in \sE^\BM_G(X/BG, \sV)$, satisfying
 \begin{equation}\label{eqn:EquivEulerClassFundClass}
 e^\BM_{G,\sE}(V)=e^\sE_G(V)\cap_X[X]_{G,\sE}=s^!s_*([X]_{G,\sE}),
 \end{equation}
 where $s:X\to V$ is the zero-section.

If now $\sE$ is $\SL$-oriented,  $X$ is lci and proper of pure dimension $d=\rnk V$ over $B$, and $V$ has a $G$-equivariant relative orientation (Definition~\ref{def:(Rel)Orientation}) 
\[
\rho:\det L_{X/B}\otimes \det\sV^\vee \xrightarrow{\sim} \sL^{\otimes 2}
\]
for some $G$-linearized invertible sheaf $\sL$ on $X$, then we may apply $\deg^{\sE,G}_B:=p_{X*}$ to 
 \begin{multline*}
  e^\BM_{G,\sE}(V)\in \sE^\BM_G(X/BG, L_{X/B}-\sV)\\\cong 
   \sE^\BM_G(X, \det L_{X/B}-\det \sV){\xymatrix{\ar[r]^\rho_\sim&}}\sE^\BM_G(X/BG)
 \end{multline*}
 giving 
 \[
 \deg^{\sE,G}_B(e^\BM_{G,\sE}(V))\in \sE^\BM(B/BG)=\sE^{0,0}_G(B).
 \]
 \end{remark} 

\subsection{Independence of the embedding}\label{subsec:Independence}
Up to now, we have fixed an embedding $G\hookrightarrow \GL_n/B$ and defined all our equivariant theories with respect to this choice. In the approach used in \cite[\S 4.2]{MorelVoevodsky}, the independence of choices is proven by showing that the colimit in $\Spc(B)$ of ind-system $\{B_jG\}_j$ canonically represents  the \'etale classifying space of $G$ in $\sH(B)$, and is thus independent of any choices made in the construction. This approach has been used in \cite{KhanRavi} and \cite{DAngelo} for such independence results, but we use here instead an approach following  \cite{EGEquivIntThy, Totaro}

  Let $\iota:G\hookrightarrow \GL_n/B$ be an embedding  of $G$ as a closed subgroup scheme of $\GL_n$ over $B$. As we will be keeping track of the embedding, we write $E_{j,\iota}G$ for the open subscheme of $n\times n+j$ matrices, with $G$ action induced by $\iota$ and the usual action of $\GL_n$. We similarly write $B_{j,\iota}G$ for the quotient $B\times^G E_{j,\iota}G$, and $\Sch^G_{q; \iota}/G$ for $\Sch^G_q/G$. Given a finite set $S$ of embeddings $\iota_j:G\to \GL_{n_j}$, $j\in S$, we set $\Sch^G_{q; \iota_S}/G=\cap_j\Sch^G_{q;\iota_j}/B$, the intersection taking place in $\Sch^G/B$. 

\begin{lemma}\label{lem:Independence} Take two  embeddings $\iota_i:G\to \GL_{n_i}/B$, $i=1,2$.   Take $X$ in $\Sch^G_{q; \iota_{\{1,2\}}}/B$.  Let $\sE\in \SH(B)$ be bounded and take $v\in \sK^G(X)$.  Let $v_{j,i}$ be the perfect complex on  $X\times^GE_{j,\iota_i}G$ given as before by descent via the $G$-linearization on $p_1^*v$. Then given $a,b\in \Z$,  there is an integer $j_0$ and canonical isomorphisms of pro-systems
\[
\{\sE^\BM_{a,b}(X\times^GE_{j,\iota_1}G/B_{j, \iota_1}, v_{j,1})\}_{j\ge j_0}
\xymatrix{\ar[r]^{\{\psi_{j,21}^v\}_{j\ge j_0}}_\sim&} \{\sE^\BM_{a,b}(X\times^GE_{j,\iota_2}G/B_{j, \iota_2}, v_{j,2})\}_{j\ge j_0},
\]
inducing a canonical isomorphism
\[
\lim_j\sE^\BM_{**}(X\times^GE_{j,\iota_1}G/B_{j, \iota_1}, v_{j,1})\xymatrix{\ar[r]^{\psi_{21}^v}_\sim&} \lim_j \sE^\BM_{**}(X\times^GE_{j,\iota_2}G/B_{j, \iota_2}, v_{j,2}).
\]
Here $j_0$ is a constant depending on $a,b$, $\sE$, $X$,  and $v$.
\end{lemma}

\begin{proof}  The proof is a direct adaption of the argument used by to prove \cite[Theorem 1.1]{Totaro} or \cite[Definition-Proposition 1]{EGEquivIntThy}, both adapted from arguments of Bogomolov \cite{Bogomolov}.

We set  $X_{j,12}:=X\times^G(E_{j+n_2,\iota_1}G\times E_{j+n_1,\iota_2}G)$, and $B_{j,12}:=   B\times^G(E_{j+n_2,\iota_1}G\times E_{j+n_1,\iota_2}G)$.

Fix $a,b\in \Z$. The $G$-scheme $E_{j+n_2,\iota_1}G$ is an open $G$-subscheme of the $n_1\times (n_1+n_2+j)$ matrices $M_{n_1, n_1+n_2+j}$, with closed complement of codimension $n_1+n_2+j+1$, and similarly $E_{j+n_1,\iota_2}G$ is an open $G$-subscheme of   $M_{n_2, n_1+n_2+j}$ with closed complement of  the same codimension $n_1+n_2+j+1$. On the other hand, using descent for $G$-linearized vector bundles, the quotient $G\backslash (X\times E_{j+n_2,\iota_1}G\times M_{n_2, n_1+n_2+j})$ exists and the projection $\tilde{q}_{j,1}:G\backslash (X\times E_{j+n_2,\iota_1}G\times M_{n_2, n_1+n_2+j})\to  X\times^G E_{j+n_2,\iota_1}G$ makes $G\backslash (X\times E_{j+n_2,\iota_1}G\times M_{n_2, n_1+n_2+j})$ a vector bundle over $X\times^G E_{j+n_2,\iota_1}G$.

Thus, letting $\tilde{X}_{j,1}:= X\times^G(E_{j+n_2,\iota_1}G\times M_{n_2, n_1+n_2+j})$, $\tilde{B}_{j,1}G:=B\times^G(E_{j+n_2,\iota_1}G\times M_{n_2, n_1+n_2+j})$, and letting $q_{j,1}:\tilde{B}_{j,1}G\to B_{j+n_2, \iota_1}G$ be the projection,  the smooth pullback
\[
q_{j,1}^!:\sE^\BM_{a,b}(X\times^G E_{j+n_2,\iota_1}G/B_{j+n_2, \iota_1}G, v_{j,1})\to \sE^\BM_{a,b}(\tilde{X}_{j,1}/\tilde{B}_{j,1}G, v_{j,1})=
\]
is an isomorphism. On the other hand,  $X_{j,12}$ is an open subscheme of $\tilde{X}_{j,1}$ and has closed complement of codimension  $n_1+n_2+j+1$. Since $\sE$ is bounded, the restriction map 
\[
\sE^\BM_{a,b}(\tilde{X}_{j,1}/\tilde{B}_{j,1}G, v_{j,1})\to
\sE^\BM_{a,b}(X_{j,12}/
B_{j,12}, v_{j,12})
\] 
is an isomorphism for  $j\gg0$, where $v_{j,12}in \sK(X_{j,12})$ is given by appying descent to $p_1^*v\in \sK(X\times_B(E_{j+n_2,\iota_1}G\times E_{j+n_1,\iota_2}G)$.  In other words, letting $p_{j,1}:B_{j,12}G\to B_{j+n_2, \iota_1}G$ be the projection,  the smooth pullback
\[
p_{j,1}^!:\sE^\BM_{a,b}(X\times^G E_{j+n_2,\iota_1}G/B_{j+n_2}G, v_{j,1})\to \sE^\BM_{a,b}(X_{j,12}/B_{j,12}G, v_{j,12})
\]
is an isomorphism for $j\gg0$. 

As above,  letting $p_{j,2}:B_{j,12}\to B_{j+n_1,\iota_2}G$ be the projection,  the smooth pullback
\[
p_{j,2}^!:\sE^\BM_{a,b}(X\times^G E_{j+n_1,\iota_2}G/B_{j+n_1,\iota_2}G, v_{j,2})\to \sE^\BM_{a,b}(X_{j,12}/B_{j,12}, v_{j,12})
\]
is an isomorphism for $j\gg0$. Since smooth pullback commutes with the refined Gysin pullback maps $i_{G,j}^!$, the families $\{p_{j,1}^!\}_j$ and $\{p_{j,2}^!\}_j$ define maps of the corresponding pro-systems,  so there is a $j_0$ such that $\{(p_{j,2}^!)^{-1}p_{j,1}^!\}_{j\ge j_0}$ gives an isomorphism of pro-systems 
\begin{multline*}
\{\sE^\BM_{a,b}(X\times^GE_{j,\iota_1}G/B_{j,\iota_1}G, v_{j,\iota_1})\}_{j\ge j_0}\\{\xymatrix{\ar[r]^{\{(p_{j,2}^!)^{-1}p_{j,1}^!\}_{j\ge j_0}}_\sim&}} \{\sE^{a,b}(X\times^GE_{j,\iota_2}G/B_{j,\iota_2}G, v_{j,\iota_2})\}_{j\ge j_0}
\end{multline*}
\end{proof}

For later use, we prove a similar indepence of choices for a product $G=G_1\times_BG_2$. Given embeddings $\iota_j:G_i\to \GL_{n_i}/B$ as closed subgroup schemes (always smooth over $B$), $i=1,2$, we the product embedding $\iota_1\times \iota_2:G_1\times_BG_2\to \GL_{n_1}\times_B\GL_{n_2}$, which we may then compose with the usual embedding $i_{n_1, n_2}:GL_{n_1}\times_B\GL_{n_2}\to GL_{n_1+n_2}$ by putting $\GL_{n_1}$ in the upper left $n_1\times n_1$ block and $\GL_{n_2}$ in the lower right block, giving the embedding $\iota_{12}:G_1\times G_2\to \GL_{n_1`+n_2}$. Given $X\in \Sch^{G_1\times G_2}/B$ we assume that the quotients $X\times^{G_1\times G_2}E_{j+n_2, \iota_1}G_1\times 
E_{j+n_1, \iota_2}G_2$, $X\times^{G_1\times G_2}E_{j,\iota_{12}}(G_1\times G_2)$,  $B_{j,\iota_1}G_1$, $B_{j,\iota_2}$ and $B_{j,\iota_{12}}(G_1\times G_2)$ all exist for all $j\ge0$. 

The isomorphism $M_{n_1, n_1+n_2+j}\times M_{n_2\times n_1+n_2+j}=M_{n_1+n_2, n_1+n_2+j}$ gives the open immersion $j:E_{n_2+j,\iota_1}G_1\times E_{n_1+j, \iota_2}G_2\to E_{j,\iota_{12}}G_1\times G_2$ and   gives us the cartesian square of quotients
\[
\xymatrix{
X\times^{G_1\times G_2}E_{j+n_2, \iota_1}G_1\times 
E_{j+n_1, \iota_2}G_2\ar[r]^-{\tilde{\alpha}_j}\ar[d]&X\times^{G_1\times G_2}E_{j,\iota_{12}}G_1\times G_2\ar[d]\\
B_{n_2+j,\iota_1}G_1\times B_{n_1+j,\iota_2}G_2\ar[r]^-{\alpha_j}&B_{j,\iota_{12}}G_1\times G_2
}
\]

\begin{lemma}\label{lem:Independence2} Take $X$ as above and take a $v\in \sK^{G_1\times G_2}(X)$. Let $\sE\in \SH(B)$ be a bounded commutative ring spectrum, and take $a,b\in \Z$. \\[2pt]
 1.  The map $\alpha_j$ is an open immersion. \\[2pt]
2.  The relative pullback map
\begin{multline*}
\sE^\BM_{a,b}(X\times^{G_1\times G_2}E_{j,\iota_{12}}\GL_{n_1+n_2}/B_{j,\iota_{12}}(G_1\times G_2), v_j)\\
\xrightarrow{\alpha_j^!}
\sE^\BM_{a,b}(X\times^{G_1\times G_2}E_{j+n_2,\iota_1}G_1\times 
E_{j+n_1,\iota_2}G_2/B_{n_2+j,\iota_1}G_1\times B_{n_1+j,\iota_2}G_2, v_j)
\end{multline*}
is an isomorphism for all $j\gg0$.
\end{lemma}

\begin{proof} (1) follows from the fact that $E_{j+n_2,\iota_1}G_1\times 
E_{j+n_1,\iota_2}G_2\to E_{j,\iota_{12}}(G_1\times G_2)$ is a $G_1\times G_2$ equivariant open immersion. Moreover the complement has codimension $\ge \min(n_1+j+1, n_2+j+1)$ so the same estimate holds for the induced open immersions $\alpha_j$ and $\tilde{\alpha}_j$ on the quotients by a free $G_1\times G_2$-action. Using the fact that $\sE$ is bounded, (2) follows.
\end{proof}

\begin{remark}\label{rem:Independence3} 1. If $X$ is in $\Sm^G/B$, the Poincar\'e duality isomorphism transports Lemmas~\ref{lem:Independence} and \ref{lem:Independence2} to the equivariant cohomology $\sE^{*,*}_G(X,v)$.

For general $X\in \Sch^G/B$, we are not able to prove the independence  of choice of embedding for  $\sE^{*,*}_G(X,v)$. In Proposition~\ref{prop:CohIndependence}  below, we have a partial result in this direction.\\[2pt]
2. Using the arguments of the proofs of Lemmas~\ref{lem:Independence} and \ref{lem:Independence2}, one shows that for $\sE$ bounded, all the structures for   equivariant $\sE$-Borel-Moore homology discussed in Remark~\ref{rem:BMEquivHomOps}  are compatible with the isomorphism of pro-systems constructed in those Lemmas. For example, given two embedding $\iota_1:G\to \GL_{n_1}/B$, $\iota_2:G\to \GL_{n_2}/B$,  a proper map $f:Y\to X$ in $\Sch^G_q/B$ and a $v\in \sK^G(X)$, we have the corresponding maps of pro-systems 
\[
\{f^1_{j*}:\sE^\BM_{a,b}(Y\times^GE_{j,\iota_1}G/B_{j,\iota_1}G, v_{j,1})\to \sE^\BM_{a,b}(X\times^GE_{j,\iota_1}G/B_{j,\iota_1}G, v_{j,1})\}_j
\]
and 
\[
\{f^2_{j*}:\sE^\BM_{a,b}(Y\times^GE_{j,\iota_2}G/B_{j,\iota_2}G, v_{j,2})\to \sE^\BM_{a,b}(X\times^GE_{j,\iota_2}G/B_{j,\iota_2}G, v_{j,2})\}_j
\]
Then letting 
\[
\psi_{Y,j}:\sE^\BM_{a,b}(Y\times^GE_{j,\iota_1}G/B_{j,\iota_1}G, v_{j,1})\to \sE^\BM_{a,b}(Y\times^GE_{j,\iota_2}G/B_{j,\iota_2}G, v_{j,2})
\]
 and  
 \[
 \psi_{X,j}:\sE^\BM_{a,b}(X\times^GE_{j,\iota_1}G/B_{j,\iota_1}G, v_{j,1})\to 
 \sE^\BM_{a,b}(Y\times^GE_{j,\iota_2}G/B_{j,\iota_2}G, v_{j,2})
 \]
 be the isomorphisms given by the proof of Lemma~\ref{lem:Independence} for $j\gg0$, the diagram
\[
\xymatrix{
\sE^\BM_{a,b}(Y\times^GE_{j,\iota_1}G/B_{j,\iota_1}G, v_{j,1})\ar[r]^{\psi_{Y,j}}\ar[d]^{f^1_*}& \sE^\BM_{a,b}(Y\times^GE_{j,\iota_2}G/B_{j,\iota_2}G, v_{j,2})\ar[d]^{f^2_*}\\
\sE^\BM_{a,b}(X\times^GE_{j,\iota_1}G/B_{j,\iota_1}G, v_{j,1})\ar[r]^{\psi_{X,j}}&\sE^\BM_{a,b}(X\times^GE_{j,\iota_2}G/B_{j,\iota_2}G, v_{j,2})
}
\]
commutes for all $j\gg0$. Indeed, the isomorphisms $\psi_{Y,j}$ and $\psi_{X,j}$ are defined as a zig-zag of relative pullback maps in cartesian squares for smooth morphisms, which by Remark~\ref{rem:RelPullback} are all compatible with the proper pushforward maps induced by $f$. 

The proof for lci pullback in $\sE$-Borel-Moore homology, refined Gysin maps,   the various purity isomorphisms (for smooth $X$),  refined intersection products,  cap products, and localization sequences, are all the same as above, using that these structures are all compatible with  relative pullback by smooth morphisms,  in cartesian squares with the base a smooth morphism.  We show in Proposition~\ref{prop:CohIndependence} below that  cap product with the cohomological $G$-equivariant  Euler class of a $G$-equivariant vector bundle is also independent of the choice of embedding. 
\\[2pt]
3. It follows directly from the construction of the isomorphism in the proof of Lemma~\ref{lem:Independence} that, given a third closed embedding  $\iota_3:G\to \GL_{n_3}/B$,  we have $\psi_{31}^v=\psi_{32}^v\circ\psi_{21}^v$.
Thus, we have a construction of $\sE^\BM_{G, **}(X/BG,v)$ that is independent of the choice of embedding $\iota:G\to \GL_n/B$, up to a unique isomorphism. In particular, all the structures discussed in (2) are well-defined on the canonically defined objects  $\sE^\BM_{G, **}(-/BG,-)$.
\end{remark}

\begin{proposition}\label{prop:CohIndependence} We retain the notation from Lemma~\ref{lem:Independence}, in particular, let
$\iota_i:G\to \GL_{n+i}/B$, $i=1,2$,  $X\in\Sch^G/B$,  $v\in \sK^G(X)$ and $\sE\in \SH(B)$ be as in Lemma~\ref{lem:Independence}. Let $\{\alpha_{j,1}\in \sE^\BM_{a,b}(X\times^GE_{j,\iota_1}G/B_{j, \iota_1}, v_{j,1})\}_{j\ge j_0}$, $\{\alpha_{j,2}\in \sE^\BM_{a,b}(X\times^GE_{j,\iota_1}G/B_{j, \iota_2}, v_{j,2})\}_{j\ge j_0}$ be families of elements such that $\alpha_{j,1}=
\psi^v_{j, 21}(\alpha_{j,1})$ for all $j\ge j_0$. 

Let $V\to X$ be a $G$-equivariant vector bundle, with sheaf of sections $\sV$. For $i=1,2$,  let $V_{j,i}$ be the vector bundle on  $X\times^GE_{j,\iota_i}$ defined from $p_1^*V\to X\times E_{j,\iota_i}$ by descent, and let $\sV_{j,i}$ be the sheaf of sections of $V_{j,i}$. We have the families of cohomological Euler classes
\[
e^\sE(V_{j,i})\in \sE(X\times^GE_{j,\iota_i}G, \sV_{j,i}),\ i=1,2
\]
defined in \S\ref{sec:FundamentalClasses}. Then for some $j_1\ge j_0$, we have
\[
\psi_{j, 21}^{v-\sV}(e^\sE(V_{j,1})\cap \alpha_{j,1})=e^\sE(V_{j,2})\cap \alpha_{j,2},\ j\ge j_1.
\]
\end{proposition}

\begin{proof} Looking to the proof of Lemmas~\ref{lem:Independence}, let 
\[
X_{j,12}:=X\times^G (E_{j+n_2,\iota_1}G\times E_{j+n_1,\iota_2}G),\ B_{j,12}G:=
B\times^G(E_{j+n_2,\iota_1}G\times E_{j+n_1,\iota_2}G),
\]
\[
X_{j,1}:=X\times^G E_{j+n_2,\iota_1}G,\ B_{j,1}G:= B\times^G E_{j+n_2,\iota_1}G,
\]
\[ X_{j,2}:=X\times^G E_{j+n_1,\iota_2}G,\ B_{j,2}G:=B\times^G E_{j+n_1,\iota_2}G
\]
giving the diagram
\[
\xymatrix{
&X_{j,12}\ar[dl]_{\tilde{p}_{j,1}}\ar[dr]^{\tilde{p}_{j,2}}\ar[rrr]&&&B_{j,12}\ar[dl]_{p_{j,1}}\ar[dr]^{p_{j,2}}\\
X_{j,1}\ar@/_15pt/[rrr]&&X_{j,2}\ar@/^15pt/[rrr]|(.45)\hole&B_{j,1}&&B_{j,2}
}
\] 
\ \\[5pt]
and the corresponding diagram of $\sE$-Borel-Moore homology
\[
\xymatrix{
& \sE^\BM_{a,b}(X_{j,12}/B_{j,12}, v_{j,12})\\
\sE^\BM_{a,b}(X_{j,1}/B_{j,1}, v_{j,1})\ar[ru]^{p_{j,1}^!}&&
\sE^\BM_{a,b}(X_{j,2}/B_{j,2}, v_{j,2})\ar[ul]_{p_{j,2}^!}
}
\]
yielding a  diagram of pro-systems. Moreover, both maps $p_{j,1}^!$, $p_{j,2}^!$ are isomorphisms on $\sE$-Borel-Moore homology for $j\ge j_0$.  

Let $V_{j,12}$ be the vector bundle on $X_{j,12}$ induced by $p_1^*V\to X\times E_{j+n_2,\iota_1}G\times E_{j+n_1,\iota_2}G$ by descent, with sheaf of sections $\sV_{j,12}$. We  have the diagram of $\sE$-cohomology, 
\[
\xymatrix{
& \sE(X_{j,12}, \sV_{j,12})\\
\sE(X_{j,1}, \sV_{j,1})\ar[ru]^{p_{1,j}^*}&&
\sE(X_{j,2},\sV_{j,2})\ar[ul]_{p_{2,j}^*}
}
\]
also giving a diagram of pro-systems.

We have the cohomological Euler class $e^\sE(V_{j,12})\in \sE(X_{j,12}, \sV_{j,12})$
with
\[
\tilde{p}_{j,1}^*(e^\sE(V_{j,1}))=e^\sE(V_{j,12})=\tilde{p}_{j,2}^*(e^\sE(V_{j,2})).
\]

For $j\ge j_0$, the assumption $\alpha_{j,2}=\psi_{j,21}^v(\alpha_{j,1}):=(p_{j,2}^!)^{-1}p_{j,1}^!(\alpha_{j,1})$ implies  $p_{j,1}^!(\alpha_{j,1})=p_{j,2}^!(\alpha_{j,2})$; define
\[
\alpha_{j,12}:=p_{j,1}^!(\alpha_{j,1})=p_{j,2}^!(\alpha_{j,2}),
\]
for $j\ge j_0$.  Then  it follows from   \eqref{eqn:ProductFunct} that  
\[
\tilde{p}_{j,1}^!(e^\sE(V_{j,1})\cap \alpha_{j,1})=e^\sE(V_{j,12})\cap \alpha_{j,12}=\tilde{p}_{j,2}^!(e^\sE(V_{j,2})\cap \alpha_{j,2})
\]
for $j\ge j_0$.  Letting $j_0'$ be the integer given by Lemma~\ref{lem:Independence} for $\sE, X, v-\sV$ and let $j_1=\max(j_0, j_0')$, we thus have
\[
\psi_{j,21}^{v-\sV}(e^\sE(V_{j,1})\cap \alpha_{j,1}):=((p_{j,2}^!)^{-1}p_{j,1}^!)((e^\sE(V_{j,1})\cap \alpha_{j,1}))=e^\sE(V_{j,2})\cap \alpha_{j,2}
\]
for $j\ge j_1$.
\end{proof}

\subsection{Change of group}\label{subsec:Change}
Suppose we have a homomorphism of smooth group schemes over $B$, $\rho:H\to G$. In this section we construct functorial and natural change of group maps $\rho^\#_{X,v}:\sE^\BM_{G,**}(X/BG, v)\to \sE^\BM_{H,**}(X/BH, v)$   and discuss compatibilities of $\rho^\#$ with the other structures we have defined on the equivariant Borel-Moore homology. 

To avoid technical problems having to do with the existence of quotients, in this section, we work over the base-scheme $B=\Spec k$, $k$ a field, so that $\Sch^G_q/B=\Sch^G/B$. We also assume throughout this section that $\sE$ is bounded

 We first discuss the case of a projection $p_2:H=G_1\times_B G_2\to G_2=G$, with $G_i$ as a closed subgroup scheme of  $\GL_{n_i}/B$, smooth over $B$, and take $X\in \Sch^{G_2}/B$. 
  
Take $v\in \sK^G(X)$.  This gives us the Tor-independent cartesian diagram
  \[
 \xymatrixcolsep{40pt}
 \xymatrix{
 X\times^{G_1\times_BG_2} E_jG_1\times_B E_jG_2\ar[d]^{\pi_{X, G_1\times G_2,j}}\ar[r]^{i_{X, G_1\times G_2,j}}&X\times^{G_1\times_BG_2} E_{j+1}G_1\times_B E_{j+1}G_2\ar[d]^{\pi_{X, G_1\times G_2,j+1}}\\
 B_jG_1\times_B B_jG_2\ar[r]^{i_{G_1\times G_2,j}}&B_{j+1}G_1\times_B B_{j+1}G_2
 }
 \]
 giving the pro-system
 \[
 \{\sE^\BM_{**}(X\times^{G_1\times_BG_2} E_jG_1\times_B E_jG_2/B_jG_1\times B_jG_2, v_j), i_{G_1\times G_2,j}^!\}_j.
 \]
By Lemma~\ref{lem:Independence2}, we may use this  pro-system to define
 $\sE^\BM_{G_1\times_B G_2, **}(X/BG_1\times BG_2, v)$.  
  
 We also have the Tor-independent cartesian diagram
 \begin{equation}\label{eqref:Diagp2}
 \xymatrixcolsep{40pt}
 \xymatrix{
 X\times^{G_1\times_BG_2} E_jG_1\times_B E_jG_2\ar[d]^{\pi_{X, G_1\times G_2,j}}\ar[r]^-{\id_X\times p_{2,j}}&X\times^{G_2} E_jG_2\ar[d]^{\pi_{X, G_2,j}}\\
 B_jG_1\times B_jG_2\ar[r]^-{p_{2,j}}&B_jG_2.
 }
 \end{equation}
 In particular, the horizontal maps are smooth. This gives the pullback maps of pro-systems
 \begin{multline*}
\{\sE^\BM_{**}(X\times^{G_2}E_jG_2/B_jG_2, v_j), i_{X,G_2,j}^!\}_j\\
\xrightarrow{\{ p_{2,j}^!\}_j}
 \{\sE^\BM_{**}(X\times^{G_1\times_BG_2} E_jG_1\times_B E_jG_2/B_jG_1\times B_jG_2, v_j), i_{X,G_1\times G_2,j}^!\}_j
\end{multline*}
 This  induces the  change of group map
 \[ p_2^*:\sE^\BM_{G_2, **}(X/BG_2, v)\to
 \sE^\BM_{G_1\times_B G_2, **}(X/B(G_1\times G_2), v).
 \]
 
 Now suppose that $H\subset G$ is a closed subgroup scheme of $G$, smooth over $B$, and we are given $G$ as a smooth closed subsgroup scheme of $\GL_n/B$. This also realizes $H$ as  a smooth closed subsgroup scheme of $\GL_n/B$, and gives us the pro-system 
 $\{\sE^\BM_{**}(X\times^GE_jG/B_jG, v_j), i_{X,G,j}^!\}_j$ defining   $\sE^\BM_{G, **}(X/BG, v)$ as the   limit, as well as the corresponding pro-systems for $H$, defining  $\sE^\BM_{H, **}(X/BH, v)$ as the  limit. As above, we have the  diagram
 \begin{equation}\label{eqn:DiagI}
 \xymatrixcolsep{40pt}
 \xymatrix{
 X\times^H E_jH\ar[d]^{\pi_{X, H,j}}\ar[r]^-{p_{X, G,H,j}}&X\times^{G} E_jG\ar[d]^{\pi_{X, G_2,j}}\\
 B_jH\ar[r]^-{p_{G,H,j}}&B_jG.
 }
 \end{equation}
 Recalling that $EH_j=EG_j$ as schemes, where we just restrict the $G$ action on $EG_j$ to $H$ to define $EH_j$, and similarly for $B_jH:=H\backslash E_jH\to G\backslash E_jG$,  we see that the maps $p_{X, G,H,j}$ and $p_{G,H,j}$ are smooth and the diagram is cartesian. 
 
 As before, this gives us the  pullback maps of pro-systems
\begin{multline*}
\{\sE^\BM_{**}(X\times^GE_jG/B_jG, v_j), i_{X,G,j}^!\}_j\\\xrightarrow{\{ p_{X, G,H,j}^!\}_j}
 \{\sE^\BM_{**}(X\times^H E_jH/B_jH, v_j), i_{X,H,j}^!\}_j
 \end{multline*}
 inducing the change of group map
 \[
 i^*:\sE^\BM_{G, **}(X/BG, v)\to \sE^\BM_{H, **}(X/BH, v).
 \]
 
For $\rho:H\to G$ an arbitrary homomorphism of smooth group schemes over $B$, we factor $\rho$ as the graph $i_\rho:H\to H\times_BG$ followed by the projection $p_2:H\times_BG\to G$ and define $\rho^\#:=i_f^*p_2^*$.

\begin{proposition}\label{prop: ChangeOfGroupProperties} Take $B=\Spec k$, $k$ a field, and let $\sE\in\SH(B)$ be a bounded commutative ring spectrum. \\[2pt]
 1. Given a homomorphism $\rho:H\to G$, the map $\rho^\#$  is independent of the choices of closed embeddings $H\hookrightarrow \GL_{n_1}$, $G\hookrightarrow\GL_{n_2}$.\\[2pt]
 2. Given homomorphisms of smooth group schemes over $B$, $\rho_1:G_1\to G_2$, $\rho_2:G_2\to G_3$, then $\rho_1^\#\rho_2^\#=(\rho_2\rho_1)^\#$.\\[2pt]
 3. Given a homomorphism  $\rho:H\to G$ of smooth group schemes over $B$, a morphism  $q:X\to Y\in \Sch^G/B$ and $v\in \sK^G(Y)$, we have
 \begin{enumerate}
\item[(i)] If $q$ is lci, then $q^!\rho^\#=\rho^\#q^!:\sE^\BM_{G, **}(Y/BG, v)\to \sE^\BM_{H, **}(X/BH, v+L_q)$. 
\item[(ii)] If $q$ is proper, then $q_*\rho^\#=\rho^\#q_*:\sE^\BM_{G, **}(X/BG, v)\to \sE^\BM_{H, **}(Y/BH, v)$. 
 \end{enumerate}
 4.   $\rho^\#:\sE^\BM_{G,**}(X/BG,v)\to \sE^\BM_{H,**}(X/BH,v)$ is natural with respect to isomorphisms in $\sK(X)$ and with respect to morphisms of commutative ring spectra in $\SH(B)$. \\[2pt]
5. For $Y\in \Sch^G/B$, and $V\to Y$ a $G$-vector bundle, $\rho_\#$ is natural with respect to cap product $e^\sE_{-}(V)\cap_Y(-)$ and external cap product $e^\sE_{-}(V)\cap_{Y,X}-$; if $Y$ is in $\Sm^G/B$, $\rho_\#$ is natural with respect to cap product $(-)\cap_Y(-)$,  external cap product $(-)\cap_{Y,X}(-)$ for $X\in \Sch^G/B$, and refined intersection product $(-)\cdot_{f,g}(-)$ for maps $f:Y\to Y_0$ in $\Sm^G_q/B$, $g:X\to Y_0$ in $\Sch^G/B$.
\\[2pt]
6. Suppose we have a closed immersion $i:Z\to X$ in $\Sch^G/B$ with open complement $j:U\to X$. Then the change of group maps $\rho^\#:\sE^\BM_{G, **}(?/BG, v)\to \sE^\BM_{H, **}(?/BH, v)$ for $?=X, Z, U$ map the long exact localization sequence
\[
\ldots\xrightarrow{i_*}\sE^\BM_{G, n,*}(X/BG, v)\xrightarrow{j^*}\sE^\BM_{G, n,*}(U/BG, v)
\xrightarrow{\del} \sE^\BM_{G, n+1,*}(Z/BG, v)\xrightarrow{i_*}\ldots
\]
to the corresponding long exact sequence in $\sE^\BM_{H, **}(?/BH, v)$, forming a commutative diagram.\\[2pt]
7. Let $f:Y\to Y_0$ be an lci morphism in $\Sch^G/B$, giving respective fundamental classes $\eta^\sE_{G,f}$, $\eta^\sE_{H,f}$. Then $\rho^\#(\eta^\sE_{G,f})=\eta^\sE_{H,f}$. In particular, for $X\in \Sch^G/B$ an lci scheme over $B$, we have $\rho^\#([X]_{G,\sE})=[X]_{H,\sE}$ \\[2pt]
8. For $X\in \Sch^G/B$ an lci scheme over $B$, and for $V\to X$ a $G$-equivariant vector bundle with sheaf of sections $\sV$, we have the equivariant Borel-Moore Euler classes 
\[
e^\BM_{G,\sE}(V)\in \sE^\BM_G(X/BG, L_{X/B}-\sV),\ e^\BM_{H,\sE}(V)\in \sE^\BM_G(X/BH, L_{X/B}-\sV). 
\]
Then 
\[
\rho^\#(e^\BM_{G,\sE}(V))=e^\BM_{H,\sE}(V).
\]
\end{proposition}

\begin{proof}  (1) follows by arguing as in Remark~\ref{rem:Independence3}.

For (2), 
we consider the following commutative diagram of homomorphisms of group schemes
\[
\xymatrix{
&&&&G_1\times G_3\ar@/^20pt/[dddl]^{p_2}\ar[dl]_{i_{\rho_1}\times \id_{G_3}}\\
G_1\ar[r]_{i_{\rho_1}}\ar@/^20pt/[rrrru]^{i_{\rho_2\rho_1}}&G_1\times G_2\ar[d]^{p_2}\ar[rr]^{\id_{G_1}\times i_{\rho_2}}&&G_1\times G_2\times G_3\ar[d]^{p_{23}}\\
&G_2\ar[rr]^{i_{\rho_2}}&&G_2\times G_3\ar[d]^{p_2}\\
&&&G_3
}
\]

By (1),  we may fix closed embeddings $G_i\hookrightarrow \GL_{n_i}$. We note that the change of group map $i^*$ for a closed embedding $i:H\to G$
 \[
i^*:\sE^\BM_{G, **}(X/BG, v)\to \sE^\BM_{H, **}(X/BH, v)
\]
 is defined on the individual terms in the pro-system $\{\sE^\BM_{**}(X\times^GE_jG/B_jG,v_j\}_j$, 
 $\{\sE^\BM_{**}(X\times^HE_j/B_jH,v_j\}_j$ by relative pullback  $i_j^!$ in the cartesian square
 \[
 \xymatrix{
 X\times^HE_j\GL_n\ar[r]^{\tilde{i}_j}\ar[d]& X\times^GE_j\GL_n\ar[d]\\
 B\times^H E_j\GL_n\ar[r]^{i_j}&B\times^G E_j\GL_n\ ;
 }
 \]
we also note that the map $i_j$ is smooth.  Thus, the functoriality of smooth relative pullback shows that  
\[
i_{\rho_1}^*(\id_{G_1}\times i_{\rho_2})^*=((\id_{G_1}\times  i_{\rho_2})\circ i_{\rho_1})^*=((i_{\rho_1}\times \id_{G_3})\circ i_{\rho_2\rho_1})^*=i_{\rho_2\rho_1}^*(i_{\rho_1}\times \id_{G_3})^*
\]
and similarly
\[
p_{23}^*\circ p_2^*=(p_2\circ p_{23})^*
\]
in  $\sE$-Borel-Moore homology. Thus, we need to show that 
 \[
p_2^*\circ i_{\rho_2}^*=(\id_{G_1}\times i_{\rho_2})^*\circ p_{23}^*
\]
and
\[
p_2^*=(i_{\rho_1}\times \id_{G_3})^*\circ (p_2p_{23})^*.
\]

The argument is again along the lines of Remark~\ref{rem:Independence3}. Note that there one new issue in understanding the composition $(\id_{G_1}\times i_{\rho_2})^*\circ p_{23}^*$, in that the different compositions are defined using different types of embeddings. 

Concretely, we first choose embeddings $G_i\hookrightarrow \GL_{n_i}$, $i,2,3$, giving block diagonal embeddings
\[
G_1\times G_2\times G_3\hookrightarrow \GL_{n_1}\times \GL_{n_2}\times \GL_{n_3}\hookrightarrow \GL_{n_1+n_2+n_3}
\]
which we may use to define $(\id_{G_1}\times i_{\rho_2})^*$, and
\[
G_2\times G_3\hookrightarrow   \GL_{n_2}\times \GL_{n_3}\hookrightarrow \GL_{n_2+n_3}
\]
which we may use to define $i_{\rho_2}^*$. However, using Lemma~\ref{lem:Independence2}, we may instead use throughout the product embeddings and relative $\sE$-Borel-Moore homology of the quotients 
\begin{multline*}
X\times^{G_1\times G_2\times G_3} E_j\GL_{n_1}\times E_j\GL_{n_2}\times E_j\GL_{n_3}\\\to
G_1\times G_2\times G_3\backslash E_j\GL_{n_1}\times E_j\GL_{n_2}\times E_j\GL_{n_3}
\end{multline*}
and
\[
X\times^{G_2\times G_3} E_j\GL_{n_1}\times E_j\GL_{n_2}\times E_j\GL_{n_3}\to
 G_2\times G_3\backslash E_j\GL_{n_1}\times E_j\GL_{n_2}\times E_j\GL_{n_3}
\]
to define $(\id_{G_1}\times i_{\rho_2})^*$, as well as using the pair 
\[
X\times^{G_2\times G_3} E_j\GL_{n_1}\times E_j\GL_{n_2}\times E_j\GL_{n_3}\to
G_2\times G_3\backslash E_j\GL_{n_1}\times E_j\GL_{n_2}\times E_j\GL_{n_3}
\]
and
\[
X\times^{G_2} E_j\GL_{n_1}\times E_j\GL_{n_2}\times E_j\GL_{n_3}\to
G_2\backslash E_j\GL_{n_1}\times E_j\GL_{n_2}\times E_j\GL_{n_3}
\]
to define $i_{\rho_2}^*$. We may then use the same collection of pairs to define $p_2^*$ and $p_{23}^*$, which then fit together to  gives us a commutative square of relative pullback maps on $\sE$-Borel-Moorre homology, showing that 
$(\id_{G_1}\times i_{\rho_2})^*\circ p_{23}^*=p_2^*\circ i_{\rho_2}^*$. A similar arguments shows that $p_2^*=(i_{\rho_1}\times \id_{G_3})^*\circ (p_2p_{23})^*$.

 The essential point in the proof of (3) is that the change of group map arises as a  relative pullback in a cartesian square whose base morphism is smooth. Then one uses: (i) lci pullback commutes with smooth relative pullback, (ii) smooth relative pull back and proper pushforward commute in cartesian squares. For (4), we use that  smooth relative pullback on $\sE^\BM(-,v)$ is natural in $\sE$ and in $v$. 

For the cap products with respect to $e^\sE_{(-)}(V)$ in (5),  we argue as in the proof of Proposition~\ref{prop:CohIndependence}. Let $\sV_{G,j}$, $\sV_{H\times G, j}$ and $\sV_{H,j}$ be the locally free sheaves on $X\times^G E_jG$, $X\times^{H\times G}E_jH\times E_jG$ and $X\times^HE_jH$  defined by descent from the respective $p_1^*\sV$, giving respective vector bundles $V_{G,j}$, $V_{H\times G,j}$ and $V_{H,j}$. We have the maps of pro-systems
\begin{align*}
\{\sE^{**}(X\times^GE_jG, \sV_{G, j})\}_j&\xrightarrow{p_{2,j}^*}
\{\sE^{**}(X\times^{H\times G}E_jH\times E_jG, \sV_{H\times G, j})\}_j 
\\&\xrightarrow{i_{\rho,j}^*}\{\sE^{**}(X\times^HE_jH, \sV_{H, j})\}_j
\end{align*}
with $p_{2,j}^*e^\sE(V_{G,j})=e^\sE(V_{H\times G, j})$, $i_{\rho,j}^*e^\sE(V_{H\times G, j})=
e^\sE(V_{H, j})$. Then for $\alpha=\{\alpha_j\in \sE^\BM_{**}(X\times^GE_jG/B_jG,w_j)\}_j\in  \sE^\BM_{G, a,b}(X, w)$, we have
\begin{align*}
\rho^\#(e^\sE_G(V)\cap \alpha)&=\{i_{\rho,j}^!p_{2,j}^!(e^\sE(V_{G,j})\cap \alpha_j)\}_j\\
&=\{i_{\rho,j}^!(e^\sE(V_{H\times G,j})\cap p_{2,j}^!(\alpha_j))\}_j\\
&=\{e^\sE(V_{H,j})\cap i_{\rho,j}^!(p_{2,j}^!(\alpha_j))\}_j=e^\sE_H(V)\cap \rho^\#(\alpha).
\end{align*} 
 
 The proof of (5) for the products involving $Y\in \Sm^G/B$ follows by using the Poincar\'e duality isomorphism $\sE^{a,b}_G(Y, v)\cong \sE^\BM_{-a,-b}(Y, \Omega_{Y/k}-v)$  and then arguing as in the proof of (1), using the naturality of these products with respect to lci (hence smooth) pullback \eqref{eqn:ProductFunct}.
 
 For (6), the proof is similar: the localization sequence is natural with respect to smooth relative pullback (Lemma~\ref{lem:RelPull-back}(3)).  For (7), we first  claim that
\[
 \rho^\#(\eta^\sE_{G, \id_{Y_0}})=\eta^\sE_{H, \id_{Y_0}}
 \]
Indeed, for $?=G, H$, $\eta^\sE_{?, \id_{Y_0}}\in \sE^\BM_?(Y_0/Y_0)$ corresponds to the unit in $\sE^{0,0}_{?}(Y_0)$ via the Poincar\'e duality isomorphism, from which our claim follows easily. For $f:Y\to Y_0$ an lci morphism in $\Sch^G_q/B$, we have $\eta_{?,f}^\sE=f^!\eta^\sE_{?, \id_{Y_0}}$, $?=G,H$, so (7) follows from (3i). (8) follows from (3i,ii), (7) and  \eqref{eqn:EquivEulerClassFundClass}. \end{proof}

\section{Virtual fundamental classes}\label{sec:VirClass} In this section, we take $B=\Spec A$, for $A$ a noetherian ring of finite Krull dimension.  In \cite{LevineVirt} we showed how,  given a commutative ring spectrum $\sE\in \SH^G(B)$,  a $G$-linearized perfect obstruction theory $\phi_\bullet\:E_\bullet\to L_{Z/B}$ on some $Z\in \Sch^G/B$  gives rise to 
a virtual fundamental class in twisted Borel-Moore homology,  
\[
[Z,\phi_\bullet]^\vir_\sE\in
\sE^\BM(Z/B,G, E_\bullet). 
\]
In the next section, we define a virtual class $[Z,\phi_\bullet]^\vir_{\sE,G}$ in the $G$-equivariant (Totaro, Edidin-Graham) Borel-Moore homology $\sE^\BM_G(Z/BG, E_\bullet)$. In this section, we use the description of the Behrend-Fantechi virtual fundamental class given by Graber-Pandharipande in \cite{GP} to simplify our earlier construction of $[Z,\phi_\bullet]^\vir_\sE$, as an aid to defining $[Z,\phi_\bullet]^\vir_{\sE,G}$. We  fix a tame group-scheme $G$ over $B$.

\begin{remark} In the next section, we will use the Graber-Pandharipande construction to define an equivariant virtual fundamental class using the algebraic Borel construction of \S\ref{sec:TEGEquivCoh}. For this, we only use the results of this section for the trivial group.
\end{remark}

As above, we let $D^\perf_G(X)$ denote the derived category of $G$-linearized perfect complexes on $X$, and $D_G(X)$ the derived category of complexes of $G$-linearized quasi-coherent sheaves on $X$. For $X\in \Sch^G/B$, the cotangent complex $L_{X/B}$ has a canonical $G$-linearization, so defines an object $L_{X/B}$ in $D_G(X)$.

\begin{definition}
A $G$-linearized perfect obstruction theory on $X$ relative to $B$ is a $G$-linearized perfect complex $E_\bullet$ and a map $\phi\:E_\bullet\to L_{X/B}$ in $D_G(X)$ that defines a perfect obstruction theory after forgetting the $G$-linearizations. In other words, $E_\bullet$ is supported in (homological) degrees $[0,1]$, $h_0(\phi)$ is an isomorphism and $h_1(\phi)$ is a surjection. 
\end{definition}
As defined in  \cite{LevineVirt}, the virtual fundamental class associated to a perfect obstruction theory only depends on the truncation $\phi\:E_\bullet\to \tau_{\le1}L_{X/B}$  (we use the homological convention for the truncation); in what follows, we will abuse notation and refer to a map  $\phi\:E_\bullet\to \tau_{\le1}L_{X/B}$ in $D_G(X)$, with $E_\bullet\in D^\perf_G(X)$ supported in $[0,1]$, and with $h_0\phi$  an isomorphism  and $h_1\phi$ a surjection, as a $G$-linearized perfect obstruction theory relative to $B$. Throughout this section, all the obstruction theories will be relative to the fixed base-scheme $B$, so we will drop the phrase ``relative to $B$''.

Given a $G$-linearized perfect obstruction theory $\phi\:E_\bullet\to \tau_{\le1}L_{X/B}$ on some $X\in \Sch^G/B$, as $X$ is quasi-projective over $B$, $G$ is tame and $B$ is affine, $\phi$ has a representatives $(E_1\to E_0)\xrightarrow{(\phi_1,\phi_0)}\tau_{\le1}L_{X/B}$ with $(E_1\to E_0)$ a $G$-linearized two-term complex of   locally free coherent sheaves on $X$. 

We first recall the construction of Graber-Pandharipande. For a $G$-linearized coherent sheaf $\sF$ on a scheme $X$, we let $C(\sF)\to X$ denote the associated abelian cone $\Spec_{\sO_X}\Sym^*\sF$; in case $\sF$ is locally free, this is just the $G$-linearized vector bundle $\V(\sF)\to X$.

Let $[\phi_\bullet]$  be a  $G$-linearized   perfect obstruction theory on $Z$. Since objects of $\Sch^G/B$ have by definition a $G$-equivariant locally closed immersion in some $\P(\sV)$, for $\sV$ a locally free $G$-linearized sheaf on $B$,   $Z$ admits a  closed immersion $i\:Z\to M$ in $\Sch^G/B$, with $M$  smooth over $B$. Thus $[\phi_\bullet]$ admits a representative as a map of  $G$-linearized complexes 
\[
\phi_\bullet\:(F_1\xrightarrow{\del} F_0)\to (\sI_Z/\sI_Z^2\xrightarrow{d}i^*\Omega_{M/B})
\]
with the $F_i$  $G$-linearized  locally free coherent sheaves on $Z$. The assumption that $\phi_\bullet$ defines a perfect obstruction theory is equivalent to the exactness of the sequence
\[
F_1\to F_0\oplus \sI_Z/\sI_Z^2\xrightarrow{\gamma}i^*\Omega_{M/B}\to 0.
\]
Let $Q=\ker(\gamma)$ with surjection $\pi\:F_1\to Q$. Then we have the  exact sequence of abelian cones
\[
0\to i^*T_{M/B}\to C(\sI_Z/\sI_Z^2)\times_Z\V(F_0)\to C(Q)\to 0,
\]
and the surjection $\pi$ gives rise to a  closed immersion
\[
i_\pi\:C(Q)\hookrightarrow \V(F_1).
\]
The surjection $\Sym^*\sI_Z/\sI_Z^2\to \oplus_{n\ge0}\sI_Z^n/\sI_Z^{n+1}$ gives the closed immersion $C_{Z/M}\hookrightarrow C(\sI_Z/\sI_Z^2)$, inducing the closed immersion 
$C_{Z/M}\times_Z\V(F_0)\to C(\sI_Z/\sI_Z^2)\times_Z\V(F_0)$ with image containing the sub-vector bundle $i^*T_{M/B}$. Let $D:= C_{Z/M}\times_Z\V(F_0)$, giving   the quotient cone
\[
D^\vir:=i^*T_{M/B}\backslash D,
\]
which is naturally a closed sub-cone of $C(Q)$. The closed immersion $i_\pi$ thus induces a closed immersion
\[
\iota_{F_\bullet}\:D^\vir\hookrightarrow \V(F_1).
\]

We have the projections
\[
\pi_1\:D\to C_{Z/M},\quad \pi_2\:D\to D^\vir.
\]
The maps $\pi_1$ and $\pi_2$  exhibit $D$ as an affine space bundle over $C_{Z/M}$ and over $D^\vir$, respectively, and are both maps over $Z$. 

The smooth pull-back maps
\[
\pi_1^!\:\sE^\BM_{a,b}(C_{Z/M}/B, v)\to \sE^\BM_{a,b}(D/B, v+ F_0),
\]
\[
\pi_2^!\:\sE^\BM_{a,b}(D^\vir/B, v)\to \sE^\BM_{a,b}(D/B, v+i^*\Omega_{M/B})
\]
are both isomorphisms, by $\A^1$-homotopy invariance. 

We have defined in \cite[Definition 4.1]{LevineVirt} the fundamental class 
\[
[C_{Z/M}]\in \sE^\BM(C_{Z/M}/B, \Omega_{M/B}); 
\]
this is the normal cone class defined here in \S\ref{sec:FundamentalClasses}.

Using the isomorphisms $\pi_1^!$, $\pi_2^!$, $[C_{Z/M}]\in \sE^\BM(C_{Z/M}/B, i^*\Omega_{M/B})$ gives rise to the class
\[
[D^\vir]:=(\pi_2^!)^{-1}(\pi_1^![C_{Z/M}])\in \sE^\BM(D^\vir/B, F_0).
\]

\begin{definition}\label{def:GPVirtualClass} Let $\sE$ be in $\SH^G(B)$ and let $\phi_\bullet\:F_\bullet\to (\sI_Z/\sI_Z^2\xrightarrow{d}i^*\Omega_{M/B})$ be a representative of a  $G$-linearized  perfect obstruction theory on a ($G$-equivariant) embedded $B$-scheme $i\:Z\to M$ with $M$ smooth over $B$.  Let $0_{\V(F_1)}\:Z\to \V(F_1)$ be the zero-section. The $\sE$-valued {\em Graber-Pandharipande virtual fundamental class} is defined as
\[
[Z,\phi,i]^\vir_{\sE, GP}:=0^!_{\V(F_1)}(\iota_{F_\bullet*}[D^\vir]) \in \sE^\BM(Z/B, F_\bullet).
\]
\end{definition}

\begin{remark} Taking $\sE=H\Z$, the motivic spectrum representing motivic cohomology, $B=\Spec k$ for some perfect field $k$, and $G$ the trivial group, we have the natural isomorphism $H\Z^\BM(Z/B, F_\bullet)\cong \CH_r(Z)$, 
where $r$ is the virtual rank of $F_\bullet$, and the Graber-Pandharipande class $[Z,\phi,i]^\vir_{H\Z, GP}\in \CH_r(Z)$ is equal to the virtual fundamental class as defined by Behrend-Fantechi \cite{BF}.
\end{remark}

We have the virtual fundamental class $[Z,[\phi]]_\sE^\vir\in \sE^\BM(Z/B, F_\bullet)$, as defined in \cite{LevineVirt}. Note that $[Z,[\phi]]_\sE^\vir$ depends only on the choice of perfect obstruction theory $[\phi]\:E_\bullet\to \tau_{\le1}L_{Z/B}$ in $D_G(Z)$. 

\begin{proposition}  Let  $\sE$, $i\:Z\to M$ and $\phi$ be as in Definition~\ref{def:GPVirtualClass}. Then
\[
[Z,[\phi]]_\sE^\vir=[Z,\phi,i]^\vir_{\sE, GP}\in \sE^\BM(Z/B, F_\bullet).
\]
\end{proposition}

\begin{proof} We prove the result by tracing through the  step-by-step construction of $[Z,]\phi]]_\sE^\vir$. This consists of starting with a representative  
\[
\phi_\bullet\:F_\bullet\to (\sI_Z/\sI_Z^2\to i^*\Omega_{M/B})
\]
 of $[\phi]$, and modifying $\phi_\bullet$ and $Z$ step-by-step, until we arrive at a particularly simple case. At each stage, there is a canonical isomorphism of the corresponding Borel-Moore homology groups and at the final stage, we have an explicit construction of a class. The class $[Z,[\phi]]^\vir_\sE\in \sE^\BM(Z/B, F_\bullet)$ is then defined the image of this explicit class through the sequence of canonical isomorphisms of Borel-Moore homology groups. The Graber-Pandharipande class, however, will be defined at each stage of the process. 

Here we will go through this procedure in reverse order, starting with the special situation needed for the explicit construction of $[Z,[\phi]]^\vir_\sE$, and going back step-by-step until we reach the general situation. We will identify the explicit construction of $[Z,[\phi]]^\vir_\sE$ with the Graber-Pandharipande construction at the initial stage of our reverse-order process, and show that the classes one constructs by the 
Graber-Pandharipande method are compatible with the canonical isomorphisms of Borel-Moore homology at each stage in the reverse process, which will prove the result. Now for the details. 

First suppose that the obstruction theory $\phi_\bullet\:F_\bullet\to (\sI_Z/\sI_Z^2\to i^*\Omega_{M/B})$ is {\em normalized}, meaning that  $\phi_1\:F_1\to \sI_Z/\sI_Z^2$ is surjective, and is
 {\em reduced},  meaning $F_0=i^*\Omega_{M/B}$ and $\phi_0=\id$. In this case we have the canonical isomorphism
\[
D^\vir:=C_{Z/M}\times_Z\V(F_0)/i^*T_{M/B}\cong
C_{Z/M}\times_Zi^*T_{M/B}
/i^*T_{M/B}\cong C_{Z/M},
\]
and via this isomorphism we have 
\[
[D^\vir]=[C_{Z/M}]\in \sE^\BM(D^\vir/B, F_0)=\sE^\BM(C_{Z/M}/B, i^*\Omega_{M/B}). 
\]
Also in this case, the surjection $\phi_1$ gives the closed immersion $i\:C_{Z/M}\to \V(F_1)$, which gets identified with the closed immersion $\iota_{F_\bullet}\:D^\vir\to \V(F_1)$. If in addition we assume that $Z$ and $M$ are affine, then $[Z,\phi]^\vir_\sE$ is defined by
\[
[Z,\phi]^\vir_\sE:=0_{\V(F_1)}^!(i_*[C_{Z/M}])\in \sE^\BM(Z/B, F_\bullet),
\]
which thus agrees with $[Z,\phi,i]^\vir_{\sE, GP}:=0^!_{\V(F_1)}(\iota_{F_\bullet*}[D^\vir])$.

Now assume that $Z$ and $M$ are affine and $\phi_1$ is surjective. In this case $\phi_0$ is also surjective; let $K_i$ be the kernel of $\phi_i$. Then $K_0$ is locally free and lifts isomorphically via the differential $\del\:F_1\to F_0$ to a locally free summand of $K_1$ (see \cite[\S 6]{LevineVirt} for proofs of these assertions). Taking the quotient of $F_\bullet$ by the subcomplex $(K_0\xrightarrow{\id}K_0)$ gives us the  normalized perfect obstruction theory $\phi_\bullet'\:F_\bullet'\to  (\sI_Z/\sI_Z^2\to i^*\Omega_{M/B})$, and our original obstruction theory is $F_\bullet'\oplus (K_0\xrightarrow{\id}K_0)$ with $\phi_0$ and $\phi_1$ the zero map on $K_0\subset F_0$ and on $K_0\subset F_1$. Also, this says that $\phi'_0\:F'_0\to i^*\Omega_{M/B}$ is an isomorphism, so  by a ``change of coordinates'', we may assume that $\phi_\bullet'\:F_\bullet'\to  (\sI_Z/\sI_Z^2\to i^*\Omega_{M/B})$ is reduced and normalized. 

The quotient map $\pi\: F_\bullet\to F'_\bullet$ is a quasi-isomorphism of perfect complexes, inducing a canonical isomorphism $\sE^\BM(Z/B, F_\bullet)\cong \sE^\BM(Z/B, F'_\bullet)$. We have our cone $D$ and closed immersion $\iota_{F_\bullet}\:D^\vir\to \V(F_1)$ and the corresponding cone $D'$ and closed immersion $\iota_{F'_\bullet}\:D^{\prime\vir}\to \V(F'_1)$. The surjection $\pi$ induces closed immersions $C(\pi)\:D^{\prime\vir}\to D^\vir$ and $\V(\pi)\:\V(F_1')\to \V(F_1)$. We claim that
\\[5pt]
i. $C(\pi)^![D^\vir]=[D^{\prime\vir}]$,\\[2pt]
ii. $\V(\pi)^!(\iota_{F_\bullet*}[D^\vir])=\iota_{F'_\bullet*}[D^{\prime\vir}]$,\\[2pt]
iii. $0^!_{\V(F_1)}(\iota_{F_\bullet*}[D^\vir])=0^!_{\V(F'_1)}(\iota_{F'_\bullet*}[D^{\prime\vir}])$.
\\[5pt]

To prove (i) and (ii), let $p_{D^{\prime\vir}}\:D^{\prime\vir}\to Z$, $p_{\V(F'_1)}\:\V(F'_1)\to Z$ be the structure morphisms. Then we have
\[
\V(F_1)=\V(F_1')\oplus  \V(K_0),\ D^\vir=D^{\prime\vir}\oplus \V(K_0),
\]
via which $C(\pi)$ and $\V(\pi)$ become the inclusion as the first summand. The respective projections onto the first summands,
\[
\pi_{\V(F_1')}\:\V(F_1)\to \V(F_1'),\ \pi_{D^{\prime\vir}}\:D^\vir\to D^{\prime\vir},
\]
thus exhibit $D^\vir$ as the vector bundle $p_{D^{\prime\vir}}^*(\V(K_0))$ over $D^{\prime\vir}$, and $\V(F_1)$ as the vector bundle $p_{\V(F_1')}^*(\V(K_0))$ over $\V(F_1')$; the inclusions $C(\pi)$ and $\V(\pi)$ are the respective zero-sections. Thus $C(\pi)^!$ is the inverse to the smooth pull-back map $\pi_{D'}^!$ and $\V(\pi)^!$ is inverse to the smooth pull-back map  $\pi_{\V(F_1')}^!$.

The identity
\[
[D^\vir]=\pi_{D^{\prime\vir}}^!([D^{\prime\vir}])
\]
follows immediately from the definition of $[D^\vir]$ and $[D^{\prime\vir}]$ and the functoriality of smooth pull-back, proving (i). The identity in (ii) is equivalent to the identity
\[
\iota_{F_\bullet*}(\pi_{D^{\prime\vir}}^!([D^{\prime\vir}]))=\pi_{\V(F_1')}^!(\iota_{F'_\bullet*}[D^{\prime\vir}]).
\]
Since the square
\[
\xymatrix{
D^\vir\ar[r]^{\iota_{F_\bullet}}\ar[d]_{\pi_{D'}}&\V(F_1)\ar[d]^{\pi_{\V(F_1')}}\\
D^{\prime\vir}\ar[r]^{\iota_{F'_\bullet}}&\V(F_1')
}
\]
is cartesian, this latter identity follows from the compatibility of proper push-forward and smooth pull-back in cartesian squares: $\iota_{F_\bullet*}\circ \pi_{D^{\prime\vir}}^!=\pi_{\V(F_1')}^!\circ\iota_{F'_\bullet*}$, Proposition~\ref{prop:LciProperties}. This proves (ii). 

The proof of (iii) is similar:  Let $p_{\V(F_1)}\:\V(F_1)\to Z$ be the structure map. The identity is equivalent to
\[
p_{\V(F'_1)}^!(p_{\V(F_1)}^!)^{-1}(\iota_{F_\bullet*}([D^\vir]))=\iota_{F'_\bullet*}[D^{\prime\vir}],
\]
and we have 
\[
\iota_{F_\bullet*}([D^\vir])= \pi_{\V(F_1')}^!(\iota_{F'_\bullet*}[D^{\prime\vir}]).
\]
By functoriality of smooth push-forward, we have 
\[
 p_{\V(F_1)}^!=p_{\V(F'_1)}^!\circ  \pi_{\V(F_1')}^!,
 \]
and thus 
\[
p_{\V(F'_1)}^!(p_{\V(F_1)}^!)^{-1}(\iota_{F_\bullet*}([D^\vir]))
=p_{\V(F'_1)}^!\circ (p_{\V(F_1)}^!)^{-1}\circ \pi_{\V(F_1')}^!(\iota_{F'_\bullet*}[D^{\prime\vir}])
=\iota_{F'_\bullet*}[D^{\prime\vir}],
\]
proving (iii). In other words, $[Z,\phi,i]^\vir_{\sE, GP}=[Z',\phi',i]^\vir_{\sE, GP}$ in $\sE^\BM(Z/B, F_\bullet)=\sE^\BM(Z/B, F'_\bullet)$. 

Now suppose that $\phi_\bullet$ is normalized, but $M$ is only quasi-projective over the affine base-scheme $B$. Let $q\:\tilde{M}\to M$ be a Jouanolou cover of $M$, that is, $q$ is a (Zariski) torsor for a $G$-linearized vector bundle on $M$ and $\tilde{M}$ is affine. Letting $\tilde{Z}=Z\times_M\tilde{M}$, we have the Jouanolou cover $q_Z\:\tilde{Z}\to Z$ and the closed immersion $\tilde{i}:=p_2\:\tilde{Z}\to \tilde{M}$. Since $\tilde{M}$ is affine, we can split the canonical surjection $\Omega_{\tilde{M}/B}\to \Omega_{\tilde{M}/M}$, giving the isomorphism $\Omega_{\tilde{M}/B}=q^*\Omega_{M/B}\oplus 
\Omega_{\tilde{M}/M}$. We define $\tilde{F}_1=q_Z^*F_1$,  $\tilde{F}_0:=q_Z^*F_0\oplus\Omega_{\tilde{M}/M}$, with differential $\del_{\tilde{F}}\:\tilde{F}_1\to 
\tilde{F}_0$ to be $q_Z^*(\del_{F})$ followed by the inclusion $q_Z^*F_0\to \tilde{F}_0$. Similarly, define $\tilde{\phi}_1:=q_Z^*\phi_1$ and $\tilde{\phi}_0=q_Z^*\phi_0\oplus\id_{\Omega_{\tilde{M}/M}}$. 

This gives us the normalized perfect obstruction theory 
\[
\tilde{\phi}_\bullet\:\tilde{F}_\bullet\to 
(\sI_{\tilde{Z}}/\sI_{\tilde{Z}}^2\to \tilde{i}^*\Omega_{\tilde{M}/B})
\]
on $\tilde{Z}$. Since $q_Z\:\tilde{Z}\to Z$ is the pull-back to $Z$ of the vector bundle torsor $\tilde{M}\to M$, we see that the relative tangent bundle $T_{q_Z}$ is canonically isomorphic to $\V(i^*\Omega_{\tilde{M}/M})$. By homotopy invariance, the smooth pull-back map
\[
q_Z^!\:\sE^\BM(Z/B, F_\bullet)\to \sE^\BM(\tilde{Z}/B, q_Z^*F_\bullet+L_{q_Z})
\]
is an isomorphism and the isomorphism $L_{q_Z}\cong  i^*\Omega_{\tilde{M}/M}$ gives us the canonical isomorphism
\[
\sE^\BM(\tilde{Z}/B, q_Z^*F_\bullet+L_{q_Z})\cong \sE^\BM(\tilde{Z}/B, \tilde{F}_\bullet).
\]
Let $\iota_{\tilde{F}_\bullet}\:\tilde{D}^\vir\to \V(\tilde{F}_1)$ be the closed immersion associated to the perfect obstruction theory $\tilde{\phi}_\bullet$.  We have  the commutative diagram
\[
\xymatrix{
&\V(\tilde{F}_1)\ar[dr]^{p_{\V(\tilde{F}_1}}\ar[dd]^(.67){q_{\V(F_1)}}|\hole\\
\tilde{D}^\vir\ar[dd]_{q_D}\ar[ur]^{\iota_{\tilde{F}_\bullet}}\ar[rr]^(.33){p_{\tilde{D}}}&&\tilde{Z}\ar[dd]^{q_Z}\\
&\V(F_1)\ar[dr]^{p_{\V(F_1)}}\\
D^\vir\ar[ur]^{\iota_{F_\bullet}}\ar[rr]_{p_D}&&Z
}
\]
with all squares cartesian.  Arguing as in the previous case, we see that
\[
q_Z^!(0^!_{\V(F_1)}(\iota_{F_\bullet*}[D^\vir]))=0^!_{\V(\tilde{F}_1)}(\iota_{\tilde{F}_\bullet*}[\tilde{D}^\vir]),
\]
in other words $q_Z^!([Z,\phi,i]^\vir_{\sE, GP}=[\tilde{Z},\tilde{\phi},\tilde{i}]^\vir_{\sE, GP}$
in $\sE^\BM(\tilde{Z}/B, \tilde{F}_\bullet)$.

We conclude the argument with the general case. If 
\[
\phi_\bullet\:F_\bullet\to (\sI_Z/\sI_Z^2\xrightarrow{d} i^*\Omega_{M/B})
\]
is an arbitrary $G$-linearized perfect obstruction theory on our quasi-projective $Z$, with $i:Z\to M$ a closed immersion in $\Sch^G/B$, $M\in \Sm^G/B$, there is a $G$-linearized locally free sheaf $\sF$ and a $G$-equivariant surjection $a\:\sF\to \sI_Z/\sI_Z^2$. Replace $\phi_\bullet$ with $\phi_\bullet\:F'_\bullet\to (\sI_Z/\sI_Z^2\xrightarrow{d} i^*\Omega_{M/B})$ with $F_i':=F_i\oplus \sF$, $\del_{F_\bullet'}=\del_{F_\bullet}\oplus\id_\sF$,  $\phi'_1=\phi_1+a$ and $\phi'_0=\phi_0+d\circ a$. Then $\phi'_\bullet$ is normalized. As above, we have a canonical isomorphism $\sE^\BM(Z, F_\bullet)=\sE^\BM(Z, F'_\bullet)$ and via this isomorphism, we have
\[
0^!_{\V(F_1)}(\iota_{F_\bullet*}[D^\vir])=0^!_{\V(F'_1)}(\iota_{\tilde{F}'_\bullet*}[D^{\prime\vir}]),
\]
that is, $[Z,\phi,i]^\vir_{GP}=[Z',\phi',i]^\vir_{GP}$
in $\sE^\BM(Z/B, F_\bullet))=\sE^\BM(Z/B, F'_\bullet)$

Thus, we have verified the identity of the two classes at the initial stage of our reverse process, and have shown that the Graber-Pandharipande classes are compatible with the canonical isomorphisms of Borel-Moore homology at each stage of the process. This completes the proof.

\end{proof}

\begin{remark} The main result of \cite{LevineVirt} is that the virtual fundamental class $[Z,[\phi]]^\vir_\sE\in \sE^\BM(Z/B, F_\bullet)$  depends only on the perfect obstruction theory $[\phi]$ on $Z$ and is independent of the various choices made in its construction. Thus, the Graber-Pandharipande class $[Z,\phi,i]^\vir_{GP}$ is also independent of the choice of embedding $i$ and the choice of representative $\phi_\bullet$ for $[\phi]$. One can also show this directly by using a small modification of the proof of this independence for the Behrend-Fantechi classes \cite[Proposition 5.3]{BF}.
\end{remark}

\section{Equivariant virtual fundamental classes via the algebraic Borel-construction}\label{sec:TEGVirClass}

In this section, we show how a $G$-linearized perfect obstruction theory gives rise to a virtual fundamental class in the equivariant Borel-Moore homology defined by the algebraic Borel-construction of \S\ref{sec:TEGEquivCoh}.  We fix a motivic commutative ring spectrum $\sE\in \SH(B)$; note that we do not use the equivariant category $\SH^G(B)$. 

We retain the assumptions on $B$ and $G$ from the previous section. In particular,   a $G$-linearized perfect obstruction $\phi$ on some $X\in \Sch^G/B$ admits a representative $(E_1\to E_0)\xrightarrow{(\phi_1,\phi_0)}\tau_{\le1}L_{X/B}$ with $(E_1\to E_0)$ a $G$-linearized two-term complex of   locally free coherent sheaves on $X$. We will assume that $X$ is in the subcategory $\Sch^G_q/B$ and  admits a $G$-equivariant closed immersion $\iota\:X\to Y$ for some $Y\in \Sm^G_q/B$. This gives us for each $j$ the $G$-equivariant closed immersion $\iota_j:=\iota\times \id_{E_jG}\:X\times E_jG\to Y\times E_jG$ and induced closed immersion $\alpha_j\: X\times^GE_jG\to  Y\times^GE_jG$.

Up to now, we have been supressing the relative nature of perfect obstruction theories and their associated virtual fundamental classes,  because we have always been working over a fixed base-scheme $B$; the relative nature comes from the fact that a perfect obstruction theory on a $B$-scheme $X$ is a map to the relative cotangent complex $L_{X/B}$. As we will be varying the base-schemes in our construction of the equivariant classes, we will include the choice of base-scheme in the notation, by writing $[X/B, \phi]^\vir_\sE\in \sE^\BM(X/B, E_\bullet)$ instead of $[X, \phi]^\vir_\sE$. We similarly denote the cone class for the closed immersion $\alpha_j$ by $[C_{\alpha_j}/B_jG]$.

In what follows, we will be defining a family of perfect obstruction theories $\phi_j:E_{\bullet,j}\to L_{X\times^GE_jG/B_jG}$, giving virtual fundamental classes  
\[
[X\times^GE_jG/B_jG, \phi_j]^\vir_\sE\in \sE^\BM(X\times^GE_jG/B_jG, E_{\bullet,j}).
\]
We then show that the family $\{[X\times^GE_jG/B_jG, \phi_j]^\vir_\sE\}_j$ is compatible with respect to the transition maps $i^!_{X,G,j}$,  giving a well-defined element in the limit which will be our equivariant virtual fundamental class
\[
[X, \phi]^\vir_{\sE,G}\in \sE^\BM_G(X/BG, E_\bullet).
\] 

In detail, letting $p_1\:X\times_BE_jG\to X$ be the projection, we have  $p_1^*L_{X/B}\cong L_{X\times_BE_jG/E_jG}$, and   the complex $L_{X/B,j}$ on $X\times^GE_jG$ induced by $p_1^*L_{X/B}$ via descent is canonically isomorphic to $L_{X\times^GE_jG/B_jG}$. Letting 
\[
\phi_j\:E_{\bullet,j}\to \tau_{\le1}L_{X/B,j}\cong \tau_{\le1}L_{X\times^GE_jG/B_jG}
\]
be the map induced by $p_1^*\phi$ by descent, $\phi_j$ thus defines a perfect obstruction theory on 
$X\times^GE_jG$ relative to $B_jG$, giving us the virtual fundamental class
\begin{equation}\label{eqn:VFCDef}
[X\times^GE_jG/B_jG, \phi_j]^\vir_{\sE}\in \sE^\BM(X\times^GE_jG/B_jG, E_{\bullet,j}).
\end{equation}

We now show that this family of classes is compatible with the transition maps $i^!_{X,G,j}$ for $\{\sE^\BM(X\times^GE_jG/B_jG, E_{\bullet,j})\}_j$.

\begin{lemma}\label{lem:VIrFundClasspull-back} Let $i\:B_0\to B_1$ be a regular closed immersion in $\Sm/B$. Let $\psi\:(F_1\to F_0)\to \tau_{\le1}L_{Z_1/B_1}$ be a perfect obstruction theory on some $Z_1\in \Sch/B_1$, relative to $B_1$.    Let $Z_0=Z_1\times_{B_1}B_0$ and assume that the cartesian square
\begin{equation}\label{eqn:TransCartSq}
\xymatrix{
Z_0\ar[r]\ar[d]&Z_1\ar[d]\\
B_0\ar[r]^i&B_1
}
\end{equation}
is Tor-independent, giving the canonical isomorphism $i^*L_{Z_1/B_1}\cong L_{Z_0/B_0}$ and the induced perfect obstruction theory $i^*\psi\:(i^*F_1\to i^*F_0)\to \tau_{\le1}L_{Z_0/B_0}$. 

Choose a  closed immersion $\iota\:Z_1\to Y_1$, with $Y_1\in \Sm_q/B_1$, giving the
 closed immersion $\iota_0\:Z_0\to Y_0:=Y_1\times_{B_1}B_0$. Suppose that the canonical map $\beta\:C_{Z_0/Y_0}\to i^*C_{Z_1/Y_1}$ is an isomorphism. Then  the relative pullback
\[
i^!\:\sE^\BM(Z_1/B_1, F_\bullet)\to  \sE^\BM(Z_0/B_0, i^*F_\bullet)
\]
satisfies
\[
i^!([Z_1/B_1,\psi]_\sE^\vir)=[Z_0/B_0,\psi_0]_\sE^\vir.
\]
\end{lemma}

\begin{proof} As in \S\ref{subsec:RelPull-back}, the purity isomorphisms for $B_0\to B$, $B_1\to B$
give  us the canonical isomorphism 
\[
\vartheta_i\:\sE^\BM(Z_0/B_1,v-\sN_i)\xrightarrow{\sim} \sE^\BM(Z_0/B_0, v)
\]
for $v\in \sK(Z_0)$. 

Letting $\sI_{Z_1}$ be the ideal sheaf of $Z_1$ in $Y_1$, we have the map $\bar\psi\:(F_1\to F_0)\to (\sI_{Z_1}/\sI_{Z_1}^2\to\iota^*\Omega_{Y_1/B_1})$ induced by $\psi$,  the corresponding closed immersion $\iota_0:Z_0\to Y_0:=Y_1\times_{B_1}B_0$ and the corresponding map $i^*\bar\psi\:(i^*F_1\to i^*F_0)\to (\sI_{Z_0}/\sI_{Z_0}^2\to\iota^*\Omega_{Y_0/B_0})$ induced by $\psi$. 

Following the Graber-Pandharipande construction given in  \S\ref{sec:VirClass}, we have the cone $D_1^\vir$ with fundamental class $[D_1^\vir/B_1]$, and closed immersion $\iota_{F_\bullet}\:D_1^\vir\hookrightarrow \V(\sF_1)$, giving the virtual fundament class $[Z_1/B_1,\psi]^\vir_\sE$ as
\[
[Z_1/B_1,\psi]^\vir_\sE=0_{\V(\sF_1)}^!(\iota_{F_\bullet*}([D_1^\vir/B_1])).
\]
We also have the corresponding  cone $D_0^\vir$ with fundamental class $[D_0^\vir/B_0]$, and closed immersion $\iota_{i^*F_{\bullet}}\:D_0^\vir\hookrightarrow \V(i^*F_1)$, giving the virtual fundamental class $[Z_0/B_0,i^*\psi]^\vir_\sE$ as
\[
[Z_0/B_0,i^*\psi]^\vir_\sE=0_{\V(i^*\sF_1)}^!(\iota_{i^*F_\bullet*}([D_0^\vir/B_0])).
\]

By our assumptions, the canonical map $\beta\:C_{Z_0/Y_0}\to i^*C_{Z_1/Y_1}$ is an isomorphism.   By Proposition~\ref{prop:FundClassIdentity}, this gives the identity
\[
[C_{Z_0/Y_0}/B_0]=i^![C_{Z_1/Y_1}/B_1]\in \sE^\BM(C_{Z_0/Y_0}/B_0, \Omega_{Y_0/B_0}).
\]
Using the zig-zag of isomorphisms of Borel-Moore homology used to define $[D_1/B_1]$, $[D_1^\vir/B_1]$,   and $[D_0/B_0]$, $[D_0^\vir/B_0]$, and using the compatibility of the relative pullback/refined Gysin map with smooth pull-back (Proposition~\ref{prop:LciProperties})  this implies
\[
[D_0^\vir/B_0]=i^![D_1^\vir/B_1]\in  \sE^\BM(D_0^\vir/B_0,i^*F_0).
\]

Thus
\begin{align*}
[Z_0/B_0,i^*\psi]^\vir_\sE&=0_{\V(i^*\sF_1)}^!(\iota_{i^*F_\bullet*}([D_0^\vir/B_0]))\\
&=0_{\V(i^*\sF_1)}^!(\iota_{i^*F_\bullet*}(i^![D_1^\vir/B_1]))\\
&=0_{\V(i^*\sF_1)}^!i^!(\iota_{i^*F_\bullet*}([D_1^\vir/B_1]))\\
&=i^!(0_{\V(\sF_1)}^!(\iota_{i^*F_\bullet*}([D_1^\vir/B_1])))\\
&=i^![Z_1/B_1,\psi]^\vir_\sE.
\end{align*}
The third identity is Proposition~\ref{prop:LciProperties}  and the fourth is Corollary~\ref{cor:commutativity}.
\end{proof} 

\begin{def/prop}\label{defprop:EquivVFC} Let  $\phi_\bullet\:E_\bullet\to \tau_{\le1}L_{X/B}$ be a $G$-linearized perfect obstruction theory on $X\in \Sch^G_q/B$. Then the family of virtual fundamental classes of \eqref{eqn:VFCDef},
\[
\{[X\times^GE_jG/B_jG, \phi_j]^\vir_\sE\in \sE^\BM(X\times^GE_jG/B_jG, E_{\bullet,j})\}_j
\]
gives a well-defined element
\[
[X,\phi]^\vir_{\sE,G}\in \sE_G^\BM(X/BG, E_\bullet).
\]
\end{def/prop}

\begin{proof}   We have the cartesian, Tor-independent diagram
\[
\xymatrix{
X\times^GE_jG\ar[r]\ar[d]&X\times^GE_{j+1}G\ar[d]\\
B_jG\ar[r]^{i_{G,j}}&B_{j+1}G,
}
\]
and the perfect obstruction theory $\phi_{j+1}$ on $X\times^GE_{j+1}G$, relative to $B_{j+1}G$. Choose a $G$-equivariant closed immersion $\iota\:X\to Y$ with $Y\in \Sm^G_q/B$, giving the closed immersions
\[
\iota_j\:X\times^GE_jG\to Y\times^GE_jG
\]
with $Y\times^GE_jG$ smooth over $B_jG$ and $\iota_j\cong \iota_{j+1}\times_{B_{j+1}G}B_jG$. This gives the isomorphism
\[
C_{\iota_j}\cong C_{X/Y}\times^GE_jG
\]
and thus the canonical map $\beta\:C_{\iota_j}\to C_{\iota_{j+1}}\times_{B_{j+1}G}B_jG$ is an isomorphism. We also have  the canonical isomorphism of perfect obstruction theories $i_{G,j}^*\phi_{j+1}\cong \phi_j$. Thus, we may apply  Lemma~\ref{lem:VIrFundClasspull-back} to show that $i_{G,j}^![X\times^GE_{j+1}G/B_{j+1}G, \phi_{j+1}]^\vir_\sE=[X\times^GE_jG/B_jG, \phi_j]^\vir_\sE$,  as desired. 
\end{proof}

\begin{proposition}\label{prop:VirtChangeOfGroup} We take $B=\Spec k$, $k$ a field and suppose $\sE$ is bounded. Let $\phi\:E_\bullet\to \tau_{\le1}L_{X/B}$ be $G$-linearized perfect obstruction theory  on some $X\in \Sch^G/B$ and let $\rho:H\to G$ be a homomorphism of smooth affine group schemes over $B$, giving the change of group map  $\rho^\#:\sE^\BM_{G,**}(-/BG,-)\to  \sE^\BM_{H,**}(-/BH,-)$. Then\\[5pt]
1. $\rho^\#([X,[\phi]]_{\sE,G}^\vir)=[X,[\phi]]_{\sE,H}^\vir$.\\[2pt]
2. Suppose $p_X:X\to B$ is proper, so the proper pushforward with respect ot $p_X$ gives the equivariant degree maps 
\[
\deg^{\sE,G}_B:\sE^\BM_G(X/BG)\to \sE^\BM_G(B/BG)=\sE^{0,0}_G(B)
\]
and
\[
\deg^{\sE,H}_B:\sE^\BM_H(X/BH)\to \sE^\BM_H(B/BH)=\sE^{0,0}_H(B).
\]
 Suppose that $\sE$ is $\SL$-oriented, $\rnk E_\bullet=0$ and we have a $G$-equivariant orientation $\tau:\det E_\bullet\xrightarrow \sL^{\otimes 2}$ for some $G$-linearized invertible sheaf $\sL$ on $X$, giving the canonical isomorphisms
\[
\sE^\BM_{G}(X/BG,E_\bullet)\cong \sE^\BM_{G}(X/BG),\ \sE^\BM_H(X/BH,E_\bullet)\cong \sE^\BM_H(X/BH)
\]
and giving the elements
\[
\deg_B^{\sE, G}([X,[\phi]]_{\sE,G}^\vir)\in \sE^{0,0}_G(B),\ \deg_B^{\sE, H}([X,[\phi]]_{\sE,H}^\vir)\in \sE^{0,0}_H(B)
\]
Then 
\[
\rho^\#(\deg_B^{\sE, G}([X,[\phi]]_{\sE,G}^\vir))=
\deg_B^{\sE, H}([X,[\phi]]_{\sE,H}^\vir).
\]
\end{proposition}

\begin{proof} Since the degree maps are just the respective pushforward maps $p_{X*}$ on $G$, respectively $H$, equivariant $\sE$-Borel-Moore homology, (2) follows from (1) and Proposition~\ref{prop: ChangeOfGroupProperties}. 

Since we will be working over different base-schemes, we include this information in the notation for the fundamental classes and virtual fundamental classes, when needed.

For (1), we refer to  the construction of $[X,[\phi]]_{\sE,G}^\vir$ as the family \eqref{eqn:VFCDef} of Graber-Pandharipande  virtual fundmental classes, the individual terms being constructed as described in Definition~\ref{def:GPVirtualClass}. We first assume that $\rho$ is a closed embedding $i:H\to G$ of group schemes over $B$, so $\rho^\#=i^*$.

 We fix a $G$-equivariant embedding $X\hookrightarrow Y$ for some $Y\in \Sm^G/B$ and a closed embedding of group schemes $G\hookrightarrow \GL_n$. As before, we let $E_jH$ and $E_jG$ denote the $B$-scheme $E_j\GL_n$ with action by $H$ and $G$ given by restriction.
 
 For each $j$, this gives us cones $C_{X\times^GE_jG/Y\times^GE_jG}$ and $C_{X\times^HE_jH/X\times^HE_jH}$. The diagram
\[
\xymatrix{
C_{X\times^HE_jH/Y\times^HE_jG}\ar[r]^-{q^H_j}\ar[d]^{p^C_j}&X\times^HE_jH\ar[d]^{p^X_j}\ar[r]&Y\times^HE_jH\ar[d]^{p^Y_j}\ar[r]&B_jH\ar[d]^{\pi_j}\\
C_{X\times^GE_jG/Y\times^GE_jG}\ar[r]^-{q^G_j}&X\times^GE_jG\ar[r]&Y\times^GE_jG\ar[r]&B_jG
}
\]
has all squares cartesian and the maps $\pi_j$,  $p^Y_j$, $p^X_j$ and ${p^C_j}$ are smooth.

By Lemma~\ref{lem:ConeClassSmoothPullback} below, we have
\begin{equation}\label{eqn:ConeClassSmoothPullback}
\pi_j^!([C_{X\times^GE_jG/Y\times^GE_jG}/B_jG])=[C_{X\times^HE_jG/Y\times^HE_jG}/B_jH]
\end{equation}
 
We have a similar cartesian diagram relating the  schemes $D^\vir_{jG}$ and $D^\vir_{jH}$ used to construct the respective virtual fundamental classes $[X\times^GE_jG/B_jG, \phi_j]^\vir_\sE$ and 
$[X\times^HE_jH/B_jH, \phi_j]^\vir_\sE$,
\[
\xymatrix{
D^\vir_{jH}\ar[r]\ar[d]^{p^D_j}&X\times^HE_jH\ar[d]^{p^X_j}\\
D^\vir_{jG}\ar[r]&X\times^GE_jG,
}
\]
and it follows from \eqref{eqn:ConeClassSmoothPullback} that the fundamental classes $[D^\vir_{jH}/B_jH]
$ and $[D^\vir_{jG}/B_jG]$ satisfy
\[
\pi^!_j([D^\vir_{jG}/B_jG])=[D^\vir_{jH}/B_jH].
\]
Similarly we have the closed immersions
\[
\iota_{E_{\bullet, j,G}}:D^\vir_{jG}\to \V(E_{1,j,G}),\ \iota_{E_{\bullet, j,H}}:D^\vir_{jH}\to \V(E_{1,j,H}),
\]
where $E_{1,j,G}$ is the locally free sheaf on  $X\times^GE_jG$  induced by $p_1^*E_1$ on $X\times E_jG$  and using descent;  the locally free sheaf $E_{1,j,H}$ on $X\times^HE_jH$  is defined similarly. This gives us  the   diagram with all squares cartesian
\[
\xymatrixcolsep{40pt}
\xymatrix{
D^\vir_{jH}\ar[r]^-{\iota_{E_{\bullet, j,H}}}\ar[d]^{p^D_j}&\V(E_{1,j,H})\ar[r]\ar[d]^{p^\V_j}&X\times^HE_jH\ar[d]^{p^X_j}\ar[r]&B_jH\ar[d]^{\pi_j}\\
D^\vir_{jG}\ar[r]^-{\iota_{E_{\bullet, j,G}}}&\V(E_{1,j,G})\ar[r]&X\times^GE_jG\ar[r]&B_jG
}
\]
Thus
\begin{align*}
\pi^!_j(0^!_{\V(E_{1,j,G})}(\iota_{E_{\bullet, j,G*}}([D^\vir_{jG}/B_jG])))&=
0^!_{\V(E_{1,j,H})}(\pi^!_j(\iota_{E_{\bullet, j,G*}}([D^\vir_{jG}/B_jG])))\\
&=0^!_{\V(E_{1,j,H})}(\iota_{E_{\bullet, j,H*}}(\pi^!_j([D^\vir_{jG}/B_jG])))\\
&=0^!_{\V(E_{1,j,H})}(\iota_{E_{\bullet, j,H*}}([D^\vir_{jH}/B_jH])).
\end{align*}
In other words,
\[
\pi^!_j([X\times^GE_jG/B_jG, \phi_j]^\vir_{\sE,G})=[X\times^HE_jH/B_jH, \phi_j]^\vir_{\sE,H}
\]

But  $[X,\phi]^\vir_{\sE, G}$ is given by the family $\{[X\times^GE_jG/B_jG, \phi_j]^\vir_{\sE,G}\}_j$ and $[X,\phi]^\vir_{\sE, H}$ is given by the family $\{[X\times^HE_jH/B_jH, \phi_j]^\vir_{\sE,H}\}_j$, and  the family of relative pullback maps $\pi^!_j$ are what we use to define $i^*$. Thus, we have
\[
i^*([X,\phi]^\vir_{\sE, G})=[X,\phi]^\vir_{\sE, H}
\]
proving (1) for the case of a closed embedding.

If $\rho$ is a projection $p_2:G_1\times G_2\to G_2$, we take embeddings $\iota_i:G_i\to \GL_{n_i}$, $i=1,2$ and  use the family of morphisms  $X\times^{G_1\times G_2}E_{j, \iota_1}G_1\times E_{j, \iota_2}G_2\to B_{j,\iota_1}G_1\times B_{j,\iota_2}G_2$ to define the $G_1\times G_2$ equivariant Borel-Moore homology, which we are allowed to do by Lemma~\ref{lem:Independence2}.

In this setup, we replace the smooth maps $\pi_j:B_jH\to B_jG$ with the smooth map $p_{2,j}:
B_{j,\iota_1}G_1\times B_{j,\iota_2}G_2\to  B_{j,\iota_2}G_2$,  replace $(-)\times^HE_jH$ with $(-)\times^{G_1\times G_2}E_{j, \iota_1}G_1\times E_{j, \iota_2}G_2$ throughout, and similarly  replace $(-)\times^GE_jG$ with $(-)\times^{G_2}E_{j, \iota_2}G_2$. Making these changes, the same proof as for a closed embedding shows that 
\[
p_2^*([X,\phi]^\vir_{\sE, G_2})=[X,\phi]^\vir_{\sE, G_1\times G_2}
\] 

For an arbitrary group scheme homomorphism $\rho$, we write $\rho=p_2\circ i_{\rho}$ and then $\rho^\#=i_{\rho}^*\circ p_2^*$, so the case of a closed embedding together with the case of a projection imply the general case. 
\end{proof}
 
 As an parallel to Proposition~\ref{prop:VirtChangeOfGroup}, we look at Euler classes rather that virtual fundamental classes; the Euler class is actually a special case of a virtual fundamental class, but the treatment of Euler classes on its own is considerably easier than going through the theory of virtual fundamental classes. 

\begin{proposition}\label{prop:EulerClassChangeOfGroup} We take $B=\Spec k$, $k$ a field and suppose $\sE$ is bounded. Let $V\to X$ be a $G$-vector on some $X\in \Sch^G/B$ with sheaf of sections $\sV$. Suppose that $X$ is an lci scheme over $B$, giving us the equivariant   Borel-Moore Euler class
$e^\BM_{G,\sE}(V)\in \sE^\BM_G(X/BG, L_{X/B}-\sV)$. 

Let $\rho:H\to G$ be a homomorphism of smooth affine group schemes over $B$, giving the change of group map  $\rho^\#:\sE^\BM_{G,**}(-/BG,-)\to  \sE^\BM_{H,**}(-/BH,-)$. Then\\[5pt]
1.  $\rho^\#(e^\BM_{G,\sE}(V))=
e^\BM_{H,\sE}(V)$.\\[2pt]
2. Suppose $p_X:X\to B$ is proper and lci, of pure dimension $d=\rnk(V)$. Then proper pushforward with respect ot $p_X$ gives the equivariant degree maps 
\[
\deg^{\sE,G}_B:\sE^\BM_G(X/BG)\to \sE^\BM_G(B/BG)=\sE^{0,0}_G(B)
\]
and
\[
\deg^{\sE,H}_B:\sE^\BM_H(X/BH)\to \sE^\BM_H(B/BH)=\sE^{0,0}_H(B).
\]
 Suppose  $\sE$ is $\SL$-oriented, $\rnk E_\bullet=0$ and we have a $G$-equivariant relative orientation $\rho:\det L_{X/B}\otimes \det\sV^\vee\xrightarrow{\sim}  \sL^{\otimes 2}$ for some $G$-linearized invertible sheaf $\sL$ on $X$, giving the canonical isomorphisms
\[
\sE^\BM_{G}(X/BG,L_{X/B}-\sV)\cong \sE^\BM_{G}(X/BG),\ \sE^\BM_H(X/BH,L_{X/B}-\sV)\cong \sE^\BM_H(X/BH)
\]
and giving the elements
\[
\deg_B^{\sE, G}(e^\BM_{G,\sE}(V))\in \sE^{0,0}_G(B),\ \deg_B^{\sE, H}(e^\BM_{H,\sE}(V)r)\in \sE^{0,0}_H(B)
\]
Then 
\[
\rho^\#(\deg_B^{\sE, G}(e^\BM_{G,\sE}(V)))=
\deg_B^{\sE, H}(e^\BM_{H,\sE}(V)).
\]
\end{proposition}
\begin{proof} (1) is just Proposition~\ref{prop: ChangeOfGroupProperties}(5) and (2) follows from (1) and the compatibility of $\rho^\#$ with proper pushforward  (Proposition~\ref{prop:VirtChangeOfGroup}(3ii)). 
\end{proof}

We will sometimes drop the subscripts $\sE, G$ from $[X,\phi]^\vir_{\sE,G}$ if there is no cause for confusion.

\section{Virtual localization}\label{sec:VirLoc}  With the machinery of \S\ref{sec:TEGEquivCoh} and \S\ref{sec:TEGVirClass}, together with our version of Vistoli's  lemma, we can essentially copy the constructions and arguments of Graber-Pandharipande to prove our virtual localization theorem. We concentrate on  the case $G=N$, with torus $T_1\subset N$, $T_1=\G_m$, and work over a perfect field $k$.

Fix an $X\in \Sch^N/k$, a $Y\in \Sm^N/k$ and a closed immersion $\iota\:X\hookrightarrow Y$ in $\Sch^N/k$. We also fix a $N$-equivariant perfect obstruction theory on $X$, represented by a two-term complex $E_\bullet:=(E_1\to E_0)$ of $N$-linearized locally free sheaves on $X$, together with a $N$-equivariant map $\phi_\bullet\:E_\bullet\to \tau_{\le1}L_{X/k}=(\sI_X/\sI_X^2\to i^*\Omega_{Y/k})$.

We recall from \cite[\S1]{GP} some facts about the $T_1$-fixed loci of $X$ and $Y$. 

We have the $T_1$-fixed subschemes $X^{T_1}\subset Y^{T_1}$, with $X^{T_1}=X\cap Y^{T_1}$. $Y^{T_1}$ is smooth; let $Y_1,\ldots, Y_s$ be the irreducible components of $Y^{T_1}$ with inclusion maps $i^Y_j\:Y_j\to Y$,  and let $X_j=X\cap Y_j$, so $X^{T_1}=\amalg_jX_j$. Let $\tilde{i}\:X^{T_1}\to X$, $\tilde{i}_j\:X_j\to X$ be the inclusions.

Let $\sF$ be a $T_1$-linearized coherent sheaf on $X_j$. Since the $T_1$-action on $X_j$  is trivial, we can decompose $\sF$ into weight spaces for the $T_1=\G_m$-action
\[
\sF=\oplus_r\sF_r
\]
where $t\in \G_m$ acts on $\sF_r$ by the character $\chi_r(t)= t^r$. 
if $\sF$ is locally free, so is each $\sF_r$. The subsheaf $\sF^\mov:=\oplus_{r\neq0}\sF_r$ is the {\em moving part} of $\sF$ and $\sF^\fix:=\sF_0$ is the {\em fixed part} of $\sF$. 

The map $\phi_\bullet$ induces the map 
\[
\tilde{\phi}\:\tilde{i}^*E_\bullet^\fix\to  \tau_{\le1}L_{X^{T_1}/k}
\]
By \cite[Proposition 1]{GP} $\tilde{\phi}$ defines a perfect obstruction theory on $X^{T_1}$. 

\begin{definition} The virtual conormal sheaf of $X^{T_1}$ in $X$, $\sN^\vir$,  is defined to be the perfect complex $\tilde{i}^*E_\bullet^\mov$.
\end{definition}

We have the order-four element $\sigma=\begin{pmatrix}0&1\\-1&0\end{pmatrix}\in N$, which together with $T_1$ generates $N$. Letting $\bar\sigma$ denote the image of $\sigma$ in $N/T_1$, $N/T_1$ is the cyclic group of order two generated by $\bar\sigma$.

\begin{lemma}\label{lem:NLin}  $\tilde{\phi}\:\tilde{i}^*E_\bullet^\fix\to  \tau_{\le1}L_{X^{T_1}/k}$ and $\sN^\vir$ inherit from $\phi$ and $E_\bullet$ natural $N$-linearizations.
\end{lemma}

\begin{proof} For $t\in T_1$, we have $t\cdot \sigma=\sigma\cdot t^{-1}$. If $\sF$ is an $N$-linearized coherent sheaf on $X$, $\tilde{i}^*\sF$ inherits an $N$-linearization on $X^{T_1}$, and this relation implies that the canonical isomorphism $\sigma\cdot\: \sigma^*(\sF)\to \sF$ given by the $N$-linearization restricts to $\sigma\cdot\:\sigma^*(\tilde{i}^*\sF_r)\to\tilde{i}^*\sF_{-r}$. Thus $\tilde{i}^*\sF_r\oplus 
\tilde{i}^*\sF_{-r}$ is an $N$-linearized subsheaf of $\tilde{i}^*\sF$ for each $r\neq0$, and $\tilde{i}^*\sF_0$ is also an $N$-linearized subsheaf of $\tilde{i}^*\sF$.   

Thus $\tilde{\phi}\:\tilde{i}^*E_\bullet^\fix\to \tau_{\le1}L_{X^{T_1}/k}$ and $\sN^\vir=\tilde{i}^*E_\bullet^\mov$ have canonical $N$-linearizations.
\end{proof}

We recall that  Witt sheaf cohomology is represented in $\SH(k)$ by the Eilenberg-MacLane spectrum $\EM(\sW_*)$, via
\[
\EM(\sW_*)^{a,b}(X)=H^{a-b}(X, \sW)
\]
for $X\in \Sm_k$. Here $\sW_*$ is the homotopy module which is $\sW$ in each degree; formally $\sW_*:=\sK^{MW}_*[\eta^{-1}]$, see \cite[\S5, \S6]{MorelICTP} for details. We use this formalism to be able to apply the machinery in $\SH(k)$ to Witt sheaf cohomology and its related constructions, such as the Borel-Moore homology and the equivariant theory. 
 
Given $X\in \Sch^N/k$, we write $H^\BM_{N,a}(X, \sW(\sL))$ for the corresponding $N$-equivariant Borel-Moore homology, as defined in \S\ref{sec:TEGEquivCoh}, via
\[
 H^\BM_{N,a}(X, \sW(\sL)):=\EM(\sW_*)^\BM_{N, a,0}(X/BN, \sL-\sO_X).
 \] 
 We similarly have the equivariant Witt sheaf cohomology $H^a_N(X, \sW(\sL))$, defined by
 \[
 H^a_N(W, \sW(\sL)):=\EM(\sW_*)^{a,0}_N(X, \sL-\sO_X).
 \]

We have shown in \cite[Theorem 8.6]{LevineAtiyahBott} that there is a positive integer $M_0$ such that,  for each $N$-linearized invertible sheaf $\sL$ on $X$,   the push-forward by the inclusion $i:X^{T_1}\to X$ induces an isomorphism
\[
i_*\:H^\BM_{N,a}(X^{T_1}, \sW(i^*\sL))[1/M_0e]\xrightarrow{\sim} H^\BM_{N,a}(X, \sW(\sL))[1/M_0e]
\]

\begin{definition}\label{def:Strict}
Let $|X|^N\subset X^{T_1}$ be the union of the $N$-stable irreducible components of $X^{T_1}$, and let  $X^{T_1}_\Ind$ 
be the union of the irreducible components $C$ of $X^{T_1}$ such that $N\cdot C=C\cup C'$ with $C'\neq C$ and $C\cap C'=\0$. We call the action {\em semi-strict} if $X^{T_1}=|X|^N\cup X^{T_1}_\Ind$ and {\em strict} if the action is semi-strict, $|X|^N\cap X^{T_1}_\Ind=\0$, and   $|X|^N$ decomposes as a disjoint union of two $N$-stable  closed  subschemes 
\[
|X|^N=X^N\amalg X^{T_1}_\fr 
\]
where the $N/T_1$-action on $X^{T_1}_\fr$ is free. 
\end{definition}

If the $N$-action is semi-strict, then the inclusion $|X|^N\to X$ induces an isomorphism 
\[
H^\BM_{N,*}(|X|^N, \sW(i^*\sL))[1/M_0e]\to H^\BM_{N,*}(X, \sW(\sL))[1/M_0e]
\]
and 
\[
H^\BM_{N,*}(X^{T_1}\setminus |X|^N, \sW(\sL))[1/M_0e]=0.
\]
This is \cite[Theorem 8.6]{LevineAtiyahBott}.

\begin{remark} Suppose that the $N$-action on $Y$ is strict. Then the $N$-action on $X$ is also strict, with $X^N=Y^N\cap X$, with $Y^{T_1}_\Ind\cap X\subset X^{T_1}_\Ind$, and with 
$Y^{T_1}_\fr\cap X=X^{T_1}_\fr\amalg X'$, with $X^{T_1}_\Ind=X'\amalg Y^{T_1}_\Ind\cap X$. 
\end{remark}

Assuming the $N$-action on $X$ is strict, we set  $|X|^N_j:=|X|^N\cap X_j$. This writes $X_j$ as a disjoint union
\[
X_j=|X|^N_j\amalg (X_j\cap X^{T_1}_\Ind)
\]
with $|X|^N_j$ and $N$-stable closed subscheme of $X^{T_1}$. This decomposes $X^{T_1}$ as a disjoint union
\[
X^{T_1}=(\amalg_j|X|^N_j)\amalg X^{T_1}_\Ind,
\]
and similarly decomposes $|X|^N$ as a disjoint union
\[
|X|^N=\amalg_j|X|^N_j.
\]
Decomposing each $|X|^N_j$ further, we have the connected components   $i_{j,\ell}\:|X|^N_{j,\ell}\to X$ of $|X|^N_j$, which are all $N$-stable.

Let $i_j\:|X|^N_j\to X$ be the inclusion,The $N$-linearized perfect obstruction theory $\tilde{\phi}$ on $X^{T_1}$ restricts to an $N$-linearized perfect obstruction theory $\phi_j\: (i_j^*E_\bullet)^\fix\to \tau_{\le1}L_{|X|^N_j/k}$ on $|X|^N_j$, and $N$-linearized the virtual conormal sheaf $\sN^\vir$ on $X^{T_1}$ restricts to the $N$-linearized  virtual conormal sheaf $\sN^\vir_j:=(i_j^*E_\bullet)^\mov$ on $|X|^N_j$. We may further restrict to $|X|^N_{j,\ell}$, giving the $N$-linearized perfect obstruction theory $\phi_{j,\ell}\: (i_{j,\ell}^*E_\bullet)^\fix\to \tau_{\le1}L_{|X|^N_{j,\ell}/k}$, and  the $N$-linearized  virtual conormal sheaf $\sN^\vir_{j,\ell}:=(i_{j,\ell}^*E_\bullet)^\mov$ on $|X|^N_{j,\ell}$. We let $N_j:=\V(\sN_j)$, $N_{j,\ell}:=\V(\sN_{j,\ell})$ denote the corresponding virtual vector bundles.

In the case of a strict $N$-action, each $N$-linearized locally free sheaf $\sV$ on a connected component of $X^{T_1}$, with $\sV=\sV^\mov$,  has a {\em generic representation type} $[\sV^\gen]$, which is an isomorphism class of an $N$-representation on a $k$-vector space (see \cite[Construction 2.7]{LevineAtiyahBott} for details). 

\begin{lemma}\label{lem:MovingEulerClass} Let $X$ be in $\Sch^N/k$. Let $\sV$ be an $N$-linearized locally free sheaf on $X^{T_1}$. Suppose  that the $N$-action is strict and that $\sV=\sV^\mov$ on $|X|^N$. Fix a $j, \ell$ and let $\sV_{j,\ell}$ denote the restriction of $\sV$ to $|X|^N_{j,\ell}$.   \\[5pt]
1. There are positive integers $M,n$ such that  the equivariant Euler class $e_N(\sV_{j,\ell}^{gen})\in H^*(BN, \sW)$ is of the form $M\cdot e^n$. If the weight space $(\sV_{j,\ell})_r$ for the $T_1$-action is zero for all even $r>0$, then $M$ is an odd integer. \\[2pt]
2. If  $k$ admits a real embedding, then  $e_N(\sV_{j,\ell}^{gen})$ is not nilpotent. If $M$ is odd, then  $e_N(\sV_{j,\ell}^{gen})$ is a non-zero divisor in $H^*(BN, \sW)$.\\[2pt]
3. $e_N(\sV_{j,\ell})$ is invertible in $H^*_N(|X|^N_{j,\ell}, \sW(\det^{-1}\sV))[1/Me]$.
\end{lemma}
\begin{proof}
This follows immediately from \cite[Lemma 9.3]{LevineAtiyahBott}.
\end{proof}

\begin{definition}\label{def:NormalEulerClass} Let $i\:X\to Y$ be a closed immersion in $\Sch^N/k$, with $Y\in \Sm^N/k$, and let $\phi\:E_\bullet\to  \tau_{\le1}L_{X/k}$ be an $N$-linearized perfect obstruction theory on $X$. Suppose that the action on $X$ is strict.  For each $i_{j,\ell}\:|X|^N_{j,\ell}\to X$ of $|X|^N$, write $e_N((i_{j,\ell}^*E_1^\mov)^{gen})\cdot e_N((i_{j,\ell}^*E_0^\mov)^{gen})$ as $M^X_{j,\ell}e^{n^X_{j,\ell}}$, following Lemma~\ref{lem:MovingEulerClass}. We define
\[
e_N(N^\vir_{j,\ell})\in H^*_N(|X|^N_{j,\ell}, \sW(\det^{-1}N_{j,\ell}))[1/M^X_{j,\ell}e]
\]
by
\[
e_N(N^\vir_{j,\ell})=e_N(\V((i_{j,\ell}^*E_0)^\mov))\cdot e_N(\V((i_{j,\ell}^*E_1)^\mov))^{-1}.
\]
Let $M^X_j=\prod_\ell M_{j,\ell}$ and let 
\begin{multline*}
e_N(N^\vir_j)=\{e_N(N^\vir_{i_{j,\ell}})\}_\ell\in 
H^*_N(|X|^N_j, \sW(\det^{-1} N_j))[1/M^X_j]\\:=\prod_\ell H^*_N(|X|^N_{j,\ell}, \sW(\det^{-1}N_{j,\ell}))[1/M^X_j].
\end{multline*}

Similarly,  noting that each irreducible component $Y_j$ of $Y^{T_1}$ is smooth, we can write $e_N((i^{Y*}_j T_Y^\mov)^{gen})$ as $M^Y_je^{n^Y_j}$.
\end{definition} 
It follows from Lemma~\ref{lem:MovingEulerClass} that $e_N(N^\vir_j)$ is well-defined and is invertible in $H^*_N(|X|^N_j, \sW(\det^{-1}N^\vir_j))[1/M^X_je]$. Similarly, the Euler class $e_N(i^{Y*}_j T_Y^\mov)$ is invertible in $H^*_N(Y_j, \sW(\det^{-1}i^{Y*}_j T_Y^\mov))[1/M^Y_je]$.

Recall the integer $M_0>0$ defined above, just before Definition~\ref{def:Strict}. 

\begin{theorem}[Virtual localization]\label{thm:VirLoc} Let $i\:X\to Y$ be a closed immersion in $\Sch^N/k$, with $Y\in \Sm^N/k$, and let $\phi\:E_\bullet\to  L_{X/k}$ be an $N$-linearized perfect obstruction theory on $X$.  Suppose the $N$-action on $X$ is strict.  Let $M=M_0\cdot \prod_{i,j}M^X_j \cdot M^Y_i$, where the $M_j^X, M_j^Y$ are as in 
Definition~\ref{def:NormalEulerClass}. Let $[|X|^N_j,\phi_j]_N^\vir\in \EM(\sW_*)^\BM_N(|X|^N_j,  i_j^*E^\fix_\bullet)$ be the equivariant virtual fundamental class for the $N$-linearized perfect obstruction theory $\phi_j\: i_j^*E^\fix_\bullet\to \tau_{\le1}L_{|X|^N_j/k}$ on $|X|^N_j$. Then
\[
[X,\phi]_N^\vir=\sum_{j=1}^s i_{j*}([|X|^N_j,\phi_j]_N^\vir\cap e_N(N^\vir_{i_j})^{-1})\in
\EM(\sW_*)^\BM_N(X,  E_\bullet)[1/Me]
\]
\end{theorem}

\begin{remark}\label{rem:Concluding}
1. In our notation $H^\BM_{N,a}(-,-)$ for the Witt-sheaf equivariant Borel-Moore homology, we have 
\[
\EM(\sW_*)^\BM_N(X/BN,  E_\bullet)=H^\BM_{N, r}(X, \sW(\sL)), 
\]
with $r$ the virtual rank of $E_\bullet$ and $\sL$ the virtual determinant, and similarly for 
$ \EM(\sW_*)^\BM_N(|X|^N_j/BN,  i_j^*E^\fix_\bullet)$.\\[2pt]
2. One is ultimately not just interested in the virtual fundamental class $[X,\phi]_N^\vir$, but rather its ``quadratic degree''. In the classical case, with $[X,\phi]^\vir\in \CH_r(X)$, $r$ is the virtual rank of the perfect obstruction theory, and if $X$ is proper over $k$ and $r=0$, one can get the numerical invariant $\deg_k[X,\phi]^\vir\in \CH_0(\Spec k)=\Z$. If one computes this by torus localization (and $r=0$), one will arrive at a torus equivariant class, whose degree will land in $\CH_0(BT)[1/P]=\Z[x_1,\ldots,x_n][1/P]$, where $n$ is the rank of the torus, and $P$ is some non-zero homogeneous element of positive degree. However, one knows that this is the image of the non-equivariant degree in $\Z$ that one wishes to compute, so the rational function $\deg_{BT}[X,\phi]_T^\vir$ will actually be the  integer $\deg_k[X,\phi]^\vir$. 

In the case of Witt sheaf Borel-Moore homology, and $N$-localization, this situation is similar, except that the virtual class $[X,\phi]^\vir_\sW$ will land in $H^\BM_r(X, \sW(\sM))$, where as before $r$ is the virtual rank of the perfect obstruction theory and now $\sM$ is the virtual determinant. In order to have a well-defined push-forward to the point, landing in $H^\BM_0(\Spec k,\sW)=H^0(\Spec k, \sW)=W(k)$, one needs $r=0$ and an isomorphism $\sM\cong \sL^{\otimes 2}$ for some invertible sheaf $\sL$ on $X$, namely, an   orientation for $E_\bullet$. If one wants to compute $\deg_k[X,\phi]^\vir_\sW\in W(k)$ by localization, one needs to have $N$ acting on all this data,  and then one will have the equivariant degree $\deg_{BN} [X,\phi]^\vir_{\sW,N}$ in $H^0(BN, \sW)$.

At this point, two problems arise. By \cite[Proposition 5.5]{LevineBG} or  \cite[Theorem 5.1]{LevineAtiyahBott}, 
\[
H^*(BN, \sW)=W(k)[e,x]/(x^2-1,x(e+1))
\]
where $x$ has degree 0 and $e$ is the degree 2 class mapping to the equivariant Euler class of the tautological rank 2 bundle on $BN$, arising from $N\subset \SL_2\subset \GL_2$. In particular
$H^0(BN, \sW)$ consists of two $W(k)$-free summands,
\[
H^0(BN, \sW)=W(k)\cdot 1\oplus W(k)\cdot x.
\]
We have the ``forget the $N$-action'' map $\pi_{BN}^*:H^0(BN, \sW)\to H^0(\Spec k, \sW)=W(k)$, and 
\[
\pi_{BN}^*(\deg_{BN} [X,\phi]^\vir_{\sW,N})=\deg_k[X,\phi]^\vir_\sW\in W(k)
\]

We also have the pullback by the structure map $p_{BN}^*:W(k)= H^0(\Spec k, \sW)$, giving a right splitting to $\pi_{BN}^*$, and defining the summand $W(k)\cdot 1$ as the image 
$p_{BN}^*(W(k))$. Moreover, $\pi_{BN}^*(x)=0$, so if $\deg_{BN} [X,\phi]^\vir_{\sW,N}=a\cdot 1+b\cdot x$, then $\deg_k[X,\phi]^\vir_\sW=a\in W(k)$.

If we use $N$-localization, we only obtain the image of $\deg_{BN} [X,\phi]^\vir_{\sW,N}$ in the degree zero part of $H^*(BN, \sW)[1/Me]$ for some positive integer $M$. We have
\[
H^*(BN, \sW)[1/Me]=W(k)[1/M][e^{\pm 1}] 
\]
and under the localization map $x$ gets sent to $-1$. Thus, if $\deg_{BN} [X,\phi]^\vir_{\sW,N}=a\cdot 1+b\cdot x$, then under localization, $\deg_{BN} [X,\phi]^\vir_{\sW,N}$ gets sent to $a-b\in W(k)[1/M]$, and if we do not know $b$, then the information given by $N$-localization tells us nothing about $\deg_k[X,\phi]^\vir_\sW=a\in W(k)$.

There is a solution to this problem, valid in all the geometric cases that we know of. See Proposition~\ref{prop:GroupExtension} below, which says that, if the $N$-action (together with the $N$-linearization of the perfect obstruction theory and the $N$-linearized orientation) comes by restriction from a group $G$ of the form $\prod_{i=1}^r \SL_{n_i}\times \prod_{j=1}^s \GL_{m_j}$ via some homomorphism $\rho:N\to G$, then indeed $\deg_{BN} [X,\phi]^\vir_{\sW,N}$ lands in the summand $W(k)\cdot 1$, in other words, $b=0$, and we recover at least the image of $\deg_k([X,\phi]^\vir_\sW)$ in $W(k)[1/M]$ by using $N$-localization.

This raises the second problem. In general   $W(k)$ will have a lot of 2-torsion. We will usually not have any information about $M$, which may be even, in which case we will lose all 2-primary torsion information once we invert $M$. What will in any case remain in $W(k)[1/M]$ are the $\Z$-valued invariants given by the signature map $W(\R)\xrightarrow{\sim} \Z$ for each real embedding $k\hookrightarrow \R$, or more generally, for each embedding of $k$ into a real closed field.  

There are some tricks one can use to improve the situation beyond this, for example, if everything is defined over $\Z[1/m]$ for some $m>0$, then one gets  a lifting of $\deg_k[X,\phi]^\vir_\sW$ to  $W(\Z[1/2m])$. As the kernel of $ W(\Z[1/2m])\to W(\R)$ is generated by the quadratic forms $x\mapsto \pm ax^2$ with $a$ a product of primes dividing $2m$, knowing just the signature recovers $\deg_k[X,\phi]^\vir_\sW$ modulo this kernel. See   \cite{KLSW, LevinePauli} for examples of this method. 
\end{remark}

The proof of Theorem~\ref{thm:VirLoc}  is accomplished by following the argument used by Graber-Pandharipande in  \cite[\S 3]{GP} for the proof of their main result. The argument there is essentially formal, relying on the various operations for the equivariant Chow groups, which as we have explained in \S\ref{sec:TEGEquivCoh}, extend to the setting of equivariant Borel-Moore homology and equivariant cohomology. The main non-formal ingredient in their proof is their use of of Vistoli's lemma, which we have extended here to the setting of  Borel-Moore homology, and which thus extends to equivariant Borel-Moore homology as defined in \S\ref{sec:TEGEquivCoh} via the algebraic Borel construction. We give here a step-by-step sketch of our extension of their proof. We set $\sE=\EM(\sW_*)$.
\\[10pt]
{\bf Step 1}  We note that the localized equivariant Borel-Moore homology of $X^{T_1}_\Ind$ is zero  \cite[Lemma 8.5(2)]{LevineAtiyahBott}.  Using the assumption that the $N$-action is strict, we have $X_j=|X|^N_j\amalg 
(X_j\cap X^{T_1}_\Ind)$, so using the localization sequence in equivariant Borel-Moore homology, we may replace $X$ with $X\setminus X^{T_1}_\Ind$, $Y$ with $Y\setminus X^{T_1}_\Ind$, and thus we may assume that $X^{T_1}_\Ind=\0$ and $|X|^N_j=X_j$.  We have the inclusions $\iota\:X\to Y$,  $\iota_j\:X_j\to Y_j$,  $i^Y_j\:Y_j\to Y$ and $i_j\:X_j\to X$.

We will use all the results of \S\ref{sec:Prelim} and \S\ref{sec:Vistoli} promoted to the equivariant setting, following Remark~\ref{rem:BMEquivHomOps}. We write $T_Y$ for $T_{Y/B}$, etc.

We have the equivariant fundamental class $[Y]_N\in \sE^\BM_N(Y/BN, \Omega_{Y/k})$, with $[Y]_N$ the image of $1\in \sE^{0,0}_N(Y)$ under the Poincar\'e duality isomorphism 
\[
\sE^{0,0}_N(Y)\cong \sE^\BM_N(Y/BN, \Omega_{Y/k}).
\] 

We recall that $Y^{T_1}$ is smooth over $k$ with irreducible components  $Y_1,\ldots, Y_s$. We have  fundamental classes $[Y_j]_N\in 
\sE^\BM_N(Y_j/BN, \Omega_{Y_j/B})$  defined as for $[Y]_N$.

 For each irreducible component $i^Y_j\:Y_j\to Y$ of $Y^{T_1}$, the normal bundle $N_{i^Y_j}$ of $i^Y_j$ is exactly $i^{Y*}_j T_Y^\mov$. Since $i_j^{Y!}[Y]_N=[Y_j]_N$,  applying the Bott residue theorem \cite[Theorem 10.5, Remark 10.6]{LevineAtiyahBott} for the  $N$-action on $Y$  gives
\begin{equation}\label{eqn:LocIdentity1}
[Y]_N=\sum_{j=1}^si^Y_{j*}\left([Y_j]_N\cap e_N(i_j^{Y*}T_Y^\mov)^{-1}\right).
\end{equation}
in $\sE^\BM_N(Y/BN, \Omega_{Y/k})[1/M_Ye]$. 

The Bott residue theorem  of \cite{LevineAtiyahBott} is written for the slightly different equivariant Borel-Moore homology $\sE^\BM_N(-/B,-)$. However, letting $\mathfrak{n}$ denote the Lie algebra of $N$, with its adjoint representation, we have 
\[
\sE^\BM_N(X/B,v-\mathfrak{n}^\vee)=\sE^\BM_N(X/BN,v),
\]
which allows us to apply the results of 
\cite{LevineAtiyahBott} to $\sE^\BM_N(-/BN,-)$ by simply replacing our  $\sE^\BM_N(X/BN,v)$ with $\sE^\BM_N(X/B,v-\mathfrak{n}^\vee)$ in all statements we need from \cite{LevineAtiyahBott}.
\\[5pt]
{\bf Step 2}. We apply the $N$-equivariant refined intersection product (\S\ref{subsec:RefinedProduct} and Remark~\ref{rem:EquivariantExtension}) over $Y$  with respect to the pair of maps $(\iota\:X\to Y,\id_Y)$, as well as the  pair $(\iota\:X\to Y,i^Y_j\:Y_j\to Y)$, and use the compatibility with proper push-forward (Lemma~\ref{lem:ProductComp}(1) and Remark~\ref{rem:EquivariantExtension})  to give the identity
\[
[X,\phi]^\vir_N=[X,\phi]^\vir_N\cdot_{\iota,\id}[Y]_N=\sum_{j=1}^si_{j*}\left([X,\phi]^\vir_N\cdot_{\iota,i^Y_j} ([Y_j]_N\cap e_N(i_j^{Y*}T_Y^\mov)^{-1})\right).
\]
in $\sE^\BM_N(X/BN, E_\bullet)[1/M_Ye]$.
 We will show that
\begin{equation}\label{eqn:MainIdentity}
[X,\phi]^\vir_N\cdot_{\iota,i^Y_j} ([Y_j]_N\cap e_N(i_j^{Y*}T_Y^\mov)^{-1})=
[X_j,\phi_j]^\vir_N\cap e_N(N^\vir_{i_j})^{-1}
\end{equation}
in $\sE^\BM_N(X_j/BN, E_\bullet)[1/Me]$ for suitable $M$,  which will yield the formula of Theorem~\ref{thm:VirLoc}. 
 \\[5pt]
{\bf Step 3}.  We have the cone $D=C_{X/Y}\times\V(E_0)$ over $X$  and its quotient cone $D^\vir:=\iota^*T_Y\backslash D$, with closed immersion $\iota_{E_\bullet}\:D^\vir\hookrightarrow \V(E_1)$, all of these being maps in $\Sch^N/k$.  

Referring to Definition/Proposition~\ref{defprop:EquivVFC}, we have 
\[
[X,\phi]^\vir_N=0^!_{\V(E_1)}(\iota_{E_{\bullet}*}[D^\vir]_N)\in \sE^\BM_N(X/BN, E_\bullet)
\]
 
For each $j$, we have a similar description of the virtual class $[X_j,\phi_j]_N^\vir$ for the $N$-linearized perfect obstruction theory $\phi_j\:i_j^*E^\fix_\bullet\to \tau_{\le1}L_{X_j/k}$  on $X_j$ with corresponding cones $D_j$ and $D_j^\vir$, namely,
\[
[X_j,\phi_j]_N^\vir=0^!_{\V(i_j^*E^\fix_1)}(\iota_{i_j^*E^\fix_{\bullet}*}[D_j^\vir]_N)\in \sE^\BM_N(X_j/BN, i_j^*E^\fix_{\bullet}).
\]
Here $D_j=C_{X_j/Y_j}\times \V(i_j^*E^\fix_0)$ and $D_j^\vir=\iota_j^*T_{Y_j}\backslash D_j$, with closed equivariant immersion $\iota_{i_j^*E^\fix_{\bullet}}\:D_j^\vir\hookrightarrow \V(i_j^*E^\fix_1)$.
\ \\[5pt]
{\bf Step 4}. For the rest of the proof, we use Remark~\ref{rem:EquivariantExtension}, without explicit mention, to extend results on $\sE$-cohomology and $\sE$-Borel-Moore homology to the equivariant setting, defined as in \S\ref{sec:TEGEquivCoh}.

We consider the closed immersions $\iota\:X \to  Y$ and $\iota_j\:X_j \hookrightarrow Y_j$, with their corresponding cones $C_\iota:=C_{X/Y}$ and $C_{\iota_j}:=C_{X_j/Y_j}$. We have the closed immersion $\beta_j\:C_{\iota_j}\to i_j^*C_\iota$. We now apply Vistoli's lemma.  Since $i^Y_j$ is a regular immersion, Proposition~\ref{prop:FundClassIdentity} and Remark~\ref{rem:VistoliLemChangeBase}  give us the identity
\begin{equation}\label{eqn:VId1}
\beta_{j*}([C_{\iota_j}]_N)=i_j^{Y!}[C_\iota]_N
\end{equation}
in $\sE^\BM_N(i_j^*C_\iota/BN)$;  we suppress the twist here and in the remainder of the proof to simplify the  notation.

The closed immersion $\beta_j$ and the inclusion $i_j^*E_0^\fix\subset i_j^*E_0^\fix\oplus i_j^*E_0^\mov=i_j^*E_0$ gives the  closed immersion $\beta^D_{j}\:D_j\times \V(i_j^*E^\mov_0)\to D$, and \eqref{eqn:VId1} gives us the relation
\begin{equation}\label{eqn:VId2}
\beta^D_{j*}([D_j\times \V(i_j^*E^\mov_0)]_N)=i_j^{Y!}[D]_N
\end{equation}
in $\sE^\BM_N(i_j^*D/BN)$. 
We also have the fundamental class $[D_j]$ induced from $[C_{\iota_j}]$ via the projection $D_j \to C_{\iota_j}$.

 We have the quotient cones 
\[
D^\vir:=\iota^*T_Y\backslash D,\ 
D_j^\vir:=i_j^*T_{Y_j}\backslash D_j,
\]
with their  fundamental classes $[D^\vir]$ and $[D_j^\vir]$ induced from $[D]$ and $[D_j]$. The closed immersion $\iota_{E_{\bullet}}\:D^\vir\hookrightarrow \V(E_1)$ induces the cartesian diagram
\[
\xymatrix{
\iota^*T_Y\ar[d]\ar[r]&D\ar[d]\\
X\ar[r]_-{0_{\V(E_1)}}&\V(E_1),
}
\]
which gives the identity
\[
[X, \phi]^\vir_N:=0_{\V(E_1)}^!(\iota_{E_{\bullet}*}[D^\vir]_N)=
0_{T_Y}^!0_{\V(E_1)}^!([D]_N).
\]
Combined with \eqref{eqn:VId2} and using the commutativity of refined Gysin pull-back,  this gives
\begin{multline*}
[X, \phi]^\vir_N\cdot_{\iota, i^Y_j}[Y_j]_N=
i_j^{Y!}[X, \phi]^\vir_N\\=
i_j^{Y!}0_{T_Y}^!0_{\V(E_1)}^!([D]_N)=
0_{i_j^{Y*}T_Y}^!0_{\V(i_j^{*}E_1)}^!i_j^{Y!}([D]_N)\\=
0_{i_j^{Y*}T_Y}^!0_{\V(i_j^{*}E_1)}^!([D_j\times  \V(i_j^*E^\mov_0)]_N).
\end{multline*}

Similarly, the closed immersion $\iota_{E_{\bullet, j}}\:D^\vir_j\hookrightarrow\V(i_j^*E_1^\fix)$ gives the cartesian diagram
\[
\xymatrix{
i_j^*T_{Y_j}\ar[d]\ar[r]&D_j\ar[d]\\
X_j\ar[r]_-{0_{\V(i_j^*E^\fix_1)}}&\V(i_j^*E^\fix_1)
}
\]
and  the identity
\[
[X_j, \phi_j]^\vir_N:=0_{\V(E_{1,j})}^!(\iota_{E_{\bullet, j}*}[D^\vir_j]_N)=
0_{i_j^*T_{Y_j}}^!0_{\V(E_{1,j})}^!([D_j]_N).
\]
We also note that 
\[
T_{Y_j}=i_j^{Y*}T^\fix_{Y}
\]
giving
\[
i_j^{Y*}T_Y=T_{Y_j}\oplus i_j^{Y*}T_Y^\mov,
\]
and in turn 
\begin{multline*}
0_{i_j^{Y*}T_Y}^!0_{\V(i_j^{*}E_1)}^![D_j \times  \V(i_j^*E^\mov_0)]_N\\
=0_{i_j^{Y*}T^\mov_Y}^!0_{\V(i_j^{*}E_1)}^![(\iota_j^*T_{Y_j}\backslash D_j)\times  \V(i_j^*E^\mov_0)]_N\\
=0_{i_j^{Y*}T^\mov_{Y/B}}^!0_{\V(i_j^{*}E_1)}^![D^\vir_j \times  \V(i_j^*E^\mov_0)]_N.
\end{multline*}

Following the discussion on \cite[pg. 498]{GP}, we have the commutative diagram 
\[
\xymatrix{
&i_j^*\iota^*T_Y\ar[d]\ar@{^(->}[r]& i_j^*D\ar[d]&\ar@{_(->}[l] D_j\times \V(i_j^*E_0^\mov)\ar[d]\\
i_j^*\iota^*T_Y^\mov\ar@{=}[r]&i_j^*\iota^*T_Y/\iota_j^*T_{Y_j}\ar[d]\ar@{^(->}[r]&i_j^*D/\iota_j^*T_{Y_j}
\ar[d]\ar@{}[dr]|-\bullet&
\ar@{_(->}[l] D^\vir_j\times_{X_j}\V(i_j^*E_0^\mov)\ar[d]^{(i_{i_j^*E_\bullet^\fix},\V(d_j^\mov))} \\
&X_j\ar@{^(->}[r]^{0_{\V(i_j^*E_1)}}&\V(i_j^*E_1)&\ar@{=}[l] \V(i_j^*E_1^\fix)\times_{X_j} \V(i_j^*E_1^\mov),
}
\]
with the ``hook'' arrows all closed immersions and all squares cartesian, except for the one marked ``$\bullet$". The map 
$\V(d_j^\mov)\:\V(i_j^*E_0^\mov)\to \V(i_j^*E_1^\mov)$ is induced from the ``moving part'' of the differential $d_j\:i_j^*E_1\to i_j^*E_0$.

Let $0_{\V(i_j^*E_1)}^{-1}(-)$ denote the scheme-theoretic pull-back  with respect to the zero-section $0_{\V(i_j^*E_1)}$.  From the above diagram, we have the closed immersion
\[
b\: 0_{\V(i_j^*E_1)}^{-1}(D^\vir_j\times_{X_j}\V(i_j^*E_0^\mov))\hookrightarrow
i_j^*\iota^*T_Y/\iota_j^*T_{Y_j}=i_j^*\iota^*T_Y^\mov
\]
Since $0_{\V(i_j^*E^\fix_1)}^{-1}(D^\vir_j)\subset X_j$, we also have the closed immersion
\[
a\:0_{\V(i_j^*E_1)}^{-1}(D^\vir_j\times_{X_j}\V(i_j^*E_0^\mov))\hookrightarrow
\V(i_j^*E_0^\mov).
\]
This gives us the commutative diagram in $\Sch^N/k$
\[
\xymatrix{
0_{\V(i_j^*E_1)}^{-1}(D^\vir_j \times_{X_j}\V(i_j^*E_0^\mov))\ar[r]^-a\ar[d]_b&\V(i_j^*E_0^\mov)\ar[d]\\
i_j^*\iota^*T_Y^\mov\ar[r]&X_j.
}
\]

Let $\gamma=0_{\V(i_j^*E_1)}^![D^\vir_j \times  \V(i_j^*E^\mov_0)]_N$, this class living in the equivariant Borel-Moore homology
$\sE^\BM_N(0_{\V(i_j^*E_1)}^{-1}(D^\vir_j \times_{X_j}\V(i_j^*E_0^\mov))/BN)$. Using the excess intersection formula as in \cite[\S3, Lemma 1]{GP}, we have
\[
0^!_{i_j^{Y*}T^\mov_Y}(a_*(\gamma))\cap e_N(\V(i_j^*E^\mov_0)) =0^!_{\V(i_j^*E^\mov_0)}(b_*(\gamma))\cap  e_N(i_j^{Y*}T^\mov_Y).
\]
This gives
\begin{multline*}
[X,\phi]_N^\vir\cdot_{\iota, i^Y_j}[Y_j]_N\cap e(\V(i_j^*E^\mov_0))\\=
0_{\V(i_j^*E^\mov_0)}^!(0_{\V(i_j^{*}E_1)}^![D^\vir_j \times  \V(i_j^*E^\mov_0)]_N)\cap e_N(i_j^{Y*}T^\mov_Y).
\end{multline*}

By using $\A^1$-homotopy invariance, it follows that the right-hand side of this identity is not changed if we replace the differential $d_j\:i_j^*E_1\to i_j^*E_0$ with the 0-map, so that the map 
\[
(\iota_{i_j^*E_\bullet^\fix},\V(d_j^\mov))\:D^\vir_j\times  \V(i_j^*E^\mov_0)\to \V(i_j^{*}E_1)
\]
factors as the projection $D^\vir_j\times  \V(i_j^*E^\mov_0)\to D^\vir_j$ followed by the canonical  closed immersion  $D^\vir_j\hookrightarrow\V(i_j^{*}E^\fix_1)\subset \V(i_j^{*}E_1)$. Using the equivariant excess intersection formula, this yields
\begin{multline*}
0_{\V(i_j^*E^\mov_0)}^!(0_{\V(i_j^{*}E_1)}^![D^\vir_j\times  \V(i_j^*E^\mov_0)]_N)\\=
(0_{\V(i_j^{*}E_1)}^![D^\vir_j \times  \V(i_j^*E^\mov_0)])\cap e_N(\V(i_j^*E^\mov_0))\\
=[X_j,\phi_j]_N^\vir\cap  e_N(\V(i_j^*E^\mov_1)),
\end{multline*}
or
\begin{multline*}
[X, \phi]_N^\vir\cdot_{\iota, i^Y_j}([Y_j]_N\cap e_N(i_j^{Y*}T^\mov_Y)^{-1})\\=
[X_j,\phi_j]_N^\vir\cap (e_N(\V(i_j^*E^\mov_1))\cdot e_N(\V(i_j^*E^\mov_0))^{-1}).
\end{multline*}
This is exactly the desired formula \eqref{eqn:MainIdentity}.

We conclude with the promised Proposition~\ref{prop:GroupExtension}. We use the change of group map described in \S\ref{subsec:Change}.

\begin{proposition}\label{prop:GroupExtension} Let $G =\prod_{i=1}^r \SL_{n_i}\times \prod_{j=1}^s \GL_{m_j}$ as group scheme over $k$. Take $X\in \Sch^G/k$, endowed with a $G$-linearized perfect obstruction theory $\phi_G:E_\bullet\to L_{X/k}$. Suppose that $X$ is proper over $k$, that $\rnk E_\bullet=0$ and that  we have a $G$-linearized invertible sheaf $\sL$ on $X$ with $G$-equivariant isomorphism $\tau_G:\det E_\bullet\xrightarrow{\sim} \sL^{\otimes 2}$.

Finally, suppose we have a homomorphism of $k$-group schemes $\rho:N\to G$,   let $\phi_N$ be the $N$-linearized perfect obstruction theory induced from $\phi_G$. We use the orientation $\tau_G$, restricted to an $N$-equivariant orientation $\tau_N$, to define the degree
\[
\deg_{BN}([X,\phi_N]^\vir_\sW)\in H^0(BN, \sW)=W(k)\cdot 1\oplus W(k)\cdot x,
\]
where the decomposition of $H^0(BN, \sW)$ is as described in Remark~\ref{rem:Concluding}. \\[5pt]
1. We have that $\deg_{BN}([X,\phi_N]^\vir_\sW)$ is in the summand $W(k)\cdot 1$ of $H^0(BN, \sW)$.
\\[2pt]
2. Let $\phi$ be the perfect obstruction theory on $X$ gotten from $\phi_N$ by forgetting the $N$-action and endow $E_\bullet$ with the orientation $\tau$ induced from $\tau_N$, again by forgetting the $N$-action. This gives us virtual fundamental class $[X,\phi]^\vir_\sW\in H^\BM_0(X)$ and its degree
\[
\deg_k([X,\phi]^\vir_\sW)\in W(k).
\]
Then the image of $\deg_{BN}([X,\phi_N]^\vir_\sW)$ in the degree 0 part $W(k)[1/M]$ of  $H^*(BN, \sW)[1/Me]$  is the same as the image of $\deg_k([X,\phi]^\vir_\sW)\in W(k)$ in the localization $W(k)[1/M]$ of $W(k)$.
\end{proposition}

\begin{proof} (2) follows from (1), following the discussion in Remark~\ref{rem:Concluding}. For (1), it follows from \cite[Theorem 10]{Anan} and \cite[Theorem 4.1]{LevineBG} that $H^*(BG,\sW)$ is a polynomial ring over $W(k)$ in finitely many variables. In particular, $H^0(BG,\sW)=W(k)$ and both the pullback map $p_{BG}^*:W(k)\to H^0(BG,\sW)$, and the forget map $\pi_{BG}^*:H^0(BG,\sW)\to W(k)$ are isomorphisms. We now apply the restriction map $\rho^\#$ of \S\ref{subsec:Change}. This gives us the commutative diagram
\[
\xymatrix{
W(k)=H^0(\Spec k, \sW)\ar@{=}[r]\ar[d]^{p_{BG}^*}_\wr&H^0(\Spec k, \sW)\ar[d]^{p_N^*}\\
H^0(BG, \sW)\ar[r]^{\rho^\#}\ar[d]^{\pi_{BG}^*}_\wr&H^0(BN, \sW)\ar[d]^{\pi_{BN}^*}\\
W(k)=H^0(\Spec k, \sW)\ar@{=}[r]&H^0(\Spec k, \sW)
}
\]
Thus $\rho^\#(H^0(BG, \sW)$ is the summand $W(k)\cdot 1$ of $H^0(BN, \sW)$. 

On the other hand, by Proposition~\ref{prop:VirtChangeOfGroup}, we have
\[
\rho^\#(\deg_{BG}([X,\phi_G]^\vir_\sW))=\deg_{BN}([X,\phi_N]^\vir_\sW)
\]
so $\deg_{BN}([X,\phi_N]^\vir_\sW)$ is in the summand $W(k)\cdot 1$ of $H^0(BN, \sW)$, as claimed.
\end{proof}

Using Proposition~\ref{prop:EulerClassChangeOfGroup} instead of  Proposition~\ref{prop:VirtChangeOfGroup}, we have a parallel result for the use of $N$-localization to compute the degree of an Euler class; the proof is the same.

\begin{proposition}\label{prop:GroupExtensionEulerClass} Let $G =\prod_{i=1}^r \SL_{n_i}\times \prod_{j=1}^s \GL_{m_j}$ as group scheme over $k$. Take $X\in \Sch^G/k$, endowed with a $G$-equivariant vector bundle $V\to X$, with sheaf of sections $\sV$.  Suppose that $X$ is proper and lci of pure dimension $d$ over $k$, that $\rnk\, V=d$ and that  we have a $G$-linearized invertible sheaf $\sL$ on $X$ with $G$-equivariant relative orientation $\tau_G:\det L_{X/B}\otimes \det \sV^\vee\xrightarrow{\sim}  \sL^{\otimes 2}$.

Suppose we have a homomorphism of $k$-group schemes $\rho:N\to G$, making $V$ an $N$-equivariant vector bundle on $X$, with an $N$-equivariant relative orientation $\tau_N:\det L_{X/B}\otimes \det \sV^\vee\xrightarrow{\sim}  \sL^{\otimes 2}$.  

 This gives us the equivariant Borel-Moore Euler class 
 \[
 e^\BM_{N}(V)\in H^\BM_{N,0}(X/BN, \sW(\det L_{X/B}\otimes\det(\sV)^\vee)). 
 \]
 Using the relative orientation $\tau_N$ to define an isomorphism 
\[
\tau_{N*}:H^\BM_{N,0}(X/BN, \sW(\det L_{X/B}\otimes\det\sV^\vee))\xrightarrow{\sim} H^\BM_{N,0}(X/BN, \sW),
\]
we have the Borel-Moore Euler class $e^\BM_{N}(V)\in H^\BM_{N,0}(X/BN, \sW)$, and the degree
\[
\deg_{BN}(e^\BM_{N}(V))\in H^0(BN, \sW)=W(k)\cdot 1\oplus W(k)\cdot x,
\]
where the decomposition of $H^0(BN, \sW)$ is as described in Remark~\ref{rem:Concluding}. \\[5pt]
1. We have that $\deg_{BN}(e^\BM_{N}(V))$ is in the summand $W(k)\cdot 1$ of $H^0(BN, \sW)$.
\\[2pt]
2. Let 
$e^\BM(V)\in H^\BM_0(X,\sW(\det L_{X/B}\otimes \det\sV^\vee))$ be Borel-Moore Euler class  of $V$,  We use the given relative orientation $\tau$ induced from $\tau_N$ by forgetting the $N$-action to give the isomorphism 
$H^\BM_0(X,\sW(\det L_{X/B}\otimes \det\sV^\vee))\cong H^\BM_0(X,\sW)$. This gives us the Borel-Moore Euler class
\[
e^\BM(V)\in H^\BM_0(X,\sW),
\]
 and its degree
\[
\deg_k(e^\BM(V))\in W(k).
\]
Then the image of $\deg_{BN}(e^\BM_{N}(V))$ in the degree 0 part $W(k)[1/M]$ of  $H^*(BN, \sW)[1/Me]$  is the same as the image of $\deg_k(e^\BM(V))\in W(k)$ in the localization $W(k)[1/M]$ of $W(k)$.
\end{proposition}

\end{document}